\documentclass[dvips,12pt]{article}
\usepackage{stmaryrd,amssymb,amsmath,amsfonts,pstricks,pstricks-add
}
\def\bt{\begin{tabular}}
\def\te{\end{tabular}}

\def\BM{\begin{pmatrix}}
\def\EM{\end{pmatrix}}

\def\cit{\text{\it I\hskip -6ptC\/}}

\def\square{\hfill\hbox{\vrule height .9ex width .8ex depth -.1ex}}
\def\rit{\text{\it I\hskip -2pt  R}}

\def\zit{\text{\it Z\hskip -4pt  Z}}

\def\nit{\text{\it I\hskip -2pt  N}}

\def\Bd{{\text B}}

\def\Ed{{\text E}}

\def\Fs{{\cal F}}

\def\Hs{{\cal H}}

\def\Ls{{\cal L}}
\def\Ms{{\cal M}}

\def\Os{{\cal O}}

\def\Ps{{\cal P}}

\def\be{\begin{equation}}
\def\ee{\end{equation}}
\def\beqn{\begin{eqnarray}}
\def\eeqn{\end{eqnarray}}
\def\nobeqn{\begin{eqnarray*}}
\def\noeeqn{\end{eqnarray*}}
\def\ba{\left(\begin{array}}
\def\ea{\end{array} \right) }

\def\bpr{\paragraph{Proof.}}

\def\epr{\square\vskip 6pt}
\def\eop{\hbox{\vrule height .9ex width .8ex depth -.1ex}}

\def\o{\overline}
\def\and{\; \mbox{and} \;}

\newcommand{\half}{\frac{1}{2}}

\def\Ker{\mathop{\rm Ker}\nolimits}

\def\Be{\begin{enumerate}}
\def\Ee{\end{enumerate}}

\def\Bena{\begin{enumerate}
\def\labelenumi{\theenumi)}
\def\theenumi{\arabic{enumi}}
\def\labelenumii{\theenumii)}
\def\theenumii{\alph{enumii}}}

\def\Bean{\begin{enumerate}
\def\labelenumii{\theenumii)}
\def\theenumii{\arabic{enumii}}
\def\labelenumi{\theenumi)}
\def\theenumi{\alph{enumi}}}

\def\Bero{\begin{enumerate}
\def\labelenumii{\theenumii)}
\def\theenumii{\arabic{enumii}}
\def\labelenumi{(\theenumi)}
\def\theenumi{\roman{enumi}}}

\def\BeRo{\begin{enumerate}
\def\labelenumii{\theenumii)}
\def\theenumii{\arabic{enumii}}
\def\labelenumi{(\theenumi)}
\def\theenumi{\Roman{enumi}}}

\def\Bi{\vskip 11pt\begin{itemize}\itemsep=18pt}%\topsep=18pt\itemsep=18pt\partopsep=18pt\parsep=18pt}
\def\Ei{\end{itemize}}

\def\Bd{\begin{description}}
\def\Ed{\end{description}}

\def\R{\right}
\def\L{\left}
\def\F{\frac}

\def\bigoplus{\mathop{\oplus}\limits}

\def\prod{\mathop{\Pi}\limits}
\def\sum{\mathop{\Sigma}\limits}

\def\Bsm{\begin{smallmatrix}}
\def\Esm{\end{smallmatrix}}

\def\resp#1{(resp. #1)}
\def\rresp#1{\qquad \mbox{(resp.} \quad #1\ )}

\def\bbf{\boldmath\bf}
\def\o{\overline}
\def\wt{\widetilde}
\def\wh{\widehat}

\def\Bi{\begin{itemize}}%\topsep=18pt\itemsep=18pt\partopsep=18pt\parsep=18pt}
\def\Ei{\end{itemize}}

\newcommand{\MM}{\mathbb{M}\,}

\def\ung{\hbox{\rm 1\hskip -3pt \rm I}\,}

\def\tr{\operatorname{tr}}
\def\det{\operatorname{det}}
\def\Gal{\operatorname{Gal}}

\def\End{\operatorname{End}}

\def\AdF{\operatorname{Ad(F)}}

\def\Irr{\operatorname{Irr}}

\def\Rep{\operatorname{Rep}}

\def\Repsp{\operatorname{Repsp}}%%%%%%%%%%%%%
\def\REPSP{\operatorname{REPSP}}%%%%%%%%%%%%%
\def\FREPSP{\operatorname{FREPSP}}%%%%%%%%%%%%%
\def\fREPSP{\operatorname{(F)REPSP}}%%%%%%%%%%%%%
\def\TAN{\operatorname{TAN}}%%%%%%%%%%%%%

\def\Int{\operatorname{Int}}
\def\Aut{\operatorname{Aut}}
\def\Out{\operatorname{Out}}

\def\ELLIP{\operatorname{ELLIP}}
\def\ellip{\operatorname{ellip}}

\def\BGL{\operatorname{BGL}}
\def\BFGL{\operatorname{BFGL}}

\def\EGL{\operatorname{EGL}}
\def\GL{\operatorname{GL}}
\def\gl{\operatorname{gl}}
\def\FGL{\operatorname{FGL}}
\def\FG{\operatorname{FG}}

\def\GD{\operatorname{GD}}

\def\Tr{\operatorname{tr}}

\def\bM{\begin{matrix}}
\def\eM{\end{matrix}}

\def\lr{left (resp. right) }

\def\bpr{\noindent{\bf{Proof\/}}:\;\;}
\def\To{\longrightarrow }

\def\RL{_{R\times L}}

\def\Rl{_{R,L}}

\def\Fvt{F_{\o v}\times F_v}
\def\FTvt{F^T_{\o v}\times F^T_v}
\def\Fnrvt{F^{nr}_{\o v}\times F^{nr}_v}
\def\Fvtrit{F_{\o v}\times \rit)\times(F_{v}\times \rit}
\def\FvTtrit{F^T_{\o v}\times \rit)\times(F^T_{v}\times \rit}
\def\tFvt{\wt F_{\o v}\times \wt F_v}

\begin{document}

\setcounter{page}{0}
{\pagestyle{empty}
\null\vfill
\begin{center}
{\LARGE Random matrices and Riemann hypothesis}
\vfill
{\sc C. Pierre\/}
\vskip 11pt

%Institut de Mathématique pure et appliquée\\
%Université de Louvain\\
%Chemin du Cyclotron, 2\\
%B-1348 Louvain-la-Neuve,  Belgium\\
%pierre@math.ucl.ac.be

\vfill
%Mathematics subject classification (1991): 11G18, 11R34, 11R37, 11R39.
\vfill

%{\LARGE $n$-dimensional global correspondences of Langlands}
%\vfill
%{\sc C. Pierre\/}
%\vfill
\begin{abstract}
The curious connection between the spacings of the eigenvalues of random matrices and the corresponding spacings of the  nontrivial zeros of the Riemann zeta function is analyzed on the basis of the geometric dynamical global program of Langlands whose fundamental structures are shifted quantized conjugacy class representatives of bilinear algebraic semigroups.  The considered symmetry behind this phenomenology is the differential bilinear Galois semigroup shifting the product, right by left, of automorphism semigroups of cofunctions and functions on compact transcendental quanta.
\end{abstract}
\vfill

\end{center}

11S37 --- 14G10 --- 34L15 --- 15A52.
\eject

\tableofcontents
\vfill\eject}

\setcounter{page}{1}
\def\thepage{\arabic{page}}
{\parindent=0pt 
\setcounter{section}{0}
%\addtocontents{toc}{\protect\vspace* {-6pt}
%\addtocontents{toc}{\protect\vspace* {-6pt}}
\section*{Introduction}
\addcontentsline{toc}{section}{Introduction}

It became more and more evident over the years that {\bbf there must exist a connection between the eigenvalues of random matrices and the nontrivial zeros of the Riemann zeta function $\zeta (s)$\/} in such a way that the distribution of the nontrivial zeros of  $\zeta (s)$ can be approached by techniques developed in random matrix theory.

And, as the distributions over ensembles of random matrices may be related to the statistical properties of waves in chaotic dynamical systems \cite{Yau}, {\bbf the Riemann dynamics should be chaotic\/}.

Studying numerically the nontrivial zeros of $\zeta (s)$, Odlyzko \cite{Odl} found that the mean spacings between their imaginary part $\gamma _j$ and the spacings between the eigenvalues of large unitary matrices correspond.  This work was based on researches of Montgomery \cite{Mon} who found that the pair correlation of the $\gamma _j$ is equal to the Gaussian unitary ensemble (GUE) pair correlation.

In their famous paper \cite{K-S}, M. Katz and P. Sarnak put the big questions: ``{\bbf Why is this so and what does it tell us about the nature (e.g. spectral) of the zeroes? Also, what is the symmetry behind this ``GUE'' law?\/}''

It is the purpose of this paper to propose a satisfying solution to these questions based {\bbf on the following considerations:\/}
\Bena
\item {\bbf The fundamental structures behind this problem are those of the geome\-tric-dynamical global program of Langlands, i.e. double symmetric towers of shifted conjugacy class representatives of bilinear algebraic semigroups over sets of increasing finite algebraic extensions of number fields\/}.

\item {\bbf The quantization of these conjugacy class representatives, which are abstract subsemivarieties \cite{Har}, by algebraic (resp. transcendental) quanta being equivalently: a) irreducible closed algebraic (resp. transcendental) subsets characterized by an extension degree $N$\/} (which is the Artin conductor); {\bbf b) unitary irreducible representations of these semigroups\/}.

\item {\bbf The considered symmetry behind this phenomenology is the one of the differential bilinear Galois semigroup shifting the product, right by left, of automorphism semigroups of cofunctions and functions on compact transcendental quanta\/}.

\item {\bbf The Riemann dynamics then corresponds to the wave chaotic dynamics described by random matrix theory\/}.
\Ee
\vskip 11pt

{\bbf The general proposed solution is organized in three parts\/} (or chapters):
\Bena
\item {\bbf Reminder of universal algebraic structures of the global program of Langlands \cite{Pie2}, \cite{Lan}\/}.

\item {\bbf Reminder of universal dynamical structures of the global program of Langlands \cite{Pie3}\/}.

\item {\bbf Responses to five questions allowing finding a solution to this problem\/}:
{\bbf \Be
\item   What is behind random matrices leading to the Gaussian orthogonal ensemble (GOE) as well as the Gaussian unitary ensemble (GUE)?

\item What is behind the partition and correlation functions between eigenvalues of random matrices?

\item What interpretation can we give to the local spacings between the eigenvalues of large random matrices?

\item What interpretation can we give to the spacings between the nontrivial zeros of the Riemann zeta function $\zeta (s)$?

\item What is the curious connection between c) and d)?
\Ee}
\Ee
\vskip 11pt

Let us analyse more concretely these three parts.
\vskip 11pt

{\bbf In chapter 1, the bilinear global version of the Langlands program\/} is recalled to be {\bbf based on:\/}

\Bean
\item {\bbf bisemiobjects $(o_R\times o_L)$ composed of the products of right semiobjetcs $o_R$\/}, localized on the lower half space or on $\rit_-$, {\bbf and of left symmetric semiobjects $o_L$\/}, localized on the upper half space or on $\rit_+$.

\item {\bbf the product, right by left, of two symmetric towers of increasing algebraic or compact transcendental extensions composed of an increasing number of algebraic or transcendental compact quanta and defining a (bisemi)lattice of biquanta\/}, i.e. the product, right by left, of two symmetric (semi)lattices of algebraic or transcendental quanta.

\item {\bbf Abstract bisemivarieties $G^{2n}(\Fvt)$ over products, right by left, $(\Fvt)$ of sets of increasing compact transcendental extensions\/} (or archimedean completions)
$F_v=\{ F_{v_1},\dots,F_{v_j},\dots,F_{v_r}\}$ and
$F_{\o v}=\{F_{\o v_1},\dots,F_{\o v_j},\dots,F_{\o v_r}\}$
in such a way that they are (functional) representation spaces of the {\bbf algebraic bilinear semigroups 
$\GL_{2n}(\tFvt)=T_{2n}^t(\wt F_{\o v})\times T_{2n}(\wt F_v)$ being the $2n$-dimensional bilinear representations\/} of the products, right by left,
$\Gal(\wt F_{\o v}/k)\times \Gal(\wt F_{v}/k)$ {\bbf of global Weil semigroups\/}.

\item {\bbf The bilinear abstract parabolic semigroups 
$P^{(2n)}(F_{\o v^1}\times F_{v^1})$\/}, over products
$(F_{\o v^1}\times F_{v^1})$ of unitary transcendental pseudoramified extensions, in such a way that they {\bbf are unitary representation spaces of the algebraic bilinear semigroup of matrices 
$\GL_{2n}(\tFvt)$\/}.
\Ee

The {\bbf abstract bisemivarieties $G^{(2n)}(\Fvt)$\/}, covered by their algebraic equivalents 
 $G^{(2n)}(\tFvt)$, {\bbf are separated bisemischemes or quasiprojectives bisemivarieties \cite{Har}\/} and constitute the corner stone of the global program of Langlands since, by an isomorphism of toroidal compactification,
 $G^{(2n)}(\tFvt)$ (and 
 $G^{(2n)}(\Fvt)$) is transformed into {\bbf the cuspidal representation
 $\Pi (\GL_{2n}(\tFvt))$ of
 $\GL_{2n}(\tFvt)$\/}.
 
 Remark that, as abstract (bisemi)varieties are quantized, their associated (bisemi)schemes, referring to quasi-projective (bisemi)varieties, are also quantized.
 
 {\bbf The differentiable functional representation space
 $\FREPSP (\GL_{2n}(\Fvt))$ of the complete bilinear semigroup
 $\GL_{2n}(\Fvt)$ is a bisemisheaf
 $(\wh M^{(2n)}_{v_R}\otimes\wh M^{(2n)}_{v_L})$ of differential bifunctions on the abstract bisemivariety
 $G^{(2n)}(\Fvt)$\/} which is given by the product, right by left, of a right semisheaf
 $\wh M^{(2n)}_{v_R}$ of differentiable cofunctions on the increasing conjugacy class representatives of the abstract right semivariety $G^{(2n)}(F_{\o v})\equiv T^{(2n)}(F_{\o v})$ by a left semisheaf
 $\wh M^{(2n)}_{v_L}$ of symmetric differentiable functions on the increasing conjugacy class representatives of the abstract left semivariety $G^{(2n)}(F_{v})\equiv T^{(2n)}(F_{v})$.
 \vskip 11pt
 
 {\bbf Chapter 2 deals with dynamical (geometric) $\GL(2n)$-bisemistructures of the global program of Langlands generated from the action of the elliptic bioperator
$(D^{(2k)}_R\otimes D^{(2k)}_L)$\/} (i.e. the product of a right linear differential elliptic operator $D^{(2k)}_R$ acting on $2k$ variables by its left equivalent  $D^{(2k)}_L$, $k\le n$) {\bbf on the bisemisheaf 
 $(\wh M^{(2n)}_{v_R}\otimes\wh M^{(2n)}_{v_L})$\/} according to:
 \[
 D^{(2k)}_R\otimes
 D^{(2k)}_L: \quad \FREPSP (\GL_{2n}( \Fvt ) )\To
 \FREPSP (\GL_{2n[2k]}(\Fvtrit))\]
where $ \FREPSP (\GL_{2n[2k]}(\Fvtrit))$ is the functional representation space of
 $ \GL_{2n}(\Fvt)$ fibered or shifted in $2k$ bilinear geometric dimensions with
 $F^{S_\rit}_v=(F_v\times \rit)$ the set of increasing left compact transcendental extensions fibered or shifted by real numbers.
 
 {\bbf $ \FREPSP (\GL_{2k}(\Fvtrit))$ is isomorphic to the total bisemispace of the tangent bibundle
 $\TAN (\wh M^{(2k)}_{v_R}\otimes\wh M^{(2k)}_{v_L})$ to the bisemisheaf
 $ (\wh M^{(2k)}_{v_R}\otimes\wh M^{(2k)}_{v_L}) \subset
 (\wh M^{(2n)}_{v_R}\otimes\wh M^{(2n)}_{v_L}) $ of which bilinear fibre
 $\Fs^{(2k)}\RL(\TAN)=\fREPSP( \GL_{2k}(\rit\times\rit))$ is   the  (functional) representation space
 of 
 $ \GL_{2k}(\rit\times\rit)$ corresponding to the action of the bioperator
 $(D^{(2k)}_R)\otimes
 (D^{(2k)}_L)$\/}.
 
 It is then proved that {\bbf the bilinear semigroup of matrices
$\GL_r(\rit^{(2k)}\times\rit^{(2k)})$ of ``algebraic'' order $r$\/}, associated with the bilinear fibre 
$\GL_{2k}(\rit\times\rit)$ and referring to the action of the bioperator
$(D^{(2k)}_R\otimes D^{(2k)}_L)$ on the bisemisheaf
$ (\wh M^{(2k)}_{v_R}\otimes\wh M^{(2k)}_{v_L})$, {\bbf corresponds to the $2k$-dimensional ``geometric'' bilinear real representation of the product, right by left, of ``differential'' Galois (or global Weil) semigroups\/} 
$\Aut_k(\phi _R(\rit))\times\Aut_k(\phi _L(\rit))$, shifting or fibering the product, right by left, of automorphism semigroups
$\Aut_k(\phi _R(F_{\o v}))\times\Aut_k(\phi _L(F_v))$ of cofunctions 
$\phi _R(F_{\o v})$ and functions $\phi _L(F_v)$ respectively on the compact transcendental real extensions $F_{\o v}$ and $F_v$ by
$\Aut_k(\phi _R(F_{\o v}\times\rit))\times\Aut_k(\phi _L(F_v\times\rit))$
{\bbf such that
\begin{multline*} \qquad \GL_r (\phi _R(\rit^{(2k)})\times\phi _L(\rit^{(2k)}))\\
=\Rep^{(2k)} (\Aut_k(\phi _R(\rit))\times\Aut_k(\phi _L(\rit))
\;, \quad 1\le j\le r\;.\end{multline*}
}

Similarly, {\bbf the unitary parabolic bilinear semigroup
$P_r(\rit^{(2k)}\times\rit^{(2k)})\subset\linebreak
\GL_r(\rit^{(2k)}\times\rit^{(2k)})$, referring to the action of the bioperator
$(D^{(2k)}_R\otimes D^{(2k)}_L)$ on the unitary bisemisheaf
$ (\wh M^{(2k)}_{\o v_R^1}\otimes\wh M^{(2k)}_{\o v^1_L})
\subset(\wh M^{(2k)}_{v_R}\otimes\wh M^{(2k)}_{v_L})$, 
corresponds to the $2k$-dimensional bilinear real representation of the product, right by left, of ``differential'' inertia Galois (or global Weil) semigroups\/} 
$\Int_k(\phi _R(\rit))\times\Int_k(\phi _L(\rit))$
{\bbf such that\/}
\[ P_r (\phi _R(\rit^{(2k)})\times\phi _L(\rit^{(2k)}))
=\Rep^{(2k)} (\Int_k(\phi _R(\rit))\times\Int_k(\phi _L(\rit))
\;.\]
Let $O_r(\rit)$ denote the orthogonal group of (algebraic) order $r$ with entries in $\rit$ and let $U_r(\cit)$ denote the unitary group with entries in $\cit$.  Then,  the orthogonal bilinear semigroup $O_r(\rit^{(2k)}\times\rit^{(2k)})$ corresponds to the real parabolic bilinear semigroup $P_r(\rit^{(2k)}\times \rit^{(2k)})$ and the unitary bilinear semigroup
$U_r(\cit^k\times\cit^k)$ corresponds to the complex parabolic bilinear semigroup
$P_r(\cit^k\times\cit^k)$.
\vskip 11pt

{\bbf Chapter 3 deals with the connection between large random matrices and the solution of the Riemann hypothesis\/} by responding to the five above mentioned questions.  From now on, {\bbf the geometric dimension $2k$ will be taken to be $1$\/}.
\Be
\item {\bbf The first question ``What is behind random matrices leading to GOE and GUE?'' leads to the eigenvalue problem in the frame of the geometric dynamical program of Langlands\/} recalled in chapter 2.

{\bbf The symmetric group at the origin of the bilinear global program of Langlands is the bilinear semigroup of automorphisms\/}
$\Aut_k(F_{\o v})\times \Aut_k(F_{v})$
\resp{Galois automorphisms $\Gal(\wt F_{\o v}/k)\times\Gal(\wt F_{v}/k)$}
{\bbf of compact transcendental (resp. algebraic) quanta generating a bisemilattice\/} of compact transcendental (resp. algebraic) quanta
{\bbf while the symmetry group at the origin of the dynamical (geometric) bilinear global program of Langlands is the bilinear semigroup of shifted automorphisms\/}
$\Aut_k(\phi _R((F_{\o v}\times\rit))\times \Aut_k(\phi _L((F_{v}\times\rit))$
{\bbf of bifunctions on compact transcendental biquanta\/} generating a bisemilattice of compact  transcendental quanta.

It then results that {\bbf the bilinear semigroup of matrices
$\GL_r(\rit\times\rit)$ constitutes the ``$r$-dimensional algebraic'' representation of the bilinear differential Galois semigroup associated with the action of the differential bioperator $(D_R\otimes D_L)$\/} on the bisemisheaf
$(\wh M^{(1)}_{v_R}\otimes\wh M^{(1)}_{v_L})$.

Let 
\[ (D_R\otimes D_L)(\phi (G^{(1)}(F_{\o v_j}\times F_{v_j}))
=E\RL(j)(\phi (G^{(1)}(F_{\o v_j}\times F_{v_j}))\;, \qquad 1\le j\le r\;,\]
be the eigenbivalue equation related to the bisemisheaf
$(\wh M^{(1)}_{v_R}\otimes\wh M^{(1)}_{v_L})$.
Then, we have that
\Be
\item {\bbf the $j$-th eigenbifunction $  \phi (G^{(1)}(F_{\o v_j}\times F_{v_j}))$ on the $j$ transcendental compact biquanta, being the $j$-th bisection of $(\wh M^{(1)}_{v_R}\otimes\wh M^{(1)}_{v_L})$, corresponds to the $j$-th eigenbivalue $E\RL(j)$ which is the shift of this bifunction and the shift of the global Hecke character\/} associated with this subbisemilattice.

\item the eigenbivalues of the matrix of $\GL_r(\rit\times\rit)$, constituting a representation of the bilinear differential Galois semigroup associated with the biaction of $(D_R\otimes D_L)$, are the eigenbivalues of the above eigenbivalue equation.
\Ee

\item {\bbf The second question ``What is behind the partition and correlation functions between eigenvalues of random matrices?'' concerns the distribution of eigenvalues \cite{D-H} of random matrices\/} of the Gaussian unitary (GUE) and orthogonal ensemble (GOE) {\bbf having $r$ quantum states and characterized by a Hamiltonian matrix of order $r$\/} whose entries are Gaussian random variables \cite{D-S}.

In the bilinear context envisaged in this paper, we are mostly interested by {\bbf the $m$-point correlation function for the bilinear Gaussian unitary (or orthogonal) ensemble BCUE (resp. BCOE) given by\/}:
\begin{align*}
R_{m_{r\RL}}(x^2_1,\dots,x^2_m)
&= \dfrac{r!}{(r-m)!}\int_{\rit^{r-m}} P_{r\RL}(x_1^2,\dots,x_r^2)\ dx_{m+1}\dots dx_r\\[11pt]
&= \det (K_r(x_k,x_\ell)^m_{k,\ell=1})\end{align*}
with
\begin{align*}
K_r(x_k,x_\ell)&= \sum^{r-1}_{i=0}\psi _i(x_k)\psi _i(x_\ell)\\
\and\qquad
\psi _i(x) &= h^{-1/2}P_i(x)\ e^{-r[(TG^T\times TG)(\rit\times\rit)]/2}\;.\end{align*}
{\bbf $P_i(x)$ being an orthogonal polynomial of degree $i$ corresponding to the weight function $e^{-r[(TG^T\times TG)(\rit\times\rit)]/2}$ where
$G(\rit\times\rit)=TG^T(\rit)\times TG(\rit)$ is  the bilinear Gauss decomposition of the matrix $G$ of BCOE\/}.

$R_{m_{r\RL}}(x^2_1,\dots,x^2_m)$ is the probability of finding a level around each of the bipoints (i.e. entries in $G$) $x_1^2,\dots,x_m^2$, the positions of the remaining levels being unobserved.

Let 
\[ K_r(x,x)=\sum^{r-1}_{i=0}\psi _i(x)\psi _i(x)\]
be the energy level density with $\psi _i(x)$ given above.

Then, we have found that:
\Be
\item {\bbf the squares of the roots of the polynomial $P_i(x)$ correspond to the eigenbivalues of the product, right by left, $(U_{r_R}\times U_{r_L})$ of Hecke operators\/}.

\item {\bbf the weight $e^{-r[(TG^T\times\rit)\times (TG\times\rit)]/2}$ is a measure of the eigenbivalues of the random matrix $G\in \GL_r(\rit\times\rit)$ being a representation of the differential bilinear Galois semigroup\/}.
\Ee

The orthogonal polynomials $P_i(x)$ satisfy the three term recurrent relation
\[ \beta _{I+1}P_{i+1}(x) = (x-\alpha _i)P_i(x)-\beta _iP_{i-1}(x)\]
leading to a tower whose matricial form is
\[ xP=JP+\beta _iP_i\]
where $J$ is the Jacobi matrix and $P$ a column vector of polynomials of increasing degrees.

It was then proven that {\bbf the $i$ roots of $P_i(x)$ are the eigenvalues of the Jacobi symmetric matrix being a representation\/} of the Hecke operator $U_{r_i}$.

We are then led to the conclusions:
\Be
\item {\bbf The probabilistic interpretation of quantum (field) theories is related to the bilinear semigroup of automorphisms
$\Aut_k(F_{\o v})\times \Aut_k(F_{v})$ of compact transcendental biquanta generating a bisemilattice of these\/}.

\item {\bbf The $m$-point correlation function for BCUE (or BCOE)\linebreak $R_{m_r}(x_1^2,\dots,x_m^2)$ constitutes a representation of the bilinear semigroup of automorphisms 
$\Aut_k(\phi _R(F_{\o v}\times\rit))\times \Aut_k(\phi _L(F_{v}\times\rit))$ of bifunctions on shifted compact transcendental biquanta\/}.
\Ee

\item {\bbf The third question ``What interpretation can we give to the local spacings \cite{Gau} between the eigenvalues of large random matrices'' depends on the dynamical global program of Langlands\/} developed in chapter 2 and summarized in the first question: it results from the $r$-dimensional algebraic representation of $\GL_r(\rit\times\rit)$ associated with the action of the differential bioperator $(D_R\otimes D_L)$ on the bisemisheaf
$(\wh M^{(1)}_{v_R}\otimes\wh M^{(1)}_{v_L})$ and leading to the above mentioned eigenbivalue equation in such a way that:
\Be
\item {\bbf the consecutive spacings
\[ \delta E\RL(j)=E\RL(j+1)-E\RL(j)\;, \qquad 1\le j\le r\le \infty \;,\]
between the eigenbivalues of the random matrix $G\in \GL_r(\rit\times\rit)$ are infinitesimal bigenerators of one biquantum of the Lie subbisemialgebra 
$\gl_1(F_{\o v^1}\times F_{v^1} )$ of the bilinear parabolic unitary semigroup 
$P_1(F_{\o v^1}\times F_{v^1} )\in \GL_1 ( F_{\o v} \times F_{v})$ and correspond to the energies of one free biquantum from subbisemilattices of $(j+1)$ biquanta\/}.

\item the $k$-th consecutive spacings
\[\delta E^{(k)}\RL(j)=E\RL(j+k)-E\RL(j)\]
between the eigenbivalues of the random matrix $G\in\GL_r(\rit\times\rit)$ are the infinitesimal bigenerators on $k$ biquanta of the Lie subbisemialgebra
$\gl_1(F_{\o v^k}\times F_{v^k} )$ of the bilinear $k$-th semigroup
$\gl_1(F_{\o v^k}\times F_{v^k} )$ and correspond to the energies of $k$ free biquanta from subbisemilattices of $(j+k)$ biquanta.
\Ee

The consecutive spacings $\delta E\RL(j)$ between the eigenbivalues of the matrix $G$ of
$\GL_r(\rit\times\rit)$ decompose into:
\[ \delta E\RL (j)=\delta EF\RL(j)+\delta EV\RL(j)\]
where {\bbf $\delta EF\RL(j)$ and $\delta EV\RL(j)$ denote respectively the fixed (or constant) and variable consecutive spacings between the $r$ eigenbivalues of $G$\/}.

Then, we have that:
\Be
\item the consecutive spacings
\[ E^{BCUE}\RL(j+1)-E^{BCUE}\RL(j)
=\delta EV^{BCUE}\RL(j)\]
between the eigenbivalues $E^{BCUE}\RL(j+1)$ and $E^{BCUE}\RL(j)$ of a unitary random matrix of $U_r(\cit\times \cit)$ (or $O_r(\rit\times\rit)$) are the variable (unitary) infinitesimal bigenerators on one biquantum  on the envisaged Lie subbisemialgebra or the variable (unitary) energies
$\delta E^{BCUE}\RL(j)$ of one biquantum in  subbisemilattices of $(j+1)$ biquanta.

\item {\bbf the $k$-th consecutive spacings
\[ E^{BCUE}\RL(j+k)-E^{BCUE}\RL(j)
=\delta EV^{(k)BCUE}\RL(j)\]
between the eigenbivalues of a unitary random matrix of $U_r(\cit\times\cit)$ (or $O_r(\rit\times\rit)$), are the variable energies
$\delta E^{(k)BCUE}\RL(j)$ on $k$ biquanta in subbisemilattices of $(j+k)$ biquanta\/}.
\Ee

\item {\bbf The fourth question: ``What interpretation can we give to the spacings between the nontrivial zeros \cite{Zag}, \cite{Pol}, of the Riemann zeta function $\zeta (s)$?'' depends on the solution of the Riemann hypothesis \cite{Bom} proposed in \cite{Pie7}\/} and briefly recalled now.

{\bbf The $1D$-pseudounramified simple global elliptic $\Gamma _{\wh M^{(1)}_{v^T\RL}}$-bisemimodule\/}
\[\phi ^{(1),(nr)}\RL(x) =\sum_n
\L( \lambda ^{(nr)}(n)\ e^{-2\pi inx}\otimes_D\lambda ^{(nr)}(n)\ e^{+2\pi inx}\R)\;, \quad x\in\rit\;,\]
where $\lambda ^{(nr)}(n)$ is a global Hecke character,

{\bbf can be interpreted as the sum of products, right by left, of semicircles of level ``$n$'' on $n$ transcendental compact quanta and constitutes a cuspidal representation of the bilinear semigroup
$\GL_2(\Fnrvt)$\/}.

Let $\zeta _R(s_-)$ and $\zeta _L(s_+)$, $s_-=\sigma -i\tau $ and $s_+=\sigma +i\tau $, be the two zeta functions defined respectively in the lower and upper half planes.

They are (distribution) inverse space functions, i.e. energy functions on the variables $s_-$ and $s_+$ conjugate to the complex space variables $z^*\in\cit$ and $z\in \cit$ of the cusp forms $f_L(z)$ and $f_R(z^*)$ submitted respectively to the transform maps $f_L(z) \to \zeta _L(s_+)$ and $f_R(z^*)\to \zeta _R(s_-)$ \cite{Pie8}.

Then, the kernel $\Ker(H_{\phi \RL\to\zeta \RL})$ of the map:
\[
H_{\phi \RL\to\zeta \RL}:\qquad 2\phi ^{(1),(nr)}\RL(x)\quad\To\quad\zeta _R(s_-)\otimes_D\zeta _L(s_+)\]
is the set of squares of the trivial zeros of $\zeta _R(s_-)$, $\zeta _L(s_+)$ and $\zeta (s)$, corresponding to the degeneracies of the products, right by left, of circles
$2\lambda ^{(nr)}(n)e^{2\pi inx}$ on  $2n$ transcendental quanta.

By this way, {\bbf we have a one-to-one correspondence between the trivial zeros of $\zeta (s)$ and the degeneracies of circles belonging to symmetric towers  of  cuspidal conjugacy class representatives of bilinear complete algebraic semigroups\/}.

Then, it is proved that {\bbf the products of the pairs of the trivial zeros of the Riemann zeta functions $\zeta _R(s_-)$ and $\zeta _L(s_+)$ are mapped into the products of the corresponding pairs of the nontrivial zeros according to\/}:
\begin{align*}
\{D_{4n^2,i^2}\cdot \varepsilon _{4n^2}\}: \qquad
\{\det (\alpha _{4n^2})\}_n
&\To \{\det (D_{4n^2,i^2}\cdot \varepsilon _{4n^2}\cdot \alpha _{4n^2})_{ss}\}_n\\
\{ (-2n)\times(-2n) \} _n
&\To \{ \lambda ^{(nr)}_+(4n^2,i^2,E_{4n^2})\times \lambda ^{(nr)}_-(4n^2,i^2,E_{4n^2})
\}_n
\end{align*}
where:
\Bi
\item $\alpha _{4n^2}$ is the split Cartan subgroup element associated with the integer $2n$;

\item $D_{4n^2,i^2}$ is the coset representative of the Lie (bisemi)algebra of the decomposition bisemigroup acting on $\alpha _{4n^2,i^2}$;

\item $\varepsilon _{4n^2}$ is the infinitesimal bigenerator of the considered bisemialgebra.
\Ei

{\bbf Every root of this Lie bisemialgebra is determined by the eigenvalues
\[ \lambda ^{(nr)}_\pm(4n^2,i^2,E_{4n^2})=
\dfrac{1\pm i\sqrt{16n^2\cdot E_{4n^2}-1}}2\]
of $D_{4n^2,i^2}\cdot\varepsilon _{4n^2}\cdot\alpha _{4n^2}$ which are the nontrivial zeros of $\zeta (s)$\/} written compactly according to
$\L(\tfrac12+ij_j\R)$ and $\L(\tfrac12-ij_j\R)$, $j\leftrightarrow n$.

Then, we have that:
\Bean
\item {\bbf the consecutive spacings
\[ \delta _{\gamma _j}=\gamma _{j+1}-\gamma _j\;, \qquad j=1,2,\dots\;, \]
between the nontrivial zeros of $\zeta (s)$ are equivalently\/}:
\Be
\item[a)] {\bbf the infinitesimal generators on one quantum of the Lie subsemialgebra $\gl_1(F^{(nr)}_{v^1})$ (or $\gl_1(F^{(nr)}_{\o v^1})$)
where $F^{(nr)}_{v^1}$ is a pseudounramified compact transcendental extension\/} (see section 1.1) of the linear\linebreak parabolic unitary semigroup
$P_1(F^{(nr)}_{v^1}) \subset \GL_1(F^{(nr)}_v)\equiv T_1(F^{(nr)}_v)$ (or\linebreak $P_1(F^{(nr)}_{\o v^1})$;

\item[b)] {\bbf the energies of one free quantum in subsemilattices of $(j+1)$ quanta\/}.
\Ee

\item The $k$-th consecutive spacings
\[ \delta _j^{(k)}=\gamma _{j+k}-\gamma _j\]
between the  nontrivial zeros of $\zeta (s)$ are equivalently:
\Be
\item[a)] the infinitesimal generators on $k$ quanta of the Lie subsemialgebra\linebreak
$\gl_1(F^{(nr)}_{v^k})$ (or $\gl_1(F^{(nr)}_{\o v^k})$ ) of the $k$-th semigroup
$\GL_1(F^{(nr)}_{v^k})$ (or\linebreak $\GL_1(F^{(nr)}_{\o v^k}) $);

\item[b)] the energies of $k$ free quanta in subsemilattices of $(j+k)$ quanta.
\Ee
\Ee

{\bbf The fifth question ``What is the curious connection between the spacings of the nontrivial zeros of $\zeta (s)$ and the corresponding spacings between the eigenvalues of random matrices?'' \cite{Ke-S} finds  response in the following statements (propositions of chapter 3)\/}:
\Bi
\item  The consecutive spacings
\[ \delta \gamma _j=\gamma _{+1}-\gamma _j\;, \qquad j\in\nit\;, \quad 1\le j\le r<\infty \;,\]
between the nontrivial zeros of the Riemann zeta function $\zeta (s)$

correspond to the consecutive spacings
\[ \delta E^{(nr)}\Rl(j) = E^{(nr)}\Rl(j+1)-E^{(nr)}\Rl(j)\;, \]
($\Rl$ means right of left),

between the square roots of the pseudounramified eigenbivalues of a (large) random matrix of $\GL_r(\rit\times\rit)$ (or of $\GL_r(\cit\times\cit)$)

and are equivalently:
\Bean
\item the infinitesimal generators on one quantum of the Lie subsemialgebra 
$\gl_1(F^{(nr)}_{v^1})$ of the linear parabolic unitary semigroup $P_1(F^{(nr)}_{v^1})$;

\item the energies of one transcendental compact pseudounramified ($N=1$) quantum in subsemilattices of $(j+1)$ transcendental compact pseudounramified quanta.\Ee

\item  The set $\{\delta E^{(nr)}\Rl(j)\}^r_{j=1}$ of consecutive spacings between the square roots of the eigenbivalues of a random matrix of $\GL_r(\rit\times\rit)$ as well as the set $\{\delta \gamma _j\}_{j=1}^r$ of consecutive spacings between the nontrivial zeros of $\zeta (s)$ constitutes a representation of the differential inertia Galois semigroup associated with the action of the differential operator $D_L$ or $D_R$.

\item Let $\{\delta E\Rl(j)\}_j$ be the set of consecutive spacings between the square roots of the eigenbivalues of a random matrix of $\GL_r(\rit\times\rit)$ or between the eigenvalues of a random matrix of $\GL_r(\rit)$.

Then, there is a surjective map:
\[ IM_{E\to \gamma }: \qquad \{ \delta E\Rl \} _j\quad\To \quad\{\delta \gamma _j\}_j\]
of which kernel $\Ker [IM_{E\to\gamma }]$ is the set $\{\delta E\Rl(j)-\delta E^{(nr)}\Rl(j)\}_j$ of differences of consecutive spacings between the square roots of the pseudoramified and pseudounramified eigenbivalues of a random matrix of $\GL_r(\rit\times\rit)$, i.e. the energies of one compact transcendental pseudoramified ($N>2$) quantum in subsemilattices of $(j+1)$ transcendental pseudoramified quanta.

\item Finally, {\bbf we can gather the results of this paper in the following proposition (see 3.33)\/}.
\vskip 11pt

{\em {\bbf
Let $\delta \gamma _j^{(k)}=\gamma _{j+k}-\gamma _j$ denote the $k$-th consecutive spacings between the  nontrivial zeros of $\zeta (s)$.

Let\/}:
\Bi
\item $\delta E\Rl^{(k)}(j)=E\Rl(j+k)-E\Rl(j)$;

\item $\delta E\Rl^{(nr),(k)}(j)=E^{(nr)}\Rl(j+k)-E^{(nr)}\Rl(j)$;

\item $\delta EV\Rl^{(k),(nr),BCOE}(j)=E^{(nr),BCOE}\Rl(j+k)-E^{(nr),BCOE}\Rl(j)$, $1\le j\le r$, $k\le j$,
\Ei
{\bbf be the $k$-th consecutive spacings between respectively\/}:
\Bi
\item   the pseudoramified eigenvalues of a random matrix of $GL_r(\rit)$;

\item  the pseudounramified eigenvalues of a random matrix of $GL_r(\rit)$;

\item the pseudounramified eigenvalues of a random unitary matrix of $O_r(\rit)$.
\Ei

{\bbf
Then, we have\/}:
\Bean
\item  {\bbf $\delta \gamma _j^{(k)}=\delta E\Rl^{(nr),(k)}(j)$ which are equivalently\/}:
\Bean
\item[a)] {\bbf  the infinitesimal generators of $k$ quanta of the Lie subsemialgebra $\gl_1(F^{(nr)}_{v^k})$ of the linear $k$-th semigroup\linebreak
$\GL_1(F^{(nr)}_{v^k})\subset\GL_1(F^{(nr)}_{v})$\/};

\item[b)] {\bbf  the energies of $k$ transcendental compact pseudounramified ($N=1$) quanta in subsemilattices in $(j+k)$ transcendental pseudounramified quanta\/};

\item[c)] {\bbf a representation of the differential Galois semigroup\/} associated with the action of the differential operator $D_L$ or $D_R$ on a function on $k$ transcendental pseudounramified quanta;
\Ee

\item a surjective map:
\[ IM^{(k)}_{E\to\gamma }: \qquad \{ \delta E^{(k)}\Rl(j) \}_j \To \{\delta \gamma ^{(k)}_j\}_j\]
of which kernel is the set $\{\delta E^{(k)}\Rl(j)-\delta E^{nr,(k)}\Rl(j)\}_j$ of difference of $k$-th consecutive spacings between the pseudoramified and pseudounramified eigenvalues of a random matrix of $\GL_r(\rit)$;

\item a bijective map:
\[ IM^{(k)}_{\gamma \to E^{(nr)}_{BCOE}}: \qquad \{\delta \gamma ^{(k)}_j\}_j \To \{ \delta EV ^{(k),(nr),BCOE}\Rl(j) \}_j\]
%of which kernel is the set $\{\delta \gamma ^{(k)}F_j=\delta \gamma ^{(k)}_j-\delta \gamma ^{(k)}V_j\}_j$ of fixed $k$-th consecutive spacings between nontrivial zeros of $\zeta (s)$ 
where $\delta \gamma _j^{(k)}$ denotes a $k$-th (variable) consecutive spacing verifying
$\delta \gamma  ^{(k)}_j=\delta E^{(k),(nr),BCOE}\Rl(j)$.
\Ee }

\item It has been inferred for a long time that the nontrivial zeros of the Riemann zeta function are probably related to the eigenvalues of some wave dynamical system. The connection between these two fields at the light of the global program of Langlands can be summarised by the proposition~3.35:
\vskip 11pt

{\em {\bbf The squares of the nontrivial zeros $\gamma _j$ of the Riemann zeta function $\zeta (s)$ are pseudounramified eigenbivalues of the eigenbivalue biwave operation:
\[
(D_R\otimes D_L)
(\phi (G^{(1)}(F^{(nr)}_{\o v_j}\times F^{(nr)}_{v_j})))
=\gamma ^2_j(\phi (G^{(1)}(F^{(nr)}_{\o v_j}\times F^{(nr)}_{v_j})))\]
of which eigenbifunctions are the sections of the bisemisheaf
$(\wh M^{(1)}_{v_R} \otimes \wh M^{(1)}_{v_L}  )$ being interpreted as the internal stringfield of an elementary (bisemi)particle}} \cite{Pie8}.
\Ei

\Ee

\section{Universal algebraic structures of the global program of Langlands}

\subsection{Symmetric structures of the program of Langlands}

In analogy with the fundamental theorem of the Galois theory which establishes a one-to-one correspondence between the set of closed intermediate fields of a given finite extension $F$ of a number field $k$ and the set of all closed normal subgroups of the Galois group $\Aut_k F$, we consider {\bbf a set of increasing finite algebraic extensions of $k$\/} in one-to-one correspondence with the corresponding Galois (sub)groups.  The envisaged global program of Langlands is then constructed on $n$-dimensional representations of such Galois (sub)groups and is recalled in this chapter.

Since symmetric semiobjects have to be considered in any generality \cite{Pie4}, we take into account {\bbf the bilinear global version of the Langlands program\/} concerning the generation of general symmetric structures, i.e. {\bbf double symmetric towers of conjugacy class representatives of (bilinear) algebraic (semi)groups\/}.
\vskip 11pt

\subsection{Algebraic and transcendental symmetric extensions}

Let then the set $\wt F$ of finite algebraic extensions of a number field $k$ of characteristic $0$ be a set of symmetric splitting fields composed of the left and right real symmetric splitting semifields $\wt F^+_L$ and $\wt F^+_R$ given respectively by the sets of positive and symmetric negative simple real roots.
\vskip 11pt

Assume that the {\bbf set of all increasing \lr splitting semifields\/} \cite{Wei}:
\begin{align*}
 \wt F_{v_1} \subset \dots \subset
  &\wt F_{v_{j,m_j}} \subset \dots \subset
 \wt F_{v_{r,m_r}} \\
\rresp{  \wt F_{\o v_1} \subset \dots \subset
  &\wt F_{\o v_{j,m_j}} \subset \dots \subset
 \wt F_{\o v_{r,m_r}}} \end{align*}
 is a set of increasing \lr real algebraic extensions characterized by degrees:
 \[ [\wt F_{v_j}:k]=[\wt F_{\o v_j}:k]=*+j\ N \;, \qquad
 1\le j\le r < \infty \;, \]
 which are integers modulo $N$, where
 \Bi
 \item $*$ denotes an integer inferior to $N$ ($*=0$ for the zero class);
\item $m_j$ labels the multiplicity of the envisaged extension.
\Ei
\vskip 11pt

These algebraic extensions are then said to be pseudoramified in contrast with the pseudounramified extensions
$\{ \wt F^{(nr)}_{v_{j,m_j}}\}_{j,m_j}$
\resp{$\{ \wt F^{(nr)}_{\o v_{j,m_j}}\}_{j,m_j}$} which are characterized by their global residue degrees 
$f_{v_j}$ \resp{$f_{\o v_j}$}:
\begin{equation}
f_{v_j}=[\wt F^{(nr)}_{v_{j,m_j}}:k]=[\wt F^{(nr)}_{\o v_{j,m_j}}:k]=j
\tag{case $N=1$}
\end{equation}
\vskip 11pt

{\bbf The smallest \lr (pseudoramified) splitting (sub)semifield\/}
$\wt F_{v_1}$ \resp{$\wt F_{\o v_1}$}
characterized by an extension degree
\begin{equation}
[\wt F_{v_{1}}:k]=[\wt F_{\o v_{1}}:k]=N
\tag{case $j=1$}
\end{equation}
{\bbf was interpreted as being a \lr algebraic quantum, i.e. an irreducible closed \lr algebraic subset\/}.
\vskip 11pt

According to the fundamental theorem of Galois, there could exist a one-to-one correspondence between the intermediate subsemifields
\begin{align*}
 \wt F_{v_r} \supseteq \dots \subseteq
  &\wt F_{v_{j}} \supseteq \dots \subseteq
 \wt F_{v_{1}} \supseteq k\\
\rresp{  \wt F_{\o v_r} \supseteq \dots \supseteq
  &\wt F_{\o v_{j}} \supseteq \dots \supseteq
 \wt F_{\o v_{1}} \supseteq k} \end{align*}
 and the corresponding closed subsemigroups of the Galois semigroup 
 $\Gal (\wt F_L^+/k)$ \resp{$\Gal (\wt F_R^+/k)$}:
 \[ \{1\}
 \subseteq \Gal (\wt F_{v_1}/k) \subseteq \dots 
 \subseteq \Gal (\wt F_{v_j}/\wt F_{v_{j-1}}) \subseteq \dots 
 \subseteq \Gal (\wt F_{r}/\wt F _{r-1}) \;.\]
 But, the extension degree $[\wt F_{v_r}:k]=[\wt F_{\o v_r}:k]$ would then be the product of the above intermediate algebraic extension degrees which is generally not verified by the given extension degree $[\wt F_{v_r}:k]=[\wt F_{\o v_r}:k]=r\cdot N$ which belongs to the zero class of integers modulo $N$ required by the searched one-to-one correspondence between the representations of the associated Galois (sub)groups and the corresponding cuspidal representations \cite{Pie2}.  However, every algebraic extension 
 $\wt F_{v_j}$ \resp{$\wt F_{\o v_j}$}, $1\le j\le r\le \infty $, characterized by the degree
$[ \wt F_{v_j}:k] =[ \wt F_{\o v_j}:k]=j\cdot N$, can be decomposed according to the fundamental theorem of Galois.
\vskip 11pt

Let then 
$\{\wt F_{v_{j,m_j}}\}_{j,m_j}$ \resp{$\{\wt F_{\o v_{j,m_j}}\}_{j,m_j}$} be the set of increasing algebraic extensions characterized by degrees
$\{[\wt F_{v_{j,m_j}}:k]=j\cdot N\}_{j,m_j}$
\resp{$\{[\wt F_{\o v_{j,m_j}}:k]=j\cdot N\}_{j,m_j}$}
and whose global Weil (or Galois) subgroups \cite{Pie2} are
$\Gal (\wt F_{v_{j,m_j}}/k)$
\resp{$\Gal (\wt F_{\o v_{j,m_j}}/k)$}.

By {\bbf an isomorphism of compactification\/}
\[ c_{v_{j,m_j}}: \quad \wt F_{v_{j,m_j}}\To F_{v_{j,m_j}} \qquad
\rresp{c_{\o v_{j,m_j}}: \quad \wt F_{\o v_{j,m_j}}\To F_{\o v_{j,m_j}}} \;,\]
{\bbf each algebraic extensions\/} 
$\wt F_{v_{j,m_j}}$ \resp{$\wt F_{\o v_{j,m_j}}$} {\bbf is sent into its compact image\/}
$ F_{v_{j,m_j}}$ \resp{$  F_{\o v_{j,m_j}}$} {\bbf which is a closed compact subset of $\rit_+$ \resp{$\rit_-$}\/}.

{\bbf Each compact image\/}
$F_{v_{j,m_j}}$ 
\resp{$F_{\o v_{j,m_j}}$} of the algebraic extension
$\wt F_{v_{j,m_j}}$
\resp{$\wt F_{\o v_{j,m_j}}$} {\bbf is thus a transcendental extension \cite{Hun} of which transcendence degree\/} % (or global residue degree)
$\tr \cdot d\cdot F_{v_{j,m_j}}$
\resp{$\tr \cdot d\cdot F_{\o v_{j,m_j}}$} is  given by:
\[ 
\tr \cdot d\cdot F_{v_{j,m_j}}=[\wt F_{v_{j,m_j}}:k]=j\cdot N \quad
\rresp{\tr \cdot d\cdot F_{\o v_{j,m_j}}=[\wt F_{\o v_{j,m_j}}:k]=j\cdot N }.
\] Every transcendental extension
$F_{v_{j,m_j}}$
\resp{$F_{\o v_{j,m_j}}$} is also an archimedean completion which can be viewed as resulting from the semigroup of automorphisms
$\Aut _k(F_{v_{j,m_j}})$
\resp{$\Aut _k(F_{\o v_{j,m_j}})$} which is a semigroup of reflections (or permutations) of a transcendental quantum
$F_{v^1_{j,m_j}} \subset F_{v_{j,m_j}}$
\resp{$F_{\o v^1_{j,m_j}} \subset F_{\o v_{j,m_j}}$} characterized by the corresponding extension degree
$[\wt F_{v^1_{j,m_j}}:k]=[\wt F_{\o v^1_{j,m_j}}:k]=N$ \cite{Pie5}.
\vskip 11pt

This implies that {\bbf the compact transcendental extension 
$F_{v_{h,m_h}}$
\resp{$F_{\o v_{h,m_h}}$}, composed of $h$ \lr transcendental quanta, of transcendental degree $h\cdot N$ over $k$\/}, is included into the transcendental extension
$F_{v_{j,m_j}}$
\resp{$F_{\o v_{j,m_j}}$}, composed of $j$ \lr compact transcendental quanta, in the sense of the fundamental theorem of Galois theory.  Indeed, the transcendental quantum
$F_{v^1_{j}} \subset F_{v_{j,m_j}}$
\resp{$F_{\o v^1_{j}} \subset F_{\o v_{j,m_j}}$} of the transcendental extension
$F_{v_{j,m_j}}$
\resp{$F_{\o v_{j,m_j}}$} is homeomorphic to the transcendental quantum
$F_{v^1_{h}} \subset F_{v_{h}}$
\resp{$F_{\o v^1_{h}} \subset F_{\o v_{h,m_h}}$}
of the transcendental extension
$F_{v_{h,m_h}}$
\resp{$F_{\o v_{h,m_h}}$}
since, by Galois automorphism and isomorphisms of compactification,
they are the compact images either of the nonunits of the algebraic quantum
$\wt F_{v^1_j}$
\resp{$\wt F_{\o v^1_j}$} of the algebraic extension
$F_{v_{j,m_j}}$
\resp{$F_{\o v_{j,m_j}}$}
or of the nonunits of the algebraic quantum
$F_{v^1_{h}}$
\resp{$F_{\o v^1_{h}}$} of the algebraic extension
$F_{v_{h,m_h}}$
\resp{$F_{\o v_{h,m_h}}$}.

Thus, in the case of transcendental extensions or compact archimedean completions, we have that:
\vskip 11pt

\subsection{Proposition}

{\em
If $F_{v_r}$ \resp{$F_{\o v_r}$} is a finite dimensional compact transcendental extension of $k$, {\bbf there is a one-to-one correspondence between the set of intermediate semifields
\begin{align*}
F_{v_r}\supset \dots \supset F_{v_j} \supset \dots \supset & F_{v_1}\supset k
\\
\rresp{F_{\o v_r}\supset \dots \supset F_{\o v_j} \supset \dots \supset & F_{\o v_1}\supset k}
\end{align*}
and the set of corresponding subsemigroups of the semigroup of automorphisms\linebreak
$\Aut_k(F_{v_r})$
\resp{$\Aut_k(F_{\o v_r})$}:\/}
\begin{align*}
\{1\} \subseteq \Aut_k(F_{v_1}) \subset \dots \subset & \Aut_{F_{v_h}}(F_{v_j}) \subset \dots \subset \Aut_{F_{v_m}}(F_{v_r}) \\
\rresp{ \{1\} \subseteq \Aut_k(F_{\o v_1}) \subset \dots \subset & \Aut_{F_{\o v_h}}(F_{\o v_j}) \subset \dots \subset \Aut_{F_{\o v_m}}(F_{\o v_r}) }
\end{align*}
in such a way that:
\Bean
\item the relative degree (or transcendental dimension) of two intermediate semifields must be an integer;

\item $F_{v_r}$ \resp{$F_{\o v_r}$} is transcendental over every intermediate semifield 
$F_{v_j}$
\resp{$F_{\o v_j}$}
if their transcendence degrees or global residue degrees verify \cite{Hun}:
\begin{align*}
\tr \cdot d\cdot F_{v_r}/k &= (\tr\cdot d\cdot F_{v_r}/F_{v_j})+(\tr\cdot d\cdot F_{v_j}/k)\\
\rresp{\tr \cdot d\cdot F_{\o v_r}/k &= (\tr\cdot d\cdot F_{\o v_r}/F_{\o v_j})+(\tr\cdot d\cdot F_{\o v_j}/k)}.
\end{align*}
\Ee
}

\bpr
In the Galois (i.e. algebraic) case, we would have that
$\Gal(\wt F_{v_r}/\wt F_{v_j})$, corresponding to $\Aut_{F_{v_j}} F_{v_r}$ in the transcendental case, is normal in
$\Gal(\wt F_{v_r}/k)$ and that
\[ \Gal (\wt F_{v_j}/k)
= \Gal (\wt F_{v_r}/k)\big/ (\wt F_{v_r}/F_{v_j})\]
implying the relative degree $\L[ \wt F_{v_r}:\wt F_{v_j}\R]=\frac rj$, being the relative index of the corresponding semigroup, is not an integer unless $j$ divides $r$.%  For that reason, the transcendental case is more general than the algebraic one taking into account that the Galois group is solvable if each maximal relative index is a prime number.
\epr
\vskip 11pt

\subsection{Bisemilattices of algebraic and transcendental quanta}

According to section~1.1. and \cite{Pie2}, only symmetric semiobjects have to be considered in any generality in such a way that {\bbf a bisemiobject $O_R\times O_L$ is composed of the product of the right semiobject $O_R$, localized in the lower half space or on $\rit_-$, and of the left symmetric semiobject $O_L$, localized in the upper half space or on $\rit_+$\/}.  So, instead of envisaging the set of \lr increasing transcendental (or algebraic) extensions:
\[ 
F_v = \{ F_{v_1},\dots,F_{v_{j,m_j}},\dots,F_{v_{r,m_r}}\}\quad
\rresp{F_{\o v} = \{ F_{\o v_1},\dots,F_{\o v_{j,m_j}},\dots,F_{\o v_{r,m_r}}\}};\]
we shall take into account their (diagonal) product
\[
F_{\o v}\times F_v :\{F_{\o v_1}\times F_{v_1},\dots,F_{\o v_{j,m_j}}\times F_{v_{j,m_j}},\dots,F_{\o v_{r,m_r}}\times F_{v_{r,m_r}}\}\]
in one-to-one correspondence with their semigroups of automorphisms:
\begin{multline*}
\Aut_k(F_{\o v})\times\Aut_k(F_{v})=\\
\{
\Aut_k F_{\o v_1}\times\Aut_k F_{v_1},\dots,
\Aut_k F_{\o v_{j,m_j}}\times\Aut_k F_{v_{j,m_j}},\dots,\\
\Aut_k F_{\o v_{r,m_r}}\times\Aut_k F_{v_{r,m_r}}\}\;.
\end{multline*}
Let
\begin{multline*}
\Gal (\wt F_{\o v}/k)\times\Gal (\wt F_{v}/k)=\\
\{\Gal (\wt F_{\o v_1}/k)\times\Gal (\wt F_{v_1}/k),\dots,
\{\Gal (\wt F_{\o v_{j,m_j}}/k)\times\Gal (\wt F_{v_{j,m_j}}/k),\dots,\\
\{\Gal (\wt F_{\o v_{r,m_r}}/k)\times\Gal (\wt F_{v_{r,m_r}}/k)\}
\end{multline*}
denote  the {\bbf products, right by left, of the Galois semigroups\/} of the sets of increasing algebraic extensions characterized by increasing degrees belonging to the zero class of integers modulo $N$: they are thus, {\bbf Weil global semigroups \cite{Pie2}\/}.

We thus have the isomosphism:
\[ \Gal (\wt F_{\o v}/k)\times \Gal (\wt F_{v}/k)
\quad \xrightarrow{\ \sim\ } \quad \Aut_k(F_{\o v})\times\Aut_k(F_{v})
\]
associated with the isomorphism:
\[ \wt F_{\o v}\times \wt F_{v} \quad \xrightarrow{\ \sim\ }\quad
 F_{\o v}\times F_{v} \]
 between the product, right by left, of sets of symmetric algebraic and transcendental extensions (or archimedean completions), which consists in the {\bbf the product, right by left, of two symmetric towers of increasing algebraic or compact transcendental extensions composed of an increasing number of algebraic or transcendental compact quanta\/}: this defines a bi(semi)lattice of biquanta, i.e. the product, right by left, of two symmetric (semi)lattices of algebraic or transcendental quanta.
 \vskip 11pt
 
\subsection{Abstract bisemivarieties}

Let $B_{\wt F_{v}}$
\resp{$B_{\wt F_{\o v}}$} be a \lr division semialgebra of real dimension $2n$ over the set $\wt F_v$ 
\resp{$\wt F_{\o v}$} of increasing real pseudoramified extensions:

$B_{\wt F_{v}}$
\resp{$B_{\wt F_{\o v}}$} is then isomorphic to the semialgebra of Borel upper (resp. lower) triangular matrices:
\[
B_{\wt F_{v}}\simeq T_{2n}(\wt F_v) \qquad 
\rresp{B_{\wt F_{\o v}}\simeq T^t_{2n}(\wt F_{\o v}} )
\] 
allowing introducing the algebraic bilinear semigroup of matrices by:
\[
B_{\wt F_{\o v}} \otimes B_{\wt F_{v}}
\simeq
T^t_{2n}(\wt F_{\o v}) \times T_{2n}(\wt F_v) 
\simeq \GL_{2n}(\wt F_{\o v}\times\wt F_v)\]
in such a way that its representation (bisemi)space is given by the tensor product
$\wt M^{(2n)}_{v_R}\otimes\wt M^{(2n)}_{v_L}$ of a right
$T^t_{2n}(\wt F_{\o v})$-semimodule $\wt M^{(2n)}_{v_R}$, localized in the lower half space, by a left
$T_{2n}(\wt F_{v})$-semimodule $\wt M^{(2n)}_{v_L}$, localized in the upper half space.
\vskip 11pt

{\bbf The $\GL_{2n}(\wt F_{\o v}\times \wt F_v)$-bisemimodule
$\wt M^{(2n)}_{v_R}\otimes\wt M^{(2n)}_{v_L}$ is an algebraic bilinear (affine) semigroup noted
$G^{(2n)}(\wt F_{\o v} \times\wt F_v)$\/} whose bilinear semigroup of Galois automorphisms is $\GL(\wh M^{(2n)}_{\o v_R}\otimes \wh M^{(2n)}_{v_L})$.

$\GL(\wh M^{(2n)}_{\o v_R}\otimes \wh M^{(2n)}_{v_L})$ constitutes the $2n$-dimensional equivalent of the product
$\Gal (\wt F_{\o v}/k)\times \Gal (\wt F_{v}/k)$ of Galois semigroups in such a way that
$G^{(2n)} (\wt F_{\o v}\times \wt F_{v})$ becomes the $2n$-dimensional (irreducible) representation space
$\Irr\Rep^{(2n)}_{\Gal_{F^+\RL}}(\Gal (\wt F_{\o v}/k)\times \Gal (\wt F_{v}/k))$
of
$\Gal (\wt F_{\o v}/k)\times \Gal (\wt F_{v}/k)$.
\[ G^{(2n)}(\wt F_{\o v}\times\wt F_{v})
= \Irr\Rep^{(2n)}_{\Gal_{F^+\RL}}(\Gal (\wt F_{\o v}/k)\times \Gal (\wt F_{v}/k))
\] implies the mononorphims:
\[
\sigma _{\wt v_R}\times\sigma _{\wt v_L}: \qquad
(\Gal (\wt F_{\o v}/k)\times \Gal (\wt F_{v}/k))\quad \To\quad
\GL(\wh M_{\o v_R}\otimes \wh M_{v_L})\approx 
G^{(2n)}(\wt F_{\o v}\otimes \wt F_{v})\;.\]
\vskip 11pt

The isomorphism
\[
\Gal (\wt F_{\o v}/k)\times \Gal (\wt F_{v}/k)
%\quad \xrightarrow{\sim} \quad
\quad \overset{\sim}{\To} \quad 
\Aut_k(F_{\o v})\times\Aut_k(F_{v})\]
between the products, right by left, of Weil global semigroups and corresponding automorphism semigroups of transcendental extensions leads to the commutative diagram
\[\begin{psmatrix}[colsep=1.2cm,rowsep=1.2cm]
\Gal (\wt F_{\o v}/k)\times \Gal (\wt F_{v}/k)
& \quad &
G^{(2n)} (\wt F_{\o v}\times\wt F_{v})\quad\\
\Aut_k(F_{\o v})\times\Aut_k(F_{v}) 
& \quad & 
G^{(2n)} (\wt F_{\o v}\times\wt F_{v})\quad
\everypsbox{\scriptstyle}
\ncline[arrows=->,nodesep=10pt]{1,1}{1,3}^{\quad\sigma _{\wt v_R}\times\sigma _{\wt v_L}}
\ncline[arrows=->,nodesep=10pt]{2,1}{2,3}^{\quad\sigma _{v_R}\times\sigma _{v_L}}
\ncline[arrows=->,nodesep=5pt]{1,1}{2,1}>{\lbag\;}
\ncline[arrows=->,nodesep=5pt]{1,3}{2,3}>{\lbag\;}
\end{psmatrix}\]
where {\bbf the monomorphism $\sigma _{v_R}\times\sigma _{v_L}$ generates the abstract bisemivariety\linebreak
$G^{(2n)} (F_{\o v}\times F_{v})$ on the product of sets of symmetric compact transcendental extensions (or archimedean completions)\/}.

The abstract bisemivariety
$G^{(2n)} ( F_{\o v}\times  F_{v})$ is covered by the algebraic bilinear (affine) semigroup 
$G^{(2n)} (\wt F_{\o v}\times\wt F_{v})$ and is thus a complete (locally) compact (algebraic) bilinear semigroup.

{\bbf At every infinite biplace
${\o v_j}\times{v}_j$ of $F_{\o v}\times F_{v}$ corresponds a conjugacy class
$g^{(2n)}_{v\RL}[j]$ of the abstract bisemivariety
$G^{(2n)} (\wt F_{\o v}\times\wt F_{v})$\/}.

The number of representatives of $g^{(2n)}_{v\RL}[j]$ corresponds to the number of equivalent extensions of $\wt F_{\o v_j}\times\wt F_{v_j}$;
\vskip 11pt

Let $\Os_{\wt F_v}$ \resp{$\Os_{\wt F_{\o v}}$} be the maximal order of
$\wt  F_{v}$
\resp{$\wt  F_{\o v}$}.

Then, $\Lambda _v=\Os_{B_{\wt F_v}}$
\resp{$\Lambda _{\o v}=\Os_{B_{\wt F_{\o v}}}$} in the division semialgebra
$B_{\wt  F_{v}}$
\resp{$B_{\wt  F_{\o v}}$} is a pseudo-ramified $\zit/N\ \zit$-lattice, in the \lr
$B_{\wt  F_{v}}$-semimodule $\wt M^{(2n)}_{v_L}$
\resp{$B_{\wt  F_{\o v}}$-semimodule $\wt M^{(2n)}_{v_R}$}.

So, we have that
\[
\Lambda _v \simeq T_{2n}(\Os_{\wt F_{v}}) \simeq T_{2n}(\Os_{F_v})\quad 
\rresp{\Lambda _{\o v} \simeq T^t_{2n}(\Os_{\wt F_{\o v}}) \simeq T^t_{2n}(\Os_{F_{\o v}})}
\]
leading to
\begin{align*}
\Lambda _{\o v}\times\Lambda _{v} 
&= T^t_{2n}(\Os_{F_{\o v}})\times T_{2n}(\Os_{F_{v}})\\
&= \GL_{2n}(\Os_{F_{\o v}}\times \Os_{F_{v}})\\
&= \GL_{2n}((\zit/N\ \zit)^2)\;.
\end{align*}
Then, the representation space
$\Repsp(\GL_{2n}(\zit/N\ \zit)^2)$ of
$\GL_{2n}((\zit/N\ \zit)^2)$ decomposes according to:
\[
\Repsp(\GL_{2n}(\zit/N\ \zit)^2)
= \bigoplus_j \bigoplus_{m_j} \L( \Lambda _{\o v_{j,m_j}}\otimes
\Lambda _{\o v_{j,m_j}}\R)\]
where
$\Lambda _{v_{j,m_j}}$
\resp{$\Lambda _{\o v_{j,m_j}}$} is the $(j,m_j)$-th subsemilattice referring to the conjugacy class representative
$g^{(2n)}_{v_L}[j,m_j]\in G^{(2n)}(F_{v})$
\resp{$g^{(2n)}_{v_R}[j,m_j]\in G^{(2n)}(F_{\o v})$}.

{\bbf The pseudo-ramified Hecke bisemialgebra $H\RL(2n)$\/} of all Hecke bioperators\linebreak
$T_R(2n;r)\otimes T_L(2n;r)$, {\bbf having a representation in the arithmetic subgroup of matrices $\GL_{2n}((\zit/N\ \zit)^2)$, generates the endomorphisms of the
$B_{F_{\o v}}\times B_{F_{v}}$-bisemimodule
$(M^{(2n)}_{v_R}\otimes M^{(2n)}_{v_L})$ decomposing it according to the bisubsemilattices
$( \Lambda _{\o v_{j,m_j}}\otimes
\Lambda _{\o v_{j,m_j}})$\/} \cite{Pie3}:
\[
M^{(2n)}_{v_R}\otimes M^{(2n)}_{v_L}=\bigoplus_{j,m_j}
(M^{(2n)}_{\o v_{j,m_j}}\otimes M^{(2n)}_{v_{j,m_j}})\;.\]
\vskip 11pt

Let
\begin{align*}
 F^{(nr)}_{v} &= \{F^{(nr)}_{v_{1}} , \dots , F^{(nr)}_{v_{j,m_j}} , \dots , F^{(nr)}_{v_{r,m_r}} \} \\
\rresp{ F^{(nr)}_{\o v} &= \{F^{(nr)}_{\o v_{1}} , \dots , F^{(nr)}_{\o v_{j,m_j}} , \dots , F^{(nr)}_{\o v_{r,m_r}} \} }\end{align*}
be the set of \lr increasing pseudounramified transcendental extensions homeomorphic to the corresponding pseudounramified algebraic extensions introduced in section~1.2.

Let $G^{(2n)}(F^{(nr)}_{\o v} \times F^{(nr)}_{v} )$ be the complete bilinear semigroup with entries in
$(F^{(nr)}_{\o v} \times F^{(nr)}_{v} )$.

Then, the kernel $\Ker (G^{(2n)}_{F\to F^{(nr)}})$ of the map:
\[
G^{(2n)}_{F\to F^{(nr)}}: \qquad
G^{(2n)}(F_{\o v}\times F_{v}) \quad \To \quad
G^{(2n)}(F^{(nr)}_{\o v}\times F^{(nr)}_{v})\]
is the smallest bilinear normal pseudoramified subgroup of
$G^{(2n)}(F_{\o v}\times F_{v})$:
\[ 
\Ker (G^{(2n)}_{F\to F^{(nr)}})=P^{(2n)}(F_{\o v^1}\times F_{v^1})\;, \]
i.e. {\bbf the parabolic bilinear subsemigroup\/} over the product
$F_{\o v^1}\times F_{v^1}$ of sets
\[
F_{\o v^1}=\{ F_{\o v^1_1} , \dots, F_{\o v^1_{j,m_j}} , \dots, 
 F_{\o v^1_{r,m_r}}\}\qquad
 \text{and} \qquad
 F_{v^1}=\{ F_{v^1_1} , \dots, F_{v^1_{j,m_j}} , \dots, 
 F_{v^1_{r,m_r}}\}\]
 {\bbf of ``unitary'' transcendental preudoramified extensions\/}.
\vskip 11pt

\subsection{Proposition}

{\em
The bilinear abstract parabolic semigroup
$P^{(2n)}(F_{\o v^1}\times F_{v^1})$ is the unitary (irreducible) representation space of the complete bilinear semigroup
$\GL_{2n}(F_{\o v}\times  F_v)$ of matrices.
}
\vskip 11pt

\bpr
\Bena
\item The abstract bisemivariety
$G^{(2n)}(F_{\o v}\times F_{v})$, covered by the algebraic (affine) bisemivariety
$G^{(2n)}(\wt F_{\o v}\times \wt F_{v})$, acts by conjugation on the bilinear parabolic subsemigroup
$P^{(2n)}(F_{\o v^1}\times F_{v^1})$ in such a way that the number of conjugates of
$P^{(2n)}(F_{\o v^1_j}\times F_{v^1_j})$ in the conjugacy class representative
$G^{(2n)}( F_{\o v_{j,m_j}}\times  F_{v_{j,m_j}})
\equiv  g^{(2n)}_{v\RL}[j,m_j]\in G^{(2n)} (  F_{\o v}\times  F_{v})$
is the index
\[ \L|
G^{(2n)}(  F_{\o v_{j,m_j}}\times   F_{v_{j,m_j}}):
P^{(2n)}(F_{\o v^1_j}\times F_{v^1_j}\R| = j\]
of the normalizer
$P^{(2n)}(F_{\o v^1_j}\times F_{v^1_j})$ in $G^{(2n)}(  F_{\o v_{j,m_j}}\times  F_{v_{j,m_j}})$.

\item Let 
$\Out (G^{(2n)}(  F_{\o v}\times   F_{v}))=
\Aut (G^{(2n)}(  F_{\o v}\times   F_{v}))\big/
\Int (G^{(2n)}(  F_{\o v}\times   F_{v}))$
be the bisemigroup \cite{Pie4} of transcendental automorphisms of the complete bilinear semigroup
$G^{(2n)}(  F_{\o v}\times   F_{v})$ where
$\Int (G^{(2n)}(  F_{\o v}\times   F_{v}))$ is the bisemigroup of transcendental inner automorphisms.

As, we have that:
\[
\Int (G^{(2n)}(  F_{\o v}\times  F_{v}))
=\Aut (P^{(2n)}(F_{\o v^1}\times F_{v^1}))\]
and considering 1), it appears that
$P^{(2n)}(F_{\o v^1}\times F_{v^1})$ is the unitary representation bisemispace with respect to the abstract bisemivariety
$G^{(2n)}( F_{\o v}\times  F_{v})$.\epr
\Ee
\vskip 11pt

\subsection{Corollary}

{\em
The rank of the bilinear parabolic semigroup
$P^{(2n)}(F_{\o v^1_j}\times F_{v^1_j})$ restricted to the $j$-th conjugacy class is
\[ r_{P^{(2n)}(F_{\o v^1_j}\times F_{v^1_j})}\simeq
(m_j\cdot N)^n\times_{(D)}(m_j\cdot N)^n\]
and the rank of the conjugacy class representative
$G^{(2n)}(  F_{\o v_{j,m_j}}\times  F_{v_{j,m_j}})$ is
\[
r_{G^{(2n)}(  F_{\o v_{j,m_j}}\times   F_{v_{j,m_j}})}
\simeq (j\cdot m_j\cdot N)^n \times_{(D)} (j\cdot m_j\cdot N)^n \]
where $\times_{D} $ is the notation for a diagonal product \cite{Pie4}.
}
\vskip 11pt

\bpr
As the bilinear parabolic semigroup
$P^{(2n)}(F_{\o v^1_j}\times F_{v^1_j})$ restricted to the $j$-th conjugacy class is the $2n$-dimensional representation space of the product, right by left, of inertia subgroups
\[
I_{F_{\o v_j}}\times I_{F_{v_j}}
= [ \Gal (F_{\o v_{j}}/k) \times \Gal (F_{v_{j}}/k) ]\ \Big/\ 
 [ \Gal (F^{(nr)}_{\o v_{j}}/k) \times \Gal (F^{(nr)}_{v_{j}}/k) ]\]
 according to \cite{Pie2} and as the order of 
 $I_{F_{\o v_j}}$ or $ I_{F_{v_j}}$ is $N$, we have that
 \[
 r_{P^{(2n)}(F_{\o v^1_j}\times F_{v^1_j})}
 = (m_j\cdot N)^n \times_{(D)} (m_j\cdot N)^n
 \]
 if it is taken into account that ``$m_j$'' real equivalent conjugacy class representatives $g^{(n)}_{v\Rl}[j]$ of dimension ``$n$'' cover one complex conjugacy class representative $g^{(2n)}_{\omega \Rl}(j)$ \cite{Pie5} of dimension $2n$:
\[g^{(2n)}_{\omega  \Rl}[j]=\{ g^{(n)}_{v \Rl}[j,m_j]\}_{m_j}
% g^{(2n)}_{v \RL}[j]\}
\;.\]
 
 It is then immediate to see that the rank of
$
G^{(2n)}(  F_{\o v_{j,m_j}}\times   F_{v_{j,m_j}})=g^{(2n)}\RL[j,m_j]$ is
\begin{equation}
 r_{G^{(2n)}(  F_{\o v_{j,m_j}}\times  F_{v_{j,m_j}})}
\simeq (j\cdot m_j\cdot N)^n \times_{(D)} (j\cdot m_j\cdot N)^n\tag*{\eop}
\end{equation}
\vskip 11pt

\subsection{Corollary}
{\em
\Bena
\item The number of transcendental quanta of the bilinear parabolic semigroup
$P^{(2n)}(F_{\o v^1_j}\times F_{v^1_j})$, $\forall\ j\in \nit$, is
\[ n(P^{(2n)}(F_{\o v^1_j}\times F_{v^1_j}))
=\ung^n \times_{(D)} \ung^n \To (m_j)^n \times_{(D)} (m_j)^n\;.\]

\item The number of transcendental quanta of the conjugacy class representative\linebreak
$G^{(2n)}( F_{\o v_{j,m_j}}\times  F_{v_{j,m_j}})$ is
\[ n(G^{(2n)}(F_{\o v_{j,m_j}}\times F_{v_{j,m_j}}))
= (j\cdot m_j)^n \times_{(D)} (j\cdot m_j)^n\;.\]
\Ee
}
\vskip 11pt

\bpr
This results from corollary~1.7 by taking into account that the degree of a transcendental quantum is $N$.  According to section~1.5, we thus have {\bbf a bisemilattice of $\sum\limits_{j}(j\times m_j)^n$ transcendental biquanta\/}.
\epr
\vskip 11pt

\subsection{Covering of complex abstract bisemivarieties by real abstract bisemivarieties}

Let $M^{(2n)}_{\omega _R}\otimes M^{(2n)}_{\omega _L}$ denote a complex
$\GL_{n}( F_{\o \omega }\times  F_\omega )$-bisemimodule being the representation space of the complete bilinear semigroup of matrices
$\GL_{n}( F_{\o \omega }\times  F_\omega )$ of the product 
$(F_{\o \omega }\times  F_\omega )$, right by left, of complex transcendental extensions covered by their real equivalents
$( F_{\o v }\times  F_v )$ as developed in \cite{Pie2}.

If each conjugacy class representative
$G^{(2n)}( F_{\o \omega_j }\times  F_{\omega _j})$ of
$G^{(2n)}( F_{\o \omega }\times  F_\omega )\equiv  M^{(2n)}_{\omega _R}\otimes M^{(2n)}_{\omega _L}$ is unique in this $j$-th class, then
$G^{(2n)}( F_{\o \omega _j}\times  F_{\omega_j} )$ is covered by $m_j$ real equivalent conjugacy class representatives
$G^{(n)}( F_{\o v_{j,m_j} }\times  F_{v_{j,m_j}} )$ of
$G^{(n)}( F_{\o v }\times  F_{v} )$ \cite{Pie2}.

So, {\bbf the complex bipoints of
$G^{(2n)}( F_{\o \omega }\times  F_\omega )$ are in one-to-one correspondence with the real bipoints of
$G^{(2n)}( F_{\o v }\times  F_v )$ and we have the inclusion:
\[\frac{G^{(2n)}( F_{\o v }\times  F_v )}{G^{(n)}( F_{\o v }\times  F_v )}
\simeq M^{(2n)}_{v _R}\otimes M^{(2n)}_{v _L} \hookrightarrow 
M^{(2n)}_{\omega _R}\otimes M^{(2n)}_{\omega _L}\]
of the real abstract bisemivariety
$M^{(2n)}_{v _R}\otimes M^{(2n)}_{v _L}$ into the complex abstract bisemivariety
$M^{(2n)}_{\omega _R}\otimes M^{(2n)}_{\omega _L}$\/}.
\vskip 11pt

\subsection{Corner stone of the global program of Langlands}

The real abstract bisemivariety
$M^{(2n)}_{v _R}\otimes M^{(2n)}_{v _L}\equiv 
G^{(2n)}( F_{\o v }\times  F_v )$, being the representation space of the bilinear algebraic semigroup
$\GL_{2n}( F_{\o v }\times  F_v )$ of matrices, constitutes the $2n$-dimension representation space of the product, right by left,
$\Gal(\wt F_{\o v}/k)\times\Gal(\wt F_{v}/k)$ of Galois semigroups according to section~1.5. So, we get the isomorphism
$\GL(\wt M^{(2n)}_{v_R}\otimes\wt M^{(2n)}_{v_L})\simeq
\GL_{2n}(\wt F_{\o v}\times\wt F_{v})$.

So, $G^{(2n)}( F_{\o v }\times  F_v )$ constitutes the corner stone of the real global  correspondence of Langlands recalled in \cite{Pie5} since, by an isomormosphism of toroidal compactification,
$G^{(2n)}( F_{\o v }\times  F_v )$ is transformed into
{\bbf the cuspidal representation
$\Pi (\GL_{2n}( \wt F_{\o v }\times \wt F_v ))$ of the algebraic bilinear semigroup of matrices
$\GL_{2n}( \wt F_{\o v }\times \wt F_v )$\/}.
\vskip 11pt

\subsection{Reducible representations of
$\GL_{2n}(   F_{\o v }\times   F_v )$}

We have envisaged until now that the representation of the general bilinear semigroup
$\GL_{2n}( F_{\o v }\times  F_v )$ was irreducible.  But, in consideration of the future developments of this paper, it is useful to take into account its reducibility.

Let $2n=2_1+2_2+ \dots + 2_k + \dots +2_\ell$  be a partition of $2n$ labeling the reducible representation of
$T_{2n}(F_v)$
\resp{$T^t_{2n}(F_{\o v})$}.

Then, we have that:
\Bena
\item {\bbf the representation 
$\Rep(\GL_{2n}(   F_{\o v }\times   F_v ))$ of the general bilinear semigroup
$\GL_{2n}(   F_{\o v }\times   F_v )=
T^t_{2n}(F_{\o v}) \times T^t_{2n}(F_{v}) $ is orthogonally completely reducible\/} if it decomposes diagonally according to the direct sum of $2$-dimensional irreducible representations
$\Rep(\GL_{2_k}(   F_{\o v }\times F_v ))$:
\[
\Rep(\GL_{2n}(   F_{\o v }\times   F_v ))=
\mathop{\boxplus}\limits_{2_k=2}^{2n}
\Rep(\GL_{2_k}( F_{\o v }\times  F_v ))\]
{\bbf and it is nonorthogonally completely reducible\/} if it decomposes diagonally according to the direct sum of irreducible $2$-dimensional bilinear representations
$\Rep(\GL_{2_k}( F_{\o v }\times  F_v ))$ and offdiagonally according to the direct sum of irreducible mixed bilinear representations
$\Rep(T^t_{2_k}( F_{\o v })\times T_{2_\ell}( F_v ))$.
\Ee

Indeed, taking into account the existence of cross products in the definition of bilinear semigroups \cite{Pie4}, we have that the representation of
$\GL_{2n}( F_{\o v }\times  F_v )$ can be reduced to:
\begin{align*}
\Rep(\GL_{2n}( F_{\o v }\times  F_v ))
&= \L( \mathop{\boxplus}\limits_{2_k=2}^{2n}
\Rep ( T^t_{2_k}( F_{\o v })) \otimes
\mathop{\boxplus}\limits_{2_\ell=2}^{2n}
\Rep ( T_{2_\ell}( F_{v })) \R)\\[11pt]
&=  \mathop{\boxplus}\limits_{2_k=2}^{2n}
\Rep ( \GL_{2_k}( F_{\o v } \times F_{v }))
\mathop{\boxplus}\limits_{2_k\neq 2\ell}^{2n}
\Rep ( T^t_{2_k}( F_{\o v })) \times T_{2_\ell}( F_{v } ) ) \;
\end{align*}
If the mixed bilinear representations
$\mathop{\boxplus}\limits_{2_k\neq 2\ell}^{2n}
\Rep ( T^t_{2_k}( F_{\o v })) \times T_{2_\ell}( F_{v })$ are equal to $0$, then the above completely reducible nonorthogonal representation of
$\GL_{2n}( F_{\o v }\times  F_v )$ reduces to the orthogonal case.
\vskip 11pt

\section[Universal dynamical structures of the global program of Langlands]{Universal dynamical structures of the global\linebreak program of Langlands}

\subsection{Dynamical functional representation spaces of abstract bi\-semivarieties}

In order to generate dynamical $\GL(2n)$-bisemistructures \cite{Pie3} referring to the geometric program of Langlands, we have to take into account the differentiable functional representation space
$\FREPSP ( \GL_{2n} ( \Fvt))$ of the complete bilinear semigroup
$\GL_{2n} ( \Fvt)$, i.e. {\bbf a bisemisheaf
$\wh M^{(2n)}_{v_R}\otimes
\wh M^{(2n)}_{v_L}$ of differentiable bifunctions on the abstract bisemivariety
$G^{(2n)}(\Fvt)$\/} which is given by the product, right by left, of a right semisheaf $\wh M^{(2n)}_{v_R}$ of differentiable (co)functions on the abstract right semivariety $G^{(2n)}(F_{\o v})\equiv T^{(2n)}(F_{\o v})$ by a left semisheaf
$\wh M^{(2n)}_{v_L}$ of symmetric differentiable functions on the abstract left semivariety $G^{(2n)}(F_{v})\equiv T^{(2n)}(F_{v})$.

This functional representation space
$\FREPSP(\GL_{2n}(\Fvt))$ of bilinear geometric dimension $2n$ splits into:
\[
\FREPSP(\GL_{2n}(\Fvt))=
\FREPSP(\GL_{2k}(\Fvt)) \oplus
\FREPSP(\GL_{2n-2k}(\Fvt)) \]
in such a way that
$\FREPSP(\GL_{2k}(\Fvt))$, $k\le n$, be the functional representation space of bilinear geometric dimension $2k$ of the bilinear semigroup
$\GL_{2k}(\Fvt)$ on which acts {\bbf the elliptic bioperator
$D^{2k}_R\otimes D^{2k}_L$\/}, i.e. the product of a right linear differential elliptic operator $D^{2k}_R$ acting on $2k$ variables by its left equivalent
$D^{2k}_L$, by its biaction:
\[ D^{2k}_R\otimes D^{2k}_L: \quad
\FREPSP(\GL_{2n}(\Fvt))\To
\FREPSP(\GL_{2n[2k]} ( \Fvtrit ))\]
where $\FREPSP(\GL_{2n[2k]} ( \Fvtrit ))$ is the functional representation space of
$\GL_{2n}(\Fvt)$ fibered or shifted into ``$2k$'' bilinear geometric dimensions.

$F^{S_\rit}_v = (F_v\times \rit)$
\resp{$F^{S_\rit}_{\o v} = (F_{\o v}\times \rit)$} denotes the set of increasing \lr transcendental extensions (or archimedean completions) fibered or shifted by real numbers
\begin{align*}
F_v\times \rit &= \{ F^{S_\rit}_{v_1},\dots, F^{S_\rit}_{v_{j,m_j}},\dots, F^{S_\rit}_{v_{r,m_r}}\}\\
\rresp{F_{\o v}\times \rit &= \{ F^{S_\rit}_{\o v_1},\dots, F^{S_\rit}_{\o v_{j,m_j}},\dots, F^{S_\rit}_{\o v_{r,m_r}}\}}\end{align*}
from their unshifted equivalents, i.e. fibre bundles with vertical fibres $\rit$ %in such a way that the
%$F^{S_\rit}_{v_{j,m_j}}$ be continuous, and deals with a {\bbf difference semifield\/} \cite{K-L-M-P} {\bbf consisting in the semifield of unshifted completions and an isomorphism of this one onto the semifield of shifted completions\/}.
and basis $F_v$ \resp{$F_{\o v}$}.

Remark that the shifted transcendental extension
$F^{S_\rit}_{v_{j,m_j}}$
\resp{$F^{S_\rit}_{\o v_{j,m_j}}$} is composed of $j$ \lr fibered or shifted transcendental quanta.
\vskip 11pt

\subsection{Bilinear fibre of tangent bibundle}

According to chapter~3 of \cite{Pie3}, the shifted functional representation space\linebreak
$\FREPSP(\GL_{2n[2k]} ( \Fvtrit ))$ decomposes into:
\begin{multline*}
\FREPSP(\GL_{2n[2k]} ( \Fvtrit ))\\
=
\FREPSP(\GL_{2k} ( \Fvtrit ))
\oplus \FREPSP(\GL_{2n-2k} ( \Fvt ))
\end{multline*}
in such a way that
$\FREPSP(\GL_{2k} ( \Fvtrit ))$ is the
{\bbf total bisemispace
$(\Delta ^{2k}_R\times \Delta ^{2k}_L)$ of the tangent bibundle\/}
$\TAN (\wh M^{(2n)}_{v_R}\otimes
\wh M^{(2n)}_{v_L})$ to the bisemisheaf\linebreak
$(\wh M^{(2n)}_{v_R}\otimes
\wh M^{(2n)}_{v_L})\equiv \FREPSP (\GL_{2k}(\Fvt))$ 
and  is isomorphic to the adjoint (functional) representation space of 
$\GL_{2k}(\Fvt)$  corresponding to the action of the bioperator
$(D ^{2k}_R\times D ^{2k}_L)$ on 
$(\wh M^{(2n)}_{v_R}\otimes
\wh M^{(2n)}_{v_L})$ which maps 
$\TAN (\wh M ^{(2n)}_{v_R}\otimes\wh M ^{(2n)}_{v_L})\simeq
\AdF\REPSP (\GL_{2k}(\Fvt)$ into itself.

We have thus that:
\begin{align*}
\Delta ^{2k}_R\times \Delta ^{2k}_L
&\simeq
\AdF\REPSP (\GL_{2k}(\Fvt) \\
&\simeq
\fREPSP (\GL_{2k}(\Fvtrit))\;,\end{align*}
where $\fREPSP (\GL_{2k}(\rit\times\rit))$ is the bilinear fibre 
$\Fs\RL(\TAN)$ of $\TAN(\wh M^{(2n)}_{v_R}\otimes \wh M^{(2n)}_{v_L})$.
\vskip 11pt

\subsection{Symbol of the bioperator $(D ^{2k}_R\times D ^{2k}_L)$}

Referring to the classical definition \cite{A-S} of the symbol $\sigma (D)$ of a differential linear operator $D$ dealing with the unit sphere bundles in the cotangent vector bundle $T^*(X)$ of the compact smooth manifold $X$, we introduce the {\bbf symbol $\sigma (D ^{2k}_R\times D ^{2k}_L)$ of the differential bioperator
$(D ^{2k}_R\times D ^{2k}_L)$ by\/}:
\[
\sigma (D ^{2k}_R\times D ^{2k}_L)\simeq
\REPSP (P_{2k}(\rit\times\rit)_{\mid F_{\o v^1}\times F_{v^1}}))\;, \]
i.e. the  unitary functional representation space of
$\GL_{2k}((\rit\times\rit)_{\mid F_{\o v}\times F_{v}})$ given by 
{\bbf the  functional representation space of the fibering or shifting bilinear parabolic semigroup
$P_{2k}(\rit\times\rit)_{\mid F_{\o v^1}\times F_{v^1}})$\/}.

So, the differential bioperator
$D ^{2k}_R\times D ^{2k}_L$ maps from the bisemisheaf
\[
(\wh M^{(2n)}_{v_R}\otimes
\wh M^{(2n)}_{v_L})=\FREPSP (\GL_{2n}(\Fvt))\]
into the corresponding perverse bisemisheaf
\[
(\wh M^{(2n)}_{v_R}[2k]\otimes
\wh M^{(2n)}_{v_L}[2k])=\FREPSP (\GL_{2n[2k]}(\Fvtrit))\;,\]
shifted in $2k$ geometric dimensions, according to:
\[
D ^{2k}_R\times D ^{2k}_L : \quad
(\wh M^{(2n)}_{v_R}\otimes\wh M^{(2n)}_{v_L})\To 
(\wh M^{(2n)}_{v_R}[2k]\otimes\wh M^{(2n)}_{v_L}[2k])\;,\]
while 
$\sigma (D ^{2k}_R\times D ^{2k}_L)$ maps the ``unitary'' bisemisheaf
\begin{align*}
(\wh M^{(2n)}_{v^1_R}\otimes\wh M^{(2n)}_{v^1_L})
&=
\FREPSP(P_{2n}(F_{\o v^1}\times F_{v^1}))\\
&\equiv
\FREPSP(\GL_{2n}(F_{\o v^1}\times F_{v^1}))
\end{align*}
into the corresponding ``unitary'' perverse bisemisheaf
\begin{align*}
(\wh M^{(2n)}_{v^1_R}[2k]\otimes
\wh M^{(2n)}_{v^1_L}[2k])
&=\FREPSP (P_{2n[2k]}(F_{\o v^1}\times \rit)\times ( F_{v^1}\times \rit))\\
&\equiv \FREPSP (\GL_{2n[2k]}(F_{\o v^1}\times \rit)\times ( F_{v^1}\times \rit))
\end{align*}
according to:
\[
\sigma (D ^{2k}_R\times D ^{2k}_L ): \quad
(\wh M^{(2n)}_{v^1_R}\otimes\wh M^{(2n)}_{v^1_L})\To 
(\wh M^{(2n)}_{v^1_R}[2k]\otimes\wh M^{(2n)}_{v^1_L}[2k])\;.\]
\vskip 11pt

\subsection{Proposition}

{\em
\Be
\item 
The bilinear semigroup of matrices
$\GL_r(\rit^{(2k)}\times\rit^{(2k)})$ with algebraic orders ``$r$'',\linebreak referring to the biaction of the bioperator
$
(D ^{2k}_R\times D ^{2k}_L)$ on the real bisemisheaf
$(\wh M^{(2n)}_{v_R}\otimes\wh M^{(2n)}_{v_L})$ over the abstract real bisemivariety
$G^{(2k)}(\Fvt)\equiv\linebreak \REPSP (\GL_{2k}(\Fvt)$, corresponds to the $2k$-dimensional bilinear real representation of the product, right by left, of ``differential'' Galois (or global Weil) semigroups
$\Aut_k(\phi _R(\rit))\times \Aut_k(\phi _L(\rit))$ fibering or ``shifting'' the product, right by left, of automorphism semigroups
$\Aut_k(\phi _R(F_{\o v}))\times \Aut_k(\phi _L(F_v))$ of cofunctions 
$\phi _R(F_{\o v})$ and functions $\phi _L(F_v)$ respectively on the compact transcendental real extensions $F_{\o v}$ and $F_v$ by $
\Aut_k(\phi _R(F_{\o v}\times \rit))\times \Aut_k(\phi _L(F_v\times \rit))$.

So we have that:
\[
\Rep^{(2k)}( \Aut_k ( \phi _R ( \rit ))\times \Aut_k(\phi _L ( \rit )))
=\GL_r(\phi _R(\rit^{(2k)})\times\phi _L(\rit^{(2k)}))\;.\]

\item 
The unitary parabolic bilinear semigroup
$P_r(\rit^{(2k)}\times\rit^{(2k)})\subset \GL_r(\rit^{(2k)}\times\rit^{(2k)})$, referring to the biaction of the bioperator on the unitary bisemisheaf
$(\wh M^{(2n)}_{v^1_R}\otimes
\wh M^{(2n)}_{v^1_L})$ over the unitary abstract bisemivariety
$P^{(2k)} (F_{\o v^1}\times F_{v^1})=\REPSP (P_{2k} (F_{\o v^1}\times F_{v^1}))$, corresponds to the $2k$-dimensional bilinear representation of the product, right by left, of ``differential'' inertia Galois (or global Weil) semigroups
$\Int_k(\phi _R(\rit))\times\linebreak \Int_k(\phi _L(\rit))$ fibering or shifting the product, right by left, of internal automorphism semigroups
$\Int_k(\phi _R(F_{\o v^1}))\times \Int_k(\phi _L(F_{v^1}))$ of cofunctions 
$\phi _R(F_{\o v^1})$ and functions $\phi _L(F_{v^1})$ respectively on the  unitary compact transcendental real extensions $F_{\o v^1}$ and $F_{v^1}$ by $
\Int_k(\phi _R(F_{\o v^1}\times \rit))\times \Int_k(\phi _L(F_{v^1}\times \rit))$.

So we have that:
\[
\Rep^{(2k)} ( \Int_k ( \phi _R ( \rit ))\times \Int_k ( \phi _L ( \rit )))
=P_r(\phi _R(\rit^{(2k)})\times\phi _L(\rit^{(2k)}))\;.\]
\Ee
}
\vskip 11pt

\bpr
\Be
\item {\bbf The abstract bisemivariety $G^{(2k)}(\Fvt)\equiv
M^{(2k)}_{\o v_R}\otimes  M^{(2k)}_{v_L})$ decomposes according to the increasing filtration of its conjugacy class representatives\/}:
\begin{multline*}
G^{(2k)} (F_{\o v_1}\times F_{v_1}) \subset \dots\subset
G^{(2k)} (F_{\o v_{j,m_j}}\times F_{v_{j,m_j}}) \subset \dots \subset
G^{(2k)} (F_{\o v_{r,m_r}}\times F_{v_{r,m_r}}) \\
\equiv
M^{(2k)}_{\o v_1}\otimes M^{(2k)}_{v_1} \subset \dots \subset
M^{(2k)}_{\o v_{j,m_j}}\otimes M^{(2k)}_{v_{j,m_j}} \dots\subset 
M^{(2k)}_{\o v_{r,m_r}}\otimes M^{(2k)}_{v_{r,m_r}}\;,\\
1\le j\le r\le \infty 
\end{multline*}
in one-to-one correspondence with the filtration of the product
$(\Fvt)$ of sets of archimedean transcendental extensions (see section~1.2).

Similarly, the unitary abstract bisemivariety
$P^{(2k)}(F_{\o v^1}\times F_{v^1})\equiv
M^{(2k)}_{\o v^1_R}\otimes  M^{(2k)}_{v^1_L})$, which is a bilinear parabolic subsemigroup, decomposes according to the increasing set of its   representatives:
\begin{multline*}
P^{(2k)} (F_{\o v^1_1}\times F_{v^1_1}) \subset \dots \subset
P^{(2k)} (F_{\o v^1_{j,m_j}}\times F_{v^1_{j,m_j}}) \subset \dots \subset
P^{(2k)} (F_{\o v^1_{r,m_r}}\times F_{v^1_{r,m_r}}) \qquad\\
\qquad \equiv
M^{(2k)}_{\o v^1_1}\otimes M^{(2k)}_{v^1_1} \subset \dots \subset
M^{(2k)}_{\o v^1_{j,m_j}}\otimes M^{(2k)}_{v^1_{j,m_j}} \dots\subset 
M^{(2k)}_{\o v^1_{r,m_r}}\otimes M^{(2k)}_{v^1_{r,m_r}}\;,
\end{multline*}
in one-to-one correspondence with the filtration of the product
$(F_{\o v^1}\times F_{v^1})$ of sets of unitary transcendental extensions.
\vskip 11pt

%%%% 2

\item Referring to the ``infinite'' general bilinear semigroup given by
$\GL(\Fvt)=\underrightarrow{\lim}\GL_m(\Fvt)$ in such a way that
$\GL_m(\Fvt)$ embeds in $\GL_{n+1}(\Fvt)$ according to the geometric dimension ``$m$'', {\bbf it is possible to introduce as in \cite{Pie5} an ``infinite quantum'' general bilinear semigroup by
\[
\GL^{(Q)}(F^{(2k)}_{\o v^1}\times F^{(2k)}_{v^1})=
\lim_{j=1\to r} \GL^{(Q)}_j(F^{(2k)}_{\o v^1}\times F^{(2k)}_{v^1})\]
}
in such a way that:
\Be
\item $\GL^{(Q)}_1(F^{(2k)}_{\o v^1}\times F^{(2k)}_{v^1})$ is the parabolic, i.e. unitary, bilinear semigroup\linebreak
$\GL_{2k}(F_{\o v^1}\times F_{v^1})\equiv
P_{2k}(F_{\o v^1}\times F_{v^1})$.
\vskip 11pt
\item $\GL^{(Q)}_j(F^{(2k)}_{\o v^1}\times F^{(2k)}_{v^1})=
\GL_{2k}(F_{\o v_j}\times F_{v_j})\simeq
(F^{(2k)}_{\o v_j}\times F^{(2k)}_{v_j})$.
\vskip 11pt
\item $\GL^{(Q)}_j(F^{(2k)}_{\o v^1}\times F^{(2k)}_{v^1}) \subset
\GL^{(Q)}_{j+1}(F^{(2k)}_{\o v^1}\times F^{(2k)}_{v^1})$ where the integer ``$j$'' denotes a global residue degree, i.e. an algebraic dimension, while the integer ``$2k$'' refers to a geometric dimension.
\Ee

It was then proved in \cite{Pie5} that the ``infinite quantum'' general bilinear semigroup
$\GL^{(Q)}(F^{(2k)}_{\o v^1}\times F^{(2k)}_{v^1})$, defined with respect to the unitary parabolic bilinear semigroup
$\GL^{(Q)}_1(F^{(2k)}_{\o v^1}\times F^{(2k)}_{v^1})=
P_{2k}(F_{\o v^1}\times F_{v^1})$, is the general bilinear semigroup
$\GL_{2k}(\Fvt)$.

So, {\bbf the set
$\{ G^{(2k)}(F_{\o v_{j,m_j}}\times F_{v_{j,m_j}})\}_{j,m_j}$ of conjugacy class representatives of the bisemivariety
$G^{(2k)}(\Fvt)$\/}, generated from the bilinear semigroup of matrices
$\GL_{2k}(\Fvt)$, {\bbf can be rewritten according to the set
$\{ G_{(Q)}^{(j,m_j)}(F^{(2k)}_{\o v^1}\times F^{(2k)}_{v^1})\}_{j,m_j}$}
of increasing conjugacy class representatives of the bisemivariety\linebreak
$G_{(Q)}^{(\o v\times v)}(F^{(2k)}_{\o v^1}\times F^{(2k)}_{v^1})=
G^{(2k)}(\Fvt)$ generated from
$\GL^{(Q)}(F^{(2k)}_{\o v^1}\times F^{(2k)}_{v^1})$.

And, thus, the filtration
\[
G^{(1)}_{(Q)}(F^{(2k)}_{\o v^1}\times F^{(2k)}_{v^1}) \subset \dots \subset
G^{(j,m_j)}_{(Q)}(F^{(2k)}_{\o v^1}\times F^{(2k)}_{v^1}) \subset \dots \subset
G^{(r,m_r)}_{(Q)}(F^{(2k)}_{\o v^1}\times F^{(2k)}_{v^1})\]
of conjugacy class representatives of
$G^{(2k)}(\Fvt)$ is now written with reference to the increasing algebraic dimensions of the bilinear subsemigroups of matrices:
\[
\GL^{(Q)}_1(F^{(2k)}_{\o v^1}\times F^{(2k)}_{v^1}) \subset \dots \subset
\GL^{(Q)}_{j,m_j}(F^{(2k)}_{\o v^1}\times F^{(2k)}_{v^1}) \subset \dots \subset
\GL^{(Q)}_{r,m_r}(F^{(2k)}_{\o v^1}\times F^{(2k)}_{v^1}) \;.\]
\vskip 11pt

%% 3
\item Referring to section~2.3, we see that the set 
$\{\phi ( G^{(2k)}(F_{\o v_{j,m_j}}\times F_{v_{j,m_j}}))\}_{j,m_j}$ of conjugacy class representatives, or {\bbf bisections, of the bisemisheaf
$(\wh M^{(2k)}_{v_R}\otimes \wh M^{(2k)}_{v_L})=\FREPSP(\GL_{2k}(\Fvt))$\/}, which can be rewritten according to
$\{\phi (G^{(j,m_j)}_{(Q)}(F^{(2k)}_{\o v^{1}}\linebreak \times F^{(2k)}_{v^{1}}))\}_{j,m_j}$
with respect to $\GL^{(Q)}(F^{(2k)}_{\o v^{1}}\times F^{(2k)}_{v^{1}})$,
{\bbf are fibered or shifted, as tangent bibundles, under the action of the bioperator
$(D^{2k}_R\otimes D^{2k}_L)$ into\/}:
\begin{multline*}
D^{2k}_R\otimes D^{2k}_L: 
\begin{aligned}[t]
&(\wh M^{(2n)}_{v_R}\otimes \wh M^{(2n)}_{v_L})
 \To (\wh M^{(2n)}_{v_R}[2k]\otimes \wh M^{(2n)}_{v_L}[2k])\\
&\{\phi (G^{(j,m_j)}_{(Q)}(F^{(2k)}_{\o v^{1}}\times F^{(2k)}_{v^{1}}))\}_{j,m_j}
\end{aligned}\\
\To
\{\phi (G^{(j,m_j)}_{(Q)}(F^{(2k)}_{\o v^{1}}\times \rit^{(2k)})\times
( F^{(2k)}_{v^{1}}))\times \rit^{(2k)}\}_{j,m_j}\;.
\end{multline*}

Similarly,
the set 
$\{\phi (P^{(2k)}(F_{\o v^1_{j,m_j}}\times F_{v^1_{j,m_j}}))\}_{j,m_j}$ of  bisections of the unitary bisemisheaf
$(\wh M^{(2k)}_{v^1_R}\otimes \wh M^{(2k)}_{v^1_L})=\FREPSP(P_{2k}(F_{\o v^1}\times F_{v^1}))$, which can be rewritten according to
$\{\phi (P^{(j,m_j)}_{(Q)}(F^{(2k)}_{\o v^{1}}\times F^{(2k)}_{v^{1}}))\}_{j,m_j}$
 are fibered or shifted under the action of the bioperator
$(D^{2k}_R\otimes D^{2k}_L)$ into:
\begin{multline*}
D^{2k}_R\otimes D^{2k}_L: 
\begin{aligned}[t]
&(\wh M^{(2n)}_{v^1_R}\otimes \wh M^{(2n)}_{v^1_L})
\To (\wh M^{(2n)}_{v^1_R}[2k]\otimes \wh M^{(2n)}_{v^1_L}[2k])\\
&\{\phi (P^{(j,m_j)}_{(Q)}(F^{(2k)}_{\o v^{1}}\times F^{(2k)}_{v^{1}}))\}_{j,m_j}
\end{aligned}
\\
\To
\{\phi (P^{(j,m_j)}_{(Q)}(F^{(2k)}_{\o v^{1}}\times \rit^{(2k)})\times
( F^{(2k)}_{v^{1}}))\times \rit^{(2k)}\}_{j,m_j}\;.
\end{multline*}
\vskip 11pt

%% 4
\item 
Let the geometric dimension $2k$ be equal to $2k=1$.  The tower of shifted real transcendental biextensions is:
\begin{multline*}
(F_{\o v_1}\times \rit)\times (F_{v_1}\times \rit) \subset \dots \subset
(F_{\o v_{j,m_j}}\times \rit)\times (F_{v_{j,m_j}}\times \rit)\\
 \subset \dots \subset
(F_{\o v_{r,m_r}}\times \rit)\times (F_{v_{r,m_r}}\times \rit) \;,\end{multline*}
i.e. a tower of (real) shifted transcendental biquanta by vertical bifibres ($\rit\times\rit$).

They are the conjugacy class representatives generated under the action of representatives of
$\GL_1(\Fvtrit)$ and they can be rewritten as the filtration
\begin{multline*}
G^{(1)}_{(Q)}((F_{\o v^1}\times \rit)\times(F_{v^1}\times \rit)\subset\dots\subset
G^{(j,m_j)}_{(Q)}((F_{\o v^1}\times \rit)\times(F_{v^1}\times \rit)\\
\subset\dots\subset
G^{(r,m_r)}_{(Q)}((F_{\o v^1}\times \rit)\times(F_{v^1}\times \rit))
\end{multline*}
of the representatives of the bilinear subsemigroups of the (infinite) quantum general bilinear semigroup
$\GL^{(Q)}((F_{\o v^1}\times \rit)\times(F_{v^1}\times \rit))$.

{\bbf The bilinear semigroup of automorphisms of these fibered or shifted transcendental extensions is\/}:
\[
\Aut_k(F_{\o v}\times\rit)\times\Aut_k(F_{v}\times\rit)=
\{\dots,\Aut_k(F_{\o v_{j,m_j}}\times\rit)\times\Aut_k(F_{v_{j,m_j}}\times\rit),\dots\}\;.\]
So
\[ 
[ \Aut_k(F_{\o v}\times\rit)\times\Aut_k(F_{v}\times\rit) ]\big/
[\Aut_k(F_{\o v})\times\Aut_k(F_{v})=
\Aut_k(\rit)_{|F_{\o v}}\times\Aut_k(\rit)_{|F_v}
\]
corresponds to the bilinear semigroup of fibering or shifting automorphisms.

The functional representatives of
$\GL^{(Q)}((F_{\o v^1}\times \rit)\times(F_{v^1}\times \rit))$, or equivalently
of $\GL_1(\Fvtrit)$, are
\[
F\ G^{(1)}_{(Q)}((F_{\o v^1}\times \rit)\times(F_{v^1}\times \rit))\subset\dots\subset
F\ G^{(j,m_j)}_{(Q)}((F_{\o v^1}\times \rit)\times(F_{v^1}\times \rit))\;, \]
i.e. bifunctions on the birepresentatives
$\{G^{(j,m_j)}_{(Q)}((F_{\o v^1}\times \rit)\times(F_{v^1}\times \rit))\}_{j,m_j}$.

And, thus,
\[
\Aut_k(\phi _R(\rit))\times \Aut_k(\phi _L(\rit))=
\{\dots,\Aut_k(\phi _{j,m_j}(\rit))\times \Aut_k(\phi _{j,m_j}(\rit)),
\dots\}_{j,m_j}\;, \]
where $\phi _{j,m_j}(\rit)$ is the $(j,m_j)$-th function over $\rit$  acting on the function\linebreak
 $\phi _{j,m_j}(F_{v_{j,m_j}})$ over the transcendental extension
 $F_{v_{j,m_j}}$ or over the conjugacy class representative of
 $\GL_1(F_{v_{j,m_j}})$.
 
 So, {\bbf
 $\Aut_k(\phi _R(\rit))\times \Aut_k(\phi _L(\rit))$ is the bilinear differential Galois semigroup where
$ \Aut_k(\phi _R(\rit))$  is the set of linear differential Galois (semi)\-subgroups \cite{Car}, \cite{And}, acting on the set of sections of the considered $1D$-differential equation and
$ \Aut_k(\phi _L(\rit))$ is the set of linear differential Galois (semi)subgroups acting on the symmetric set of sections\/}.

Referring to the fundamental theorem of Galois theory, we see that
{\bbf the bilinear differential Galois semigroup
$\Aut_k(\phi _R(\rit))\times \Aut_k(\phi _L(\rit))$ corresponds to the upper bilinear differential Galois semigroup\/}
$\Aut_k(\phi _R(\rit))_{|F_{\o v_{r,m_r}}}\times\linebreak \Aut_k(\phi _L(\rit))_{|F_{v_{r,m_r}}}$
with respect to the upper algebraic dimension ``$r$''.

Similarly, the bilinear semigroup of unitary, i.e. internal, shifting automorphisms:
\[ 
\Int_k(\rit)_{|F_{\o v^1}}\times\Int_k(\rit)_{|F_{v^1}}=
[ \Int_k(F_{\o v^1}\times\rit)\times\Int_k(F_{v^1}\times\rit) ]\big/
[\Int_k(F_{\o v^1})\times\Int_k(F_{v^1})
\]
has for differential Galois bilinear semigroup
$\Int_k(\psi _R(\rit )) \times\Int_k(\psi _L(\rit))$ where\linebreak
{\bbf $\Int_k(\psi _L(\rit))$  is the unitary linear differential Galois (semi)group\/}, i.e. the inertia linear differential (semi)group, {\bbf acting on a ``unitary section'' of the envisaged $1D$-differential equation, and
$\Int_k(\psi _R(\rit))$  is the inertia linear differential (semi)group acting on the symmetric section\/}.
\vskip 11pt

%% 5
\item
If the geometric dimension is $2k$, then we have a tower of fibered or shifted real conjugacy class representatives
\[
G^{(2k)}((F_{\o v_1}\times \rit)\times(F_{v_1}\times \rit )) \subset\dots\subset
G^{(2k)}((F_{\o v_{j,m_j}}\times \rit)\times(F_{v_{j,m_j}}\times \rit))
\subset\dots
\]
 generated under the (bi)action of
 $\GL_{2k}(\Fvt)$ or, equivalently, a tower of class representatives
 \begin{multline*}
G^{(1)}_{(Q)}((F^{(2k)}_{\o v^1}\times \rit^{(2k)})\times(F^{(2k)}_{v^1}\times \rit^{(2k)} )) \\
\subset\dots\subset
G^{(j,m_j)}_{(Q)}((F^{(2k)}_{\o v^1}\times \rit^{(2k)})\times(F^{(2k)}_{v^1}\times \rit^{( 2k )}))
\subset\dots
 \end{multline*}
of the fibered or shifted bisemivariety
$G^{(\o v\times v)}_{(Q)}((F^{(2k)}_{\o v^1}\times \rit^{(2k)})\times(F^{(2k)}_{v^1}\times \rit^{(2k)} )) $ generated from
$\GL^{(Q)}((F^{(2k)}_{\o v^1}\times \rit^{(2k)})\times(F^{(2k)}_{v^1}\times \rit^{(2k)} )) $.

Referring to 4), we see that the $2k$-dimensional representation of
$(\Aut_k(\phi _R(F_{\o v}\times \rit )) \times \Aut_k(\phi _L(F_v\times\rit) )) $
given by
$\Rep^{(2k)}(\Aut_k(\phi _R(F_{\o v}\times \rit )) \times \Aut_k(\phi _L(F_v\times\rit)))$ decomposes into:
\begin{multline*}
\Rep^{(2k)} ( \Aut_k (\phi _R ( F_{\o v}\times \rit ) )\times \Aut_k ( \phi _L ( F_v\times\rit )))\\
=\Rep^{(2k)} ( \Aut_k ( \phi _R ( F_{\o v}))\times \Aut_k ( \phi _L ( F_v )))\\
\times
\Rep^{(2k)} ( \Aut_k ( \phi _R ( \rit )) \times \Aut_k ( \phi _L ( \rit )))
\end{multline*}
which generates (and is isomorphic to) the fibered or shifted bisemivariety
$G^{(2k)}(\Fvtrit)$, rewritten according to
$G^{(\o v\times v)}_{(Q)}((F^{(2k)}_{\o v^1}\times \rit^{(2k)})\times(F^{(2k)}_{v^1}\times \rit^{(2k)})$.

$\Rep^{(2k)}(\Aut_k(\phi _R( \rit )) \times \Aut_k(\phi _L(\rit)))$ is thus the $2k$-dimensional representation of the differential bilinear Galois semigroup
$(\Aut_k(\phi _R( \rit))\times \Aut_k(\phi _L(\rit)))$ and is isomorphic to the shifting bisemivariety
$G^{(2k)}(\phi _R( \rit)\times\phi _L(\rit))$ or
$G^{(\o v\times v)}_{(Q)}(\phi _R( \rit^{(2k)})\times\phi _L(\rit^{(2k)}))$.

So, we have {\bbf a tower of $2k$-dimensional representations of the differential bilinear Galois subsemigroups\/}
\begin{multline*}
\Rep^{(2k)}(\Aut_k(\phi _R( \rit)_{|F_{\o v_1}})\times \Aut_k(\phi _L(\rit)_{|F_{v_1}}))\\ \subset\dots\subset
\Rep^{(2k)}(\Aut_k(\phi _R( \rit)_{|F_{\o v_{j,m_j}}})\times \Aut_k(\phi _L(\rit)_{|F_{v_{j,m_j}}})) \\ \subset\dots\subset
\Rep^{(2k)}(\Aut_k(\phi _R( \rit)_{|F_{\o v_{r,m_r}}})\times \Aut_k(\phi _L(\rit)_{|F_{v_{r,m_r}}}))
\end{multline*}
which are respectively in one-to-one correspondence with the bilinear semigroups:
\begin{multline*}
\GL_1(\phi _R( \rit^{(2k)})\times \phi _L(\rit^{(2k)}))=\GL_{2k}(\phi _R( \rit)\times\phi _L(\rit))_{|F_{\o v_{1}}\times F_{v_{1}}} \\ 
\subset \dots \subset
\GL_j(\phi _R( \rit^{(2k)})\times \phi _L(\rit^{(2k)}))=\GL_{2k}(\phi _R( \rit)\times\phi _L(\rit))_{|F_{\o v_{j,m_j}}\times F_{v_{j,m_j}}} \\
\subset \dots \subset
\GL_r(\phi _R( \rit^{(2k)})\times \phi _L(\rit^{(2k)}))=\GL_{2k}(\phi _R( \rit)\times\phi _L(\rit))_{|F_{\o v_{r,m_r}}\times F_{v_{r,m_r}}} \;.
\end{multline*}
So, we get the thesis by considering the upper algebraic dimension ``$r$'':
\[
\Rep^{(2k)}(\Aut_k(\phi _R( \rit)_{|F_{\o v_r}})\times \Aut_k(\phi _L(\rit)_{|F_{v_r}}))=
\GL_r(\phi _R( \rit^{(2k)})\times \phi _L(\rit^{(2k)}))
\] 

where
$\phi _R( \rit^{(2k)})\simeq \rit^{(2k)}$ and $ \phi _L(\rit^{(2k)})\simeq \rit^{(2k)}$ are generally constant functions.

The unitary case referring to the $2k$-dimensional representation of the bilinear inertia semigroup
$\Int_k(\phi _R(\rit )) \times \Int_k(\phi _L(\rit )) $ can be reached similarly, i.e.
\begin{align}
\Rep^{(2k)}(\Int_k(\phi _R( \rit)_{|F_{\o v^1_r}})&\times \Int_k(\phi _L(\rit)_{|F_{v^1_r}}))\notag\\
&=
P_r(\phi _R( \rit^{(2k)})\times \phi _L(\rit^{(2k)}))\notag\\
&=
\GL_r(\phi _R( \rit^{(2k)})\times \phi _L(\rit^{(2k)}))_{|F_{\o v^1}\times F_{v^1}}
\tag*{\eop}
\end{align}
\Ee
\vskip 11pt

\subsection{Corollary}

{\em
The complex bilinear representation of differential Galois (semi)groups can be found similarly as it was done for the real case.

Let
$F_\omega =\{F_{\omega _1},\dots,F_{\omega _{j,m_j}},\dots,F_{\omega _{r,m_r}}\}$
\resp{$F_{\o\omega} =\{F_{\o \omega _1},\dots,F_{\o \omega _{j,m_j}},\dots,F_{\o \omega _{r,m_r}}\}$} be the set of complex transcendental extensions (or infinite complex archimedean completions) covered by its real equivalent $F_v$ \resp{$F_{\o v}$}.

Then, {\bbf the complex bilinear semigroup\/} of matrices
$\GL_r(\cit^k\times\cit^k)$, of algebraic order ``$r$'', referring to the action of the bioperator
$(D^{2k}_R\otimes D^{2k}_L)$ on the complex bisemisheaf
$(\wh M^{(2k)}_{\o\omega _R}\otimes\wh M^{(2k)}_{\omega _L})$ on the abstract complex bisemivariety
$G^{(2k)}(F_{\o\omega }\times F_\omega )\equiv \REPSP (\GL_k(F_{\o\omega }\times F_\omega ))$,
{\bbf corresponds to the $k$-dimensional complex bilinear representation of the product, right by left, of ``differential'' Galois (or global Weil) semigroups
$(\Aut_k(\psi _R(\cit))\times \Aut_k(\psi _L(\cit) )) $ shifting the product, right by left, of automorphism semigroups
$(\Aut_k(\psi _R(F_{\o\omega }))\times \Aut_k(\psi _L(F_\omega ) )) $ of cofunctions 
$\psi _R(F_{\o\omega })$ and $\psi _L(F_\omega )$ by 
$(\Aut_k(\psi _R(F_{\o\omega }\times\cit))\times \Aut_k(\psi _L(F_\omega )\times\cit )) $\/}.

So, we have that:
\[
\GL_r(\phi _R(\cit^k)\times \phi _L(\cit^k))=
\Rep^{(2k)}    (\Aut_k(\psi _R(\cit)_{|F_{\o\omega_r}}\times \Aut_k(\psi _L(\cit)_{|F_{\omega _R}})))\;.\]
And, in the unitary case, we have:
\[
P_r(\psi _R(\cit^k)\times \psi _L(\cit^k))=
\Rep^{(2k)}    (\Int_k(\psi _R(\cit)_{|F_{\o\omega^1_r}} )\times \Int_k(\psi _L(\cit)_{|F_{\omega^1 _R}}))    \]
where:
\Bi
\item $P_r(\dots)$ is the bilinear parabolic subsemigroup;
\item $\Rep^{(2k)}(\dots)$ is the $k$-dimensional complex representation of $(\dots)$.
\Ei}
\vskip 11pt

\subsection{Corollary}

{\em
Let $O_r(\rit)$ denote the orthogonal group of algebraic order $r$ with entries in the reals $\rit$ and let $U_r(\cit)$ denote the unitary group of algebraic order $r$ in the complexes $\cit$.

Then, the orthogonal bilinear semigroup $O_r(\rit^{(2k)}\times\rit^{(2k)})$ corresponds to the real parabolic bilinear semigroup
 $P_r(\rit^{(2k)}\times\rit^{(2k)})$ and the unitary bilinear semigroup $U_r(\cit^k\times\cit^k)$ corresponds to the complex parabolic bilinear semigroup
 $P_r(\cit^k\times\cit^k)$.
 }
 \vskip 11pt
 
 \bpr
 This results from the definition of a bilinear semigroup recalled in section~1.5 and from proposition~2.4 and corollary~2.5.
 
 We have more particularly that:
 \begin{align}
 O_r(\rit^{(2k)}\times\rit^{(2k)})=\Rep^{(2k)}(\Int_k(\phi _R(\rit))_{|F_{\o v^1_r}}\times\Int_k(\phi _L(\rit))_{|F_{v^1_r}})\notag\\
 \noalign{and that}
 U_r(\cit^{k}\times\cit^{k})=\Rep^{(2k)}(\Int_k(\psi _R(\cit))_{|F_{\o \omega ^1_r}}\times\Int_k(\psi _L(\cit))_{|F_{\omega ^1_r}})\;.\tag*{\eop}
 \end{align}
 \vskip 11pt
 
 \subsection{Corollary}
 
 {\em
 In the one-dimensional geometric case, i.e. when $2k=1$, we have that
 \Bena
 \item {\bbf the orthogonal bilinear semigroup of algebraic order $r$
 \[ O_r(\rit\times\rit)=P_r(\rit\times\rit)\]
 corresponds to the product, right by left, of differential inertia Galois semigroups
 $(\Int_k(\phi _R(\rit))_{|F_{\o v^1_r}}\times\Int_k(\phi _L(\rit))_{|F_{v^1_r}})$\/} shifting the product, right by left,\linebreak
 $(\Int_k(\phi _R(F_{\o v^1_r}))\times\Int_k(\phi _L(F_{v^1_r}) )) $ of internal automorphism semigroups of cofunctions 
 $\phi _R(F_{\o v^1_r})$ and functions $\phi _L(F_{v^1_r})$ respectively on the unitary transcendental lower and upper real extensions 
 $F_{\o v^1_r}$ and $F_{v^1_r}$.
 
 \item {\bbf the unitary bilinear semigroup of algebraic order $r$
  \[ U_r(\cit\times\cit)=P_r(\cit\times\cit)\]
corresponds to the product, right by left, of differential inertia Galois semigroups
 $(\Int_k(\psi _R(\cit))_{|F_{\o \omega ^1_r}}\times\Int_k(\psi _L(\cit))_{|F_{\omega ^1_r}})$\/} shifting the product, right by left,
 $(\Int_k(\psi _R(F_{\o \omega ^1_r}))\times\Int_k(\psi _L(F_{\omega ^1_r}) )) $ of internal automorphism semigroups of cofunctions 
 $\psi _R(F_{\o \omega ^1_r})$ and functions $\psi _L(F_{\omega ^1_r})$ respectively on the unitary complex transcendental  upper  extensions 
 $F_{\o \omega ^1_r}$ and $F_{\omega ^1_r}$.
 \Ee
 }
 \vskip 11pt
 
 \bpr
 \Bena
\item $(\Int_k(\phi _R(F_{\o v^1_r}))\times\Int_k(\phi _L(F_{v^1_r}) ))$ corresponds to the bilinear semigroup of internal automorphisms of bifunctions on ``unitary'' biquanta in a bisection
$\phi _R(F_{\o v_r})\times \phi _L(F_{v^1_r})$ at ``$r$'' biquanta
$(F_{\o v_r}\times F_{v_r})$ while
$(\Int_k(\phi _R(F_{\o v^1_r}\times\rit))\times\Int_k(\phi _L(F_{v^1_r}\times\rit) ))$ corresponds to the bilinear semigroup of shifted internal automorphisms (under the action of a bioperator
$(D_R\otimes D_L)$) of bifunctions on shifted ``unitary'' biquanta
$((F_{\o v^1_r}\times \rit)\times (F_{v^1_r}\times \rit))$ in a bisection
$\phi _R(F_{\o v_r}\times\rit)\times \phi _L(F_{v_r}\times\rit)$ at ``$r$'' shifted biquanta
$(F_{\o v_r}\times\rit)\times (F_{v_r}\times\rit)$.

\item The complex unitary case referring to $U_r(\cit\times\cit)$ can be handled similarly as the real unitary case referring to $O_r(\rit\times\rit)$ by taking into account that a complex bisection
$\psi _R(F_{\o \omega _r})\times\psi _L(F_{\omega _r})$ is covered by real bisections 
$\{\phi _R(F_{\o v _{r,m_r}})\times\phi _L(F_{v _{r,m_r}})\}_{m_r}$ as developed in section~1.9.\epr
\Ee
\vskip 11pt

\subsection[Bilinear Hilbert semispaces and Von Neumann bisemialgebras]{Bilinear Hilbert semispaces and Von Neumann bisemi-\linebreak algebras}

Let $(\wh M^{(2n)}_{v_R}\otimes \wh M^{(2n)}_{v_L})$ be the bisemisheaf of differentiable bifunctions on the abstract bisemivariety
$G^{(2n)}(\Fvt)$.

The set
\[ \{\phi ^{(2n)}_{j,m_j}(g^{(2n)}_{v\RL}[j,m_j])\}_{j,m_j}=
  \{\phi ^{(2n)}_{j,m_j}(g^{(2n)}_{v_R}[j,m_j])\otimes
  \phi ^{(2n)}_{j,m_j}(g^{(2n)}_{v_L}[j,m_j])\}_{j,m_j}\]
  of differentiable bifunctions, i.e. bisections of
  $(\wh M^{(2n)}_{v_R}\otimes \wh M^{(2n)}_{v_L})$,
 \\
  forms an increasing filtration with respect to the algebraic dimension ``$j$'':
\[
\phi ^{(2n)}_{1}(g^{(2n)}_{v\RL}[1])\subset \dots \subset
\phi ^{(2n)}_{j,m_j}(g^{(2n)}_{v\RL}[j,m_j]) \subset \dots \subset
\phi ^{(2n)}_{r,m_r}(g^{(2n)}_{v\RL}[r,m_r]) \;, \qquad 
j\le r\le \infty \;, \]
on the corresponding filtration of conjugacy class representatives of the abstract bisemivariety $G^{(2n)}(\Fvt)$:
\[ g^{(2n)}_{v\RL}[1] \subset \dots \subset
  g^{(2n)}_{v\RL}[j,m_j] \subset \dots \subset
  g^{(2n)}_{v\RL}[r,m_r] \;.\]
  \vskip 11pt
  
 This bisemisheaf $(\wh M^{(2n)}_{v_R}\otimes \wh M^{(2n)}_{v_L})$ is transformed into
 {\bbf an extended internal left bilinear Hilbert semispace
 $H^+_{(\wh M^{(2n)}_{v_R}\otimes \wh M^{(2n)}_{v_L})}$\/} by taking into account
 \Bena
 \item a map \cite{Pie6}:
 \[ B_L\circ p_L: \qquad \wh M^{(2n)}_{v_R}\otimes \wh M^{(2n)}_{v_L}
 \quad \To \quad \wh M^{(2n)}_{v_{L_R}}\otimes \wh M^{(2n)}_{v_L})
 \equiv H^+_{(\wh M^{(2n)}_{v_R}\otimes \wh M^{(2n)}_{v_L})}\]
 where:
 \Bi
 \item $p_L$ is a projective linear map projecting the right semisheaf
 $\wh M^{(2n)}_{v_R}$ onto the left semisheaf $\wh M^{(2n)}_{v_L}$;
 \item $B_L$ is a bijective linear isometric map from the projected right semisheaf\linebreak
 $\wh M^{(2n)}_{v_{L(P)_R}}$ to $\wh M^{(2n)}_{v_{L_R}}$ mapping each covariant element of 
 $\wh M^{(2n)}_{v_{L(P)_R}}$ into a contravariant element.
 \Ei
 
 \item an internal bilinear form defined from
 $H^+_{(\wh M^{(2n)}_{v_R}\otimes \wh M^{(2n)}_{v_L})}$ into $\cit$ for every bisection
 $\phi ^{(2n)}_{j,m_j}(g^{(2n)}_{v\RL}[j,m_j])$ by:
 \[
 (\phi ^{(2n)}_{j,m_j}(g^{(2n)}_{v_{L_R}}[j,m_j],
    \phi ^{(2n)}_{j,m_j}(g^{(2n)}_{v_{L}}[j,m_j] )))\To \cit\;.\]
 \Ee
 \vskip 11pt
 
 This bilinear Hilbert semispace
 $H^+_{(\wh M^{(2n)}_{v_R}\otimes \wh M^{(2n)}_{v_L})}$, noted in condensed form
 $H^+_{\wh M^{(2n)}_{v\RL}}$, is a natural representation (bisemi)space for the bialgebra of elliptic bioperators
 $(D^{2k}_R\otimes D^{2k}_L)$ as noticed in \cite{Pie6}.

 {\bbf A bisemialgebra of von Neumann $\MM\RL(H^+_{\wh M^{(2n)}_{v\RL}})$ in 
$H^+_{\wh M^{(2n)}_{v\RL}}$\/} is an involutive subbisemialgebra of the bisemialgebra
$(\Ls^B_R\otimes \Ls^B_L)(H^+_{\wh M^{(2n)}_{v\RL}})$ of bounded bioperators having a closed norm topology.
\vskip 11pt

Due to the structure of the bisemisheaf
$(\wh M^{(2n)}_{v_R}\otimes \wh M^{(2n)}_{v_L})$, the bilinear Hilbert semispace
$H^+_{\wh M^{(2n)}_{v\RL}}$ is ``solvable'' in the sense that we have a tower of embedded bilinear Hilbert subsemispaces
\[ H^+_{\wh M^{(2n)}_{v\RL}}(1) \subset \dots \subset
 H^+_{\wh M^{(2n)}_{v\RL}}(j) \subset \dots \subset
 H^+_{\wh M^{(2n)}_{v\RL}}(r) \]
where $H^+_{\wh M^{(2n)}_{v\RL}}$ is given by the set
$ \{  \phi ^{(2n)}_{j,m_j}(g^{(2n)}_{v_{L_R}}[j,m_j])\otimes
  \phi ^{(2n)}_{j,m_j}(g^{(2n)}_{v_{L}}[j,m_j] )\}_{j,m_j}$ of bisections.
  
  But, we can also construct a tower of direct sums of embedded extended bilinear Hilbert subsemispaces:
  \[ H^+_{\wh M^{(2n)}_{v\RL}}\{1\} \subset \dots \subset
 H^+_{\wh M^{(2n)}_{v\RL}}\{j\} \subset \dots \subset
 H^+_{\wh M^{(2n)}_{v\RL}}\{r\} \]
 where $ H^+_{\wh M^{(2n)}_{v\RL}}\{j\}$ is defined by: 
 \[
 H^+_{\wh M^{(2n)}_{v\RL}}\{j\}=\bigoplus_{\nu =1}^j   H^+_{\wh M^{(2n)}_{v\RL}}(\nu )\;.\]
 \vskip 11pt
 
 \subsection{Random bioperators}
 
 Let $(D^{2k}_R\otimes D^{2k}_L)$ be the differential (elliptic) bioperator acting on the set
 $\{\phi ( G^{(2n)}(F_{\o v_{j,m_j}}\times F_{v_{j,m_j}}))\}_{j,m_j}$  (also written
 $\{\phi^{(2n)}_{j,m_j}(g^{(2n)}_{v\RL}[j,m_j])\}_{j,m_j}$) of differentiable bifunctions or bisections of the bisemisheaf
 $(\wh M^{2n}_{v_R}\otimes \wh M^{2n}_{v_L})$ according to:
\begin{multline*}
D^{2k}_R\otimes D^{2k}_L: 
\begin{aligned}[t]
&(\wh M^{(2n)}_{v_R}\otimes \wh M^{(2n)}_{v_L})
\To (\wh M^{(2n)}_{v_R}[2k]\otimes \wh M^{(2n)}_{v_L}[2k])\\
&\{\phi_{j,m_j} (G^{(2n)}(F_{\o v_{j,m_j}}\times F_{v_{j,m_j}}))\}_{j,m_j}
\end{aligned}
\\
\To
\{\phi_{j,m_j} (G^{(2n)}_{[2k]}((F_{\o v_{j,m_j}}\times \rit)\times (F_{v_{j,m_j}}\times\rit)))  \}_{j,m_j}
\end{multline*}
where
$\phi_{j,m_j} (G^{(2n)}_{[2k]}((F_{\o v_{j,m_j}}\times \rit)\times (F_{v_{j,m_j}}\times\rit) )) $ is the bifunction on the $(j,m_j)$-th conjugacy class representative
$G^{(2n)}_{[2k]}((F_{\o v_{j,m_j}}\times \rit)\times (F_{v_{j,m_j}}\times\rit))$ fibered or shifted in $2k$ bilinear geometric dimensions.

The bioperator
$(D^{2k}_R\otimes D^{2k}_L)$ is a random bioperator in the sense that, for every bifunction
\[
\phi_{j,m_j} (G^{(2n)}_{[2k]}(F_{\o v_{j,m_j}}\times F_{v_{j,m_j}} )) =
\phi_{j,m_j} (G^{(2n)}_{[2k]}(F_{\o v_{j,m_j}} )) \times
\phi_{j,m_j} (G^{(2n)}_{[2k]}(F_{v_{j,m_j}}))
\] belonging to the bilinear Hilbert semispace
$H^+_{\wh M^{(2n)}_{v\RL}}$ the bilinear form
\[
(D^{2k}_R(\phi_{j,m_j} (G^{(2n)}_{[2k]}(F_{\o v_{j,m_j}} ))) ,
(D^{2k}_L(\phi_{j,m_j} (G^{(2n)}_{[2k]}(F_{v_{j,m_j}})) ))) \]
is measurable.

{\bbf The random bioperator\/}
$(D^{2k}_R\otimes D^{2k}_L)$ acting on
$H^+_{\wh M^{(2n)}_{v\RL}}$ is a set
$\{ (D^{2k}_R(j,m_j)\otimes D^{2k}_L(j,m_j )) \}_{j,m_j}$ of bioperators acting on the bisemisheaf
$\wh M^{(2n)}_{v\RL}$.
\vskip 11pt

\subsection{Towers of embedded von Neumann subbisemialgebras}

Referring to the tower of embedded bilinear Hilbert subsemispaces associated with the bisemisheaf
$\wh M^{(2n)}_{v\RL}$ and to the definition of a bisemialgebra of von Neumann\linebreak
$\MM\RL(H^+_{\wh M^{(2n)}_{v\RL}})$ given in section~2.8, it appears that there exists a {\bbf tower of embedded von Neumann subbisemialgebras\/}:
\[
\MM\RL(H^+_{\wh M^{(2n)}_{v\RL}}(1)) \subset \dots \subset
\MM\RL(H^+_{\wh M^{(2n)}_{v\RL}}(j)) \subset \dots \subset
\MM\RL(H^+_{\wh M^{(2n)}_{v\RL}}(r)) \]
according to the algebraic dimensions $1\le j\le r$, as well as a 
{\bbf tower of sums of embedded von Neumann subbisemialgebras\/}:
\[
\MM\RL(H^+_{\wh M^{(2n)}_{v\RL}}\{1\}) \subset \dots \subset
\MM\RL(H^+_{\wh M^{(2n)}_{v\RL}}\{j\}) \subset \dots \subset
\MM\RL(H^+_{\wh M^{(2n)}_{v\RL}}\{r\}) \]
where 
\[
\MM\RL(H^+_{\wh M^{(2n)}_{v\RL}}\{j\})=
\bigoplus_{\nu =1}^j
\MM\RL(H^+_{\wh M^{(2n)}_{v\RL}}(\nu ))\;.\]

The bisemisheaf $(\wh M^{2n}_{v_R}\otimes \wh M^{2n}_{v_L})$ gives rise to the extended internal left bilinear Hilbert semispace
$H^+_{\wh M^{2n}_{v_R}\otimes \wh M^{2n}_{v_L}}$ according to section~2.8.

Similarly, the diagonal bisemisheaf
$(\wh M^{2n}_{v_R}\otimes_D \wh M^{2n}_{v_L})$, whose bisections are diagonal bisections
 \[
\phi_{j,m_j} ( G^{(2n)}(F_{\o v_{j,m_j}}\times_D F_{v_{j,m_j}}))=
 \phi_{j,m_j} ( G^{(2n)}(F_{\o v_{j,m_j}}))\times_D  \phi_{j,m_j} ( G^{(2n)}( F_{v_{j,m_j}}))\]
 characterized by a diagonal bilinear basis (the offdiagonal bilinear basis elements being null) \cite{Pie6}, gives rise to the {\bbf diagonal internal left bilinear Hilbert semispace\/}  $\Hs^+_{\wh M^{2n}_{v_R}\otimes_D \wh M^{2n}_{v_L}}$
 by taking into account a $(B_L\circ p_L)$ map and the existence of an internal diagonal bilinear form, i.e. a scalar product, as in section~2.8.
 \vskip 11pt
 
 \subsection{Proposition}
 
 {\em
 Let $\MM\RL(H^+_{\wh M^{(2n)}_{v\RL}} )$ be the bisemialgebra of von Neumann on the extended internal left bilinear Hilbert semispace
 $H^+_{\wh M^{2n}_{v_R}\otimes \wh M^{2n}_{v_L}}$ and let
 $\MM\RL(\Hs^+_{\wh M^{2n}_{v_R}\otimes_D \wh M^{2n}_{v_L}})$ be the von Neumann bisemialgebra on the diagonal internal left bilinear Hilbert semispace
 $\Hs^+_{\wh M^{2n}_{v_R}\otimes_D \wh M^{2n}_{v_L}}$.
 
 Then, {\bbf the discrete spectrum $\Sigma (D^{2k}_R\otimes D^{2k}_L)$ of the bioperator $(D^{2k}_R\otimes D^{2k}_L)$ is obtained by the composition of morphisms\/}:
 \begin{multline*}
 i_{\{ j \} ^D\RL}\circ i_{\{ j \} \RL} : 
 \begin{array}[t]{ccc}
 \MM\RL(H^+_{\wh M^{(2n)}_{v\RL}})& \To & \L[ \MM\RL(\Hs_{\wh M^{(2n)}_{v_R}[2k]\otimes_D\Hs_{\wh M^{(2n)}_{v_L}[2k]}})\R]_j\\
 (D^{2k}_R\otimes D^{2k}_L) & \To & \Sigma  (D^{2k}_R\otimes D^{2k}_L) 
 \end{array}\end{multline*}
 where $i_{\{j\}\RL}$ and $ i_{\{j\}^D\RL}$ are given by:
 \begin{align*}
 i_{\{j\}\RL}: \quad &\MM\RL(H^+_{\wh M^{(2n)}_{v\RL}})& \To & \L[ \MM\RL(H^+_{\wh M^{(2n)}_{v\RL[2k]}}\{j\})\R]_j \\
 i_{\{j\}^D\RL}: \quad &\L[ \MM\RL(H^+_{\wh M^{(2n)}_{v\RL[2k]}}\{j\} )\R]_j & \To 
& \L[ \MM\RL(\Hs^+_{\wh M^{(2n)}_{v_R[2k]}\otimes_D
  \wh M^{(2n)}_{v_L[2k]}}\{j\} )\R]_j\;.
 \end{align*}
 $\L[ \MM\RL(\Hs^+_{\wh M^{(2n)}_{v_R[2k]}\otimes_D
 \wh M^{(2n)}_{v_L[2k]}}\{j\})\R]_j$ is the increasing tower, over the running algebraic index ``$j$'', of sums of von Neumann subbisemialgebras:
 \[
 \MM\RL(\Hs^+_{\wh M^{(2n)}_{v_R[2k]}\otimes_D
 \wh M^{(2n)}_{v_L[2k]}}\{j\})=
\bigoplus_{\nu =1}^j \MM\RL(\Hs^+_{(\wh M^{(2n)}_{v_R[2k]}\otimes\wh M^{(2n)}_{v_L[2k]})}(\nu )\]
over respectively the tower:
 \[
 \Hs^+_{\wh M^{(2n)}_{v_R[2k]} \otimes_D
 \wh M^{(2n)}_{v_L[2k]}}\{1\}\subset \dots \subset 
 \Hs^+_{\wh M^{(2n)}_{v_R[2k]} \otimes_D
 \wh M^{(2n)}_{v_L[2k]}}\{j\}  \subset \dots \subset 
  \Hs^+_{\wh M^{(2n)}_{v_R[2k]} \otimes_D
 \wh M^{(2n)}_{v_L[2k]}}\{r\}\]
 of sums of diagonal internal left bilinear Hilbert subsemispaces
 \[
\Hs^+_{\wh M^{(2n)}_{v_R[2k]} \otimes_D
 \wh M^{(2n)}_{v_L[2k]}}\{j\} )=
 \bigoplus_{\nu =1}^j
 (\Hs^+_{(\wh M^{(2n)}_{v_R[2k]}\otimes_D\wh M^{(2n)}_{v_L[2k]})}(\nu  )) \]
 shifted in $(2k)$ bilinear geometric dimensions.
 }
 \vskip 11pt
 
 \bpr
 The morphism $ i_{\{j\}\RL}$
 \[
  i_{\{j\}\RL}: \qquad 
  \MM\RL(H^+_{\wh M^{(2n)}_{v\RL}})\To  \L[ \MM\RL(H^+_{\wh M^{(2n)}_{v\RL[2k]}}\{j\} )\R]_j\]
  is in fact implicit depending on the decomposition of the bisemisheaf
 $(\wh M^{(2n)}_{v_R} \otimes \wh M^{(2n)}_{v_L})$ into bisections on the conjugacy class representatives
 $\{g^{(2n)}_{v\RL} [j,m_j]\}_{j,m_j}$ of the abstract bisemivariety
 $G^{(2n)}(\Fvt)$.
 
 $ \MM\RL(\Hs^+_{\wh M^{(2n)}_{v\RL[2k]}}\{j\} )$ is the subbisemialgebra of von Neumann on the {\bf shifted\/} extended internal left bilinear Hilbert subsemispace
 $ \Hs^+_{\wh M^{(2n)}_{v\RL[2k]}}\{j\}$.
 
 The morphism
 \[
  i_{\{j\}^D\RL}: \quad \L[ \MM\RL(H^+_{\wh M^{(2n)}_{v\RL[2k]}}\{j\})\R]_j  \To 
  \L[ \MM\RL(\Hs^+_{\wh M^{(2n)}_{v_R[2k]}\otimes_D
\wh M^{(2n)}_{v_L[2k]}}\{j\} )\R]_j \]
 sends each von Neumann subbisemialgebra
 \[
\MM\RL(H^+_{\wh M^{(2n)}_{v\RL[2k]}} )\{j\} =\bigoplus_{\nu =1}^j
\MM\RL(H^+_{\wh M^{(2n)}_{v\RL[2k]}}(\nu ))\]
on the shifted extended internal left bilinear Hilbert subsemispace
\[
H^+_{\wh M^{(2n)}_{v\RL[2k]}}\{j\} =\bigoplus_{\nu =1}^j
H^+_{\wh M^{(2n)}_{v\RL[2k]}}(\nu )\]
onto the corresponding von Neumann diagonal subbisemialgebra
\[
 \MM\RL(\Hs^+_{\wh M^{(2n)}_{v_R[2k]} \otimes_D \wh M^{(2n)}_{v_L[2k]}}\{j\})
 =\bigoplus_{\nu =1}^j
\MM\RL(\Hs^+_{\wh M^{(2n)}_{v_R[2k]} \otimes_D \wh M^{(2n)}_{v_L[2k]}}\{j\})
\]
on the diagonal internal left bilinear Hilbert subsemispace:
\[
 \Hs^+_{\wh M^{(2n)}_{v_R[2k]} \otimes_D \wh M^{(2n)}_{v_L[2k]}}\{j\}
 =\bigoplus_{\nu =1}^j
\Hs^+_{\wh M^{(2n)}_{v_R[2k]} \otimes_D \wh M^{(2n)}_{v_L[2k]}}\{\nu \}
\]
shifted in $(2k)$ bilinear geometric dimensions.\epr
\vskip 11pt
 
 \subsection{Cuspidal representations of the global program of Langlands}
 
 The differential bioperator 
$(D^{2k}_R\otimes D^{2k}_L)$ maps the bisemisheaf 
$(\wh M_{v_R}^{(2n)}\otimes \wh M_{v_L}^{(2n)} )$ into the corresponding perverse bisemisheaf
$(\wh M_{v_R}^{(2n)}[2k]\otimes \wh M_{v_L}^{(2n)}[2k])$ according to:
\[
(D^{2k}_R\otimes D^{2k}_L): \qquad
\wh M_{v_R}^{(2n)}\otimes \wh M_{v_L}^{(2n)} \quad\To\quad
\wh M_{v_R}^{(2n)}[2k]\otimes \wh M_{v_L}^{(2n)}[2k]
\]
in such a way that
$\wh M_{v_R}^{(2n)}[2k]\otimes \wh M_{v_L}^{(2n)}[2k]$ decomposes into a tower of fibered or shifted bisections or bifunctions
$\{\phi _{j,m_j}(G^{(2n)}_{[2k]}(F_{\o v_{j,m_j}})) \times
\phi _{j,m_j}(G^{(2n)}_{[2k]}(F_{v_{j,m_j}}))\}$ (see section~2.10).
\vskip 11pt

On the other hand, referring to the global program of Langlands \cite{Pie2}, {\bbf there is a one-to-one correspondence between the bisemisheaf
$\wh M_{v_R}^{(2n)}\otimes \wh M_{v_L}^{(2n)}$ over the abstract bisemivariety
$G^{(2n)}(\Fvt)$ and its cuspidal counterpart\/}
$(\wh M_{v^T_R}^{(2n)}\otimes \wh M_{v^T_L}^{(2n)})$
on the toroidal abstract bisemivariety
$G^{(2n)}(F^T_{\o v}\times F^T_v)$ over the sets
$F^T_v=\{F^T_{v_1},\dots,F^T_{v_{j,m_j}},\dots,\linebreak F^T_{v_{,m_r}}\}$ and
$F^T_{\o v}=\{F^T_{\o v_1},\dots,F^T_{\o v_{j,m_j}},\dots,F^T_{\o v_{r,m_r}}\}$
of toroidal real archimedean completions or transcendental extensions.

The toroidal bisemisheaf has for bisections the bifunctions
\[
\phi _{j,m_j}(G^{(2n)}(F^T_{\o v_{j,m_j}}) )   \times
\phi _{j,m_j}(G^{(2n)}(F^T_{v_{j,m_j}})  ) 
=\L( \lambda (2n,j,m_j)\ e^{-2\pi ijx}\R)\times
\L( \lambda (2n,j,m_j)\ e^{+2\pi ijx}\R)
\]
where:
\Bi
\item $\vec x=\sum\limits_{c=1}^{2n}x_c\ \vec e_c$, $x\in\rit^{2n}$;

\item $\lambda ^2(2n,j,m_j)=\prod\limits_{c=1}^{2n}\lambda ^2_c(2n,j,m_j)$ is a product of eigenbivalues
$\lambda ^2_c(2n,j,m_j)$ of the Hecke bioperator
$(T_R(2n;r) \otimes T_L(2n;r) )$ whose representation is $\GL_{2n}(\Os_{F^T_{\o v}}\times \Os_{F^T_{v}})$ referring to section~1.5.
\Ei

Let 
$\Gamma _{\wh M^{(2n)}_{v^T\RL}}= \{
\phi _{j,m_j}(G^{(2n)}(F^T_{\o v_{j,m_j}}))\times
\phi _{j,m_j}(G^{(2n)}(F^T_{v_{j,m_j}}))\}_{j,m_j}$ denote the set of increasing bisections of the bisemisheaf
$\widehat M_{v^T_R}^{(2n)}\otimes \widehat M_{v^T_L}^{(2n)}$.
\vskip 11pt

Then, {\bbf
a global elliptic $(\Gamma _{\widehat M_{v^T\RL}^{(2n)}})$-bisemimodule\/}
$\phi ^{(2n)}\RL(x)$, referring to the bihomomorphism
\[
\phi ^{(2n)}\RL(x) : \qquad
\Gamma _{\wh M^{(2n)}_{v^T\RL}} \quad \To \quad \End (
\Gamma _{\wh M^{(2n)}_{v^T\RL}} )\;, \]
is given by:
\[
\phi ^{(2n)}\RL(x)=
\sum_j\sum_{m_j} \L (\lambda (2n,j,m_j)\ e^{-2\pi ijx}\R)\times
\sum_j\sum_{m_j} \L (\lambda (2n,j,m_j)\ e^{+2\pi ijx}\R)
\]
in such a way that
$\phi ^{(2n)}\RL(x)$
{\bbf constitutes a real cuspidal representation of bilinear geometric dimension $2n$, of the product, right by left, of Weil global semigroups
$\Gal(\wt F_{\o v}/k)\times\Gal(\wt F_{v}/k) $\/} according to the global program of Langlands.

Remark that $\phi ^{(2n)}\RL(x)$ covers the corresponding ``complex'' cuspidal representation \cite{Pie2}.
\vskip 11pt

\subsection{Proposition}

{\em
{\bbf The global elliptic  $(\Gamma _{\widehat M_{v^T\RL}^{(2n)}})$-bisemimodule\/}
$\phi ^{(2n)}\RL(x)$ is the functional representation space
$\FREPSP(\GL_{2n}(\FTvt))$ of the bilinear semigroup 
$\GL_{2n}(\FTvt)$, over the product, right by left, of toroidal real archimedean completions
$F^T_{\o v}$ and
$F^T_{v}$, under the (bi)action of the (bi)monomorphisms:
\[
\sigma _{v^T_R}\times\sigma _{v^T_L}: \qquad
\Aut_k(F^T_{\o v})\times\Aut_k(F^T_{v} )   \quad \To \quad
G^{(2n)}(F^T_{\o v}\times F^T_{v})\;,\]
where $\Aut_k(F^T_{\o v}\times\Aut_k(F^T_{v} )) $ is the bilinear semigroup of automorphisms of toroidal transcendental extensions associated with a $1D$-bisemilattice of transcendental biquanta.
}

\bpr
This results from the bimonomorphism
\[
\sigma _{v_R}\times\sigma _{v_L}: \qquad 
\Aut_k(F_{\o v})\times\Aut_k(F_{v}) \quad \To \quad
G^{(2n)}(F_{\o v}\times F_{v})\]
introduced in section~1.5 and generating the abstract bisemivariety
$G^{(2n)}(\Fvt)$ from 
$\GL_{2n}(\Fvt)$ as well as from the definition of the global elliptic
$\Gamma _{\wh M^{(2n)}_{v^T\RL}}$-bisemimodule $\phi ^{(2n)}\RL(x)$ given in section~2.12.\epr
\vskip 11pt

\subsection{Proposition}

{\em 
The shifted global elliptic bisemimodule
$\ELLIP\RL(2n[2k],r)$ resulting from the action
\[
(D^{2k}_R\otimes D^{2k}_L): \qquad
\ELLIP\RL(2n,r) \quad \To \quad\ELLIP\RL(2n[2k],r) \]
of the bioperator $(D^{2k}_R\otimes D^{2k}_L)$ on the global elliptic
$\Gamma _{\wh M^{(2n)}_{v^T\RL}}$-bisemimodule $\phi ^{(2n)}\RL(x)$, noted here
$\ELLIP\RL(2n,r)$, and generated under the (bi)monomorphism:
\begin{multline*}
\sigma _{v^T_R\otimes\rit}\times\sigma _{v^T_L\otimes\rit}:\qquad
\Aut_k(F^T_{\o v}\times \rit)\times \Aut_k(F^T_{v}\times \rit) \\
\quad \To \quad
G^{(2n)}((\FvTtrit))\end{multline*}
gives rise to {\bbf the eigenbivalue equation:
\[
(D^{2k}_R\otimes D^{2k}_L)(\ELLIP\RL(2n,r))=
E_{2k\RL}(2n,j)(\ELLIP\RL(2n,r))\ , \; 1\le j\le r\;, \]
where the eigenbivalues $E_{2k\RL}(2n,j)$ are shifts in $2k$ real dimensions of the global Hecke bicharacters
$\lambda (2n,j,m_j)$ associated with the subbisemilattices\/} characterized by the global residue degrees $j$.
}
\vskip 11pt

\bpr
Similarly as in proposition~2.13, the (bi)map:
\[
\sigma _{v^T_R\otimes\rit}\times\sigma _{v^T_L\otimes\rit}: \qquad
\Aut_k(F^T_{\o v}\times\rit)\times\Aut_k(F^T_{v}\times\rit) \quad \To \quad
G^{(2n)}((F^T_{\o v}\times\rit)\times (F^T_{v}\times\rit))\;,\]
where $\Aut_k(F^T_{\o v}\otimes\rit)\times\Aut_k(F^T_{v}\otimes\rit)$ is the bilinear semigroup of automorphisms of fibered or shifted toroidal transcendental extensions, is responsible for the generation of the toroidal abstract fibered or shifted bisemivariety
$G^{(2n)}((F^T_{\o v}\times \rit)\times(F^T_{v}\times \rit))$, referring to section~2.1, of which functional representation space is the shifted global elliptic bisemimodule
$\ELLIP\RL(2n[2k],r)$ obtained from the global elliptic bisemimodule
$\ELLIP\RL(2n,r)$ under the action of the bioperator
$(D^{2k}_R\otimes D^{2k}_L)$ giving rise to the eigenbivalue equation:
\[
(D^{2k}_R\otimes D^{2k}_L)(\ELLIP\RL(2n,r))
=E^{2k}\RL(2n,j)(\ELLIP\RL(2n,r))\;, \qquad \forall\ j\;, \quad 1\le j\le r\;.\]
The functional representation space
$\FREPSP(\GL_{2n}(F^T_{\o v}\times\rit) \times (F^T_{v}\times\rit))$ of the bilinear semigroup of matrices
$\GL_{2n} ((  F^T_{\o v}\times\rit) \times (F^T_{v}\times\rit )) $ responsible for the generation of the abstract bisemivariety
$G^{2n} ((  F^T_{\o v}\times\rit) \times (F^T_{v}\times\rit )) $ is given by the set of embedded bisemifunctions:
\begin{multline*}
\phi_1 ( G^{(2n)}(F^T_{\o v_{1}}\times\rit))\times \phi_1 ( G^{(2n)}(F^T_{v_{1}} \times\rit))  \\ \subset \dots\subset
\phi_{j,m_j} ( G^{(2n)}(F^T_{\o v_{j,m_j}}\times\rit))\times \phi_{j,m_j} ( G^{(2n)}(F^T_{v_{j,m_j}} \times\rit)) \\
\subset \dots\subset
\phi_{r,m_r} ( G^{(2n)}(F^T_{\o v_{r,m_r}}\times\rit))\times \phi_{r,m_r} ( G^{(2n)}(F^T_{v_{r,m_r}}\times\rit))
\end{multline*}
introduced in section~5.12. of \cite{Pie3}.

Each bisemifunction
$\phi_{j,m_j} ( G^{(2n)}(F^T_{\o v_{j,m_j}}\times\rit))  \times \phi_{j,m_j} ( G^{(2n)}(F^T_{v_{j,m_j}}\times\rit )) $ is the product, right by left,
$T^{(2n)}_R([2k],(j,m_j))\times T^{(2n)}_L([2k],(j,m_j))$ of a $2n$-dimensional real semitorus\linebreak
$T^{(2n)}_R([2k],(j,m_j))$ shifted in $2k$ real dimensions and localized in the lower half space by its symmetric equivalent $ T^{(2n)}_L([2k],(j,m_j))$  localized in the upper half space.

They have for analytic development:
\begin{align*}
T^{(2n)}_L([2k],(j,m_j)) \simeq E_{2k_L}(2n,j,m_j))\ \lambda (2n,j,m_j)\  e^{2\pi ijx}\\[11pt]
\rresp{T^{(2n)}_R([2k],(j,m_j)) \simeq E_{2k_R}(2n,j,m_j))\ \lambda (2n,j,m_j)\ e^{-2\pi ijx}}
\end{align*}
referring to section~2.12, where $E_{2k}(2n,j,m_j)$ is the shift in $2k$ real dimensions of the global Hecke character 
$\lambda (2n,j,m_j)$ being also a product of eigenvalues of the Hecke operator as described in section~2.12.

On the other hand, referring to proposition~2.11, the toroidal spectral representation of
$(D^{2k}_R\otimes D^{2k}_L)$ is given by the set of $r$-tuples:
\[
\ellip\RL(2n,1) 
\subset \dots \subset\ellip\RL(2n,j)  
\subset \dots \subset\ellip\RL(2n,r) \]
where
$\ellip\RL(2n,j) $ is given by
\[
\ellip\RL(2n,j) =
(\lambda (2n,j,m_j)\ e^{-2\pi ijx})\times
(\lambda (2n,j,m_j)\ e^{+2\pi ijx})\]
and to which corresponds the set of increasing eigenbivalues
\[
E_{2k\RL}(2n,1) \subset\dots\subset
E_{2k\RL}(2n,j) \subset\dots\subset
E_{2k\RL}(2n,r)\]
where $E_{2k\RL}(2n,j)$ is the shift in $2k$ real dimensions of the Hecke bicharacter\linebreak
$\lambda ^2(2k,j,m_j)    \subset \lambda ^2(2n,j,m_j)    $ taking into account that
\[
 \lambda ^2(2n,j,m_j)    =
 \prod_{c=1}^{2k}\lambda ^2_c(2k,j,m_j)\times
 \prod_{d=2n-2k}^{2n}\lambda ^2_d(2n,j,m_j)\]
 \begin{align*}
 & \text{with} & \lambda ^2(2k,j,m_j)=&\prod\limits_{c=1}^{2k} \lambda ^2_c(2k,j,m_j)\\[11pt]
&\text{and with} & \lambda ^2(2n-2k,j,m_j)=&\prod\limits_{d=2n-2k}^{2n} \lambda ^2_d(2n,j,m_j)\tag*{\eop}
\end{align*}
\vskip 11pt

\section{Large random matrices and Riemann zeta function}

\subsection{Five questions to find a solution to this problem}

Chapters 1 and 2 have been devoted to the mathematical tools necessary to clarify the conceptual framework behind the random matrices and, more particularly, the closed numerical connection between the spacings of the nontrivial zeros of the Riemann zeta function and the spacing of the eigenvalues of typical large random matrices.

It will be shown in this chapter that {\bbf the symmetry behind the Gaussian unitary ensemble (GUE) is the symmetric (bisemi)group of ``Galois'' automorphisms fibered or ``shifted'' algebraic and transcendental (bi)quanta.

The constant reference to the global program of Langlands\/} in chapter~1 {\bbf and to the geometric-shifted global program of Langlands, as well as to von Neumann (bisemi)algebras\/}, in chapter~2, {\bbf is thus not fortuitous\/}.

In order to find a solution to this problem, it will be answered in this chapter to {\bbf the five following questions\/}:
\Bena
\item {\bbf What is behind random matrices leading to GOE (Gaussian orthogonal ensemble), as well as GUE (Gaussian unitary ensemble)?\/}

\item {\bbf What is behind the partition and correlation function(s) between eigenvalues of random matrices?\/}

\item {\bbf What interpretation can we give to the local spacings between the eigenvalues of large random matrices?\/}

\item {\bbf What interpretation can we give to the  spacings between the nontrivial zeros of $\zeta  (s)$?\/}

\item {\bbf What is the curious connection between 3) and 4)?\/}
\Ee

But, first, we would like to outline that {\bbf the geometric dimension\/} envisaged in this chapter {\bbf will be one real, and, possibly, one complex\/}, i.e. that $n=1$, in the sense of the reducible global program of Langlands developed in \cite{Pie2}.  Thus, only curves, covering possibly surfaces, will be considered in order to meet the conditions of question~4).
\vskip 11pt

\subsection{The first question ``What is behind random matrices leading to GOE and GUE?''}

This question reflects the importance of large random matrices in the present development of mathematics and of physics.

These random matrices lead to the eigenvalue problem in the frame of the geometric shifted global program of Langlands, recalled in section~1.11 and in chapter~2, in order to find a response to question~3).
\vskip 11pt

\subsection{Bilinear differential Galois semigroup}

The symmetry group behind or at the origin of the bilinear global program of Langlands is
{\bbf the bilinear semigroup of automorphisms
$\Aut_k(F_{\o v})\times\Aut_k(F_{v})$
\resp{Galois automorphisms
$\Gal(F_{\o v}/k)\times\Gal(F_{v}/k)$}
of compact transcendental \resp{algebraic} quanta generating a bisemilattice of compact transcendental \resp{algebraic} quanta\/} referring to section~1.4.

Similarly, the symmetry group at the origin of the geometric bilinear global program of Langlands is the bilinear semigroup of fibered or shifted automorphisms
$\Aut_k(F_{\o v}\times\rit)\times\Aut_k(F_{v}\times\rit)$ of compact transcendental quanta generating a bisemilattice of compact shifted transcendental quanta according to proposition~2.4.

Referring  to this same proposition, the one-dimensional shifted functional representation space
$\FREPSP(\GL_1(\Fvtrit))$ of
$\GL_1(\Fvtrit)$ is given by the shifted bisemisheaf
$\wh M^{(1)}_{v_R}[1]\otimes\wh M^{(1)}_{v_L}[1])$ whose set
$\Gamma (\wh M^{(1)}_{v_R}[1]\otimes\wh M^{(1)}_{v_L}[1])$ of bisections is the set
$\{ \phi   G^{(j,m_j)}_{(Q)}(F_{\o v^1}\times\rit)\times (F_{v^1}\times\rit))\}_{j,m_j}$ of fibered or shifted differentiable bifunctions obtained from the set
$\{ \phi (G^{(j,m_j)}_{(Q)}(F_{\o v^1}\times F_{v^1}))\}_{j,m_j}$ of differentiable bifunctions under the action of the elliptic bioperator
$D_R\otimes D_L$, $2k=1$.

The set 
\[ 
\{ \phi (G^{(j,m_j)}_{(Q)}(F_{\o v^1}\times\rit)\times (F_{v^1}\times\rit))\}_{j,m_j}
=
\{ \phi_R (G^{(j,m_j)}_{(Q)}(F_{\o v^1}\times\rit))
\otimes
\phi_L (G^{(j,m_j)}_{(Q)}(F_{\o v^1}\times\rit))\}_{j,m_j}
\]
of differentiable bifunctions, being equivalent to the set
$\{ \phi (G^{(1)}(F_{\o v_{j,m_j}}\times\rit)\times (F_{v_{j,m_j}}\times\rit))\}_{j,m_j}$ of differentiable bifunctions,
is isomorphic to the set
$\{ \phi (G^{(j,m_j)}_{(Q)}(F_{\o v^1}\times F_{v^1}))\times
\phi (G^{(j,m_j)}_{(Q)}(\rit\times \rit)) \}_{j,m_j}$ of bifunctions in such a way that
{\bbf
$\{ \phi (G^{(j,m_j)}_{(Q)}(\rit\times \rit))\}_{j,m_j=1}^r$ is the fibering or shifting functional representation bisemispace obtained under(and being isomorphic to)  the biaction of the bilinear differential Galois semigroup\/}
\begin{align*}
\Aut_k(\phi _R(\rit))\times
\Aut_k(\phi _L(\rit))
& \simeq\{\Aut_k(\phi _R(\rit))_{|F_{\o v_{j,m_j}}}\times
\Aut_k(\phi _L(\rit))_{|F_{v_{j,m_j}}} \}_{j=1}^r\\
&\simeq \GL_r(\phi _R(\rit)\otimes\phi _L(\rit))\\
&\simeq \GL_r(\rit\times\rit)\end{align*}
according to proposition~2.4.

Thus, the bilinear semigroup $\GL_r(\rit\times\rit)$ of matrices constitutes the representation of the bilinear differential Galois semigroup associated with the action of the differential bioperator $(D_R\otimes D_L)$ on the bisemisheaf
$(\wh M^{(1)}_{v_R}\otimes \wh M^{(1)}_{v_L})$ whose bisections are the set
$\{ \phi (G^{(j,m_j)}_{(Q)}(F_{\o v^1}\times F_{v^1})) \}_{j=1}^r=
\{ \phi (G^{(1)}(F_{\o v_{j,m_j}}\times F_{v_{j,m_j}})) \}_{j,m_j}$ of differentiable bifunctions.

Let us set 
\[
\phi _j(G_{(Q)}(F_{\o v^1}\times F_{v^1})) =\sum_{j=1}^j
 \phi (G^{(j,m_j)}_{(Q)}(F_{\o v^1}\times F_{v^1})) \]
 used in the following proposition.
 \vskip 11pt
 
 \subsection{Proposition}
 
 {\em
 If
 \[
\phi _r(G_{(Q)}(F_{\o v^1}\times F_{v^1})) =\sum_{j=1}^r
 \phi (G^{(j,m_j)}_{(Q)}(F_{\o v^1}\times F_{v^1}))\;, \]
 let
 \[
 (D_R\otimes D_L)
 (
\phi _r(G_{(Q)}(F_{\o v^1}\times F_{v^1})) = E\RL(j)\
 \phi_R (G_{(Q)}(F_{\o v^1}\times F_{v^1}))\;, \qquad 1\le j\le r\;, \]
 be the eigenbivalue equation related to the bisemisheaf
 $(\wh M^{(1)}_{v_R}\otimes \wh M^{(1)}_{v_L})$ and associated with the tower of shifted differentiable bifunctions
 $
\{\phi (G^{(j,m_j)}_{(Q)}((F_{\o v^1}\times \rit)\times (F_{v^1}\times\rit)))\}_{j=1}^r
$.

Then, {\bbf
the eigenbivalues of the matrix of $\GL_r(\rit\times\rit)$, constituting a representation of the bilinear differential Galois semigroup associated with the biaction of $(D_R\otimes D_L)$, are the eigenbivalues 
$\{E\RL(j)\}_{j=1}^r$ of the above eigenbivalue equation\/}.}
\vskip 11pt

\bpr
\Bena
\item %The matrix of 
$\GL_r(\rit\times\rit)$, being:
\Bean
\item the bilinear fibre $\Fs\RL(\TAN)$ of the tangent bibundle
$\TAN (\wt M^{(1)}_{v_R}\otimes\wt M^{(1)}_{v_L})
\simeq \AdF\REPSP(\GL_1(F_{\o v}\times F_v)$
according to section~2.2,
\item a representation of the bilinear differential Galois semigroup associated with the biaction of $(D_R\otimes D_L)$,
\Ee
constitutes a representation of the bioperator $(D_R\otimes D_L)$ because it generates endomorphisms of $\TAN (\wt M^{(1)}_{v_R}\otimes\wt M^{(1)}_{v_L})$.

\item Referring to propositions~2.11 and 2.14, we see that {\bbf the set of differentiable bifunctions\/}
$\{\phi (G^{(j,m_j)}_{(Q)}(F_{\o v^1}\times F_{v^1})) \}_j$ constituting the
$r$-bituple
\[ \langle \phi (G^{(1)}_{(Q)}(F_{\o v^1}\times F_{v^1})),\dots,
 \phi (G^{(j,m_j)}_{(Q)}(F_{\o v^1}\times F_{v^1})),\dots,
 \phi (G^{(r,m_r)}_{(Q)}(F_{\o v^1}\times F_{v^1}))\rangle
 \]
 {\bbf is the spectral representation of $(D_R\otimes D_L)$\/} (and the basis)
 {\bbf associated with the $r$-tuple of eigenbivalues\/}:
 \be
 \langle E\RL(1),\dots, E\RL(j),\dots, E\RL(r)\rangle\;.\tag*{\eop}\ee
 \Ee
 \vskip 11pt
 
 \subsection{Corollary}
 
 {\em
 As the $j$-th eigenbifunction
 $\phi (G^{(j,m_j)}_{(Q)}(F_{\o v^1}\times F_{v^1}))$ is a bifunction on ``$j$'' transcendental compact biquanta, the corresponding eigenbivalue $E\RL(j)$ will be the shift of this bifunction corresponding to the biaction of the bioperator $(D_R\otimes D_L)$.
 }
 \vskip 11pt
 
 \bpr
 We refer to section~2.10 and proposition~2.11 concerning;
 \Bena
 \item the bisemialgebra of von Neumann 
 $\MM\RL(H^+_{\wh M^{(1)}_{v\RL}})$ on the tower of embedded bilinear Hilbert semispaces associated with the bisemisheaf
 $\wh M^{(1)}_{v\RL}\equiv 
 \wh M^{(1)}_{v_R}\otimes\wh M^{(1)}_{v_L}$
 
 and
 
 \item the discrete spectrum $\sum(D_R\otimes D_L)$ of the differential bioperator $(D_R\otimes D_L)$.\epr
 \Ee
 \vskip 11pt
 
 \subsection[Bilinear Gaussian orthogonal (resp. unitary) ensemble BGOE (resp. BGUE)]{Bilinear Gaussian orthogonal (resp. unitary) ensemble\linebreak BGOE (resp. BGUE)}
 
 The {\bbf Gaussian orthogonal \resp{unitary} ensemble GOE \resp{GUE} is defined in the space of real symmetric \resp{hermitian} matrices by two requirements \cite{Meh}\/}:
 \Bean
 \item the ensemble is invariant under every transformation:
 \[ H \To W^T\ H\ W  \qquad \rresp{H \To U^{-1}\ H\ U}\]
 where:
 \Bi
 \item $W$ \resp{$U$} is any real orthogonal \resp{unitary} matrix;
 \item $H$ is a real symmetric \resp{hermitian} matrix, generally related to the hamiltonian matrix invariant \resp{not invariant} under time reversal.
 \Ei
 
 \item the various elements $H_{ij}$ are statistically independent.
 \Ee
 \vskip 11pt
 
 In its simplified bilinear version, the hamiltonian $H$ corresponds to the differential (random) bioperator $(D_R\otimes D_L)$ acting on the bisemisheaf
 $\wh M^{(1)}_{v\RL}$ and belonging to the bisemialgebra of von Neumann
 $\MM\RL(H^+_{\wh M^{(1)}_{v\RL}})$ according to sections~2.8 to 2.11 and \cite{Pie1}.
 
 Let then:
 \[ D_R\otimes D_L : \qquad \wh M^{(1)}_{v_R}\otimes \wh M^{(1)}_{v_L}
 \quad \To \quad \wh M^{(1)}_{v_R}[1]\otimes \wh M^{(1)}_{v_L}[1]\]
 be the biaction of $(D_R\otimes D_L)$ on the bisemisheaf
 $(\wh M^{(1)}_{v_R}\otimes \wh M^{(1)}_{v_L})$ generating the perverse bisemisheaf
 $(\wh M^{(1)}_{v_R}[1]\otimes\wh M^{(1)}_{v_L}[1])$ whose bisections
 $\{ \phi_{j,m_j} (G^{(1)}_{[1]}(F_{\o v_{j,m_j}}\times\rit)\times (F_{v_{j,m_j}}\times\rit))\}_{j,m_j}$
 are differentiable bifunctions on the $(j,m_j)$-th conjugacy class representatives of the shifted or fibered bilinear semigroup
 $G^{(1)}((F_{\o v}\times\rit)\times (F_{v}\times\rit))$.
 
 Referring to \cite{Pie1}, {\bbf the perverse bisemisheaf
 $(\wh M^{(1)}_{v_R}[1]\otimes \wh M^{(1)}_{v_L}[1])$ is the operator-valued stringfield of an elementary (bisemi)particle\/}.
 
 According to section~3.3, {\bbf
 the bilinear semigroup of matrices $\GL_r(\rit\times\rit)$ constitutes the representation of the bilinear differential Galois semigroup associated with the biaction of 
 $(D_R\otimes D_L)$ on
 $(\wh M^{(1)}_{v_R}\otimes\wh M^{(1)}_{v_L})$\/}.
 
 And, thus, {\bbf $\GL_r(\rit\times\rit)$ (and $O_r(\rit\times\rit)$), or its complex equivalent
$\GL_r(\cit\times\cit)$ (and $U_r(\rit\times\rit)$), is the new bilinear Gaussian real (orthogonal) \resp{complex (``unitary'')} ensemble labeled BGOE \resp{BGUE} corresponding to GOE \resp{GUE}\/}.
%
%Remark the complex entries $(\cit\times\cit)$ correspond at a sign to the real entries $(\rit\times\rit)$ because we have in fact that $\cit\times\cit$ is $\cit^*\times\cit$ and $\rit\times\rit$ is $\rit_-\times\rit_+$ verifying $\cit^*\times \cit \simeq \rit_-\times\rit_+$.
 \vskip 11pt
 
 \subsection{Mixed higher bilinear $KK$-theory}
 
 To be complete, the deformations of random matrices have to be envisaged in the light of
 {\bbf the new interpretation of homotopy-cohomotopy\/} \cite{Pie5}
 {\bbf viewed as deformations of Galois representations\/}
 in the context of mixed higher bilinear $KK$-theory related to the Langlands dynamical bilinear global program \cite{Pie5}.
 \vskip 11pt
 
 Referring to proposition~2.4, we see that the $1D$-geometric infinite quantum general bilinear semigroup
 $\GL^{(Q)}(F_{\o v^1}\times F_{v^1})$ is defined by:
 \[
 \GL^{(Q)}(F_{\o v^1}\times F_{v^1})=\lim_{j=1\to r} \GL_j^{(Q)}(F_{\o v^1}\times F_{v^1})
 \]
 in such a way that $\GL_1^{(Q)}(F_{\o v^1}\times F_{v^1})$ is the parabolic, i.e. unitary bilinear semigroup
 $P_1(F_{\o v^1}\times F_{v^1})$, and that its shifted equivalent 
 $\GL^{(Q)}((F_{\o v^1}\times \rit)\times (F_{v^1}  \times\rit))  $ is given by
 \[
\GL^{(Q)}((F_{\o v^1}\times \rit)\times (F_{v^1}  \times\rit))=
\lim_{j=1\to r} \GL^{(Q)}_j((F_{\o v^1}\times \rit)\times (F_{v^1}  \times\rit))
\]
leading to a filtration
\[
G^{(1)}_{(Q)}((F_{\o v^1}\times \rit)\times (F_{v^1}  \times\rit))\subset\dots\subset
G^{(r,m_r)}_{(Q)}((F_{\o v^1}\times \rit)\times (F_{v^1} \times\rit))  \]
of its representatives.
\vskip 11pt

So, {\bbf the bilinear version of the algebraic $K$-theory restricted to $1D$-geometric dimension is given by\/}:
\[
{K^1(G^{(1)}(F_{\o v^1}\times F_{v^1}))=
\Pi _1(\BGL^{(Q)}(F_{\o v^1}\times F_{v^1})^+)}\]
where the quantum classifying bisemispace
$\BGL^{(Q)}(F_{\o v^1}\times F_{v^1})$ is the base bisemispace of all equivalence classes of deformations of the Galois representations of
$\GL^{(Q)}(\wt F_{\o v^1}\times \wt F_{v^1})$ given by the kernels 
$\GL^{(Q)}(\delta F_{\o {v^1+\ell} }\times \delta F_{v^1+\ell })$ of the maps:
\[
\GD^{(Q)}_\ell: \qquad
\EGL^{(Q)}(F_{\o {v^1+\ell }}\times F_{v^1+\ell }) \quad \To \quad
\BGL^{(Q)}(F_{\o v^1 }\times F_{v^1 }) \]
where $\GD^{(Q)}_\ell$ is a universal principal $\GL^{(Q)}(F_{\o v^1}\times F_{v^1})$-bibundle.

Referring to chapter~3 of \cite{Pie5}, it is easy to see that
\[
\BGL^{(Q)}(F_{\o v^1}\times F_{v^1})\equiv
\GL_1(F_{\o v}\times F_{v})\]
 and that the maps $\GD^{(Q)}_\ell$ become
 \[ \GD^{(1)}_\ell : \qquad 
 \GL_1(F_{\o {v+\ell} }\times F_{v+\ell }) \quad \To \quad
\GL_1(F_{\o v}\times F_{v }) \]
having the same interpretation as $\GD^{(Q)}_\ell$.
\vskip 11pt
In order to recall the {\bbf bilinear version of the mixed higher $KK$-theory of Quillen adapted to the Langlands dynamical global program in $1D$-geometric dimension\/} \cite{Pie5},  we have to take into account:
\Bean
\item the bisemisheaf
$(\wh M^{(1)}_{v_R}\otimes \wh M^{(1)}_{v_L})$, noted here
$\FG^{(1)}(F_{\o v }\times F_{v })$ and being the functional representation space of
$\GL_1(F_{\o v }\times F_{v }) $;

\item the ``plus'' classifying bisemisheaf
$\BFGL^{(Q)}(\rit\times\rit)^+= \BFGL_1(\rit\times\rit)^+$, being the base bisemisheaf of all equivalence classes of one-dimensional inverse deformations of the Galois differential representation of
$\FGL_1(\rit\times\rit)$ due to the action of the bioperator $(D_R\otimes D_L)$ on
$\FG^{(1)}(F_{\o v}\times F_v)$, and corresponding to one-dimensional deformations of the Galois representation of
$\GL_1(F_{\o v}\times F_v)$ given by the kernel
$\{\GL_1(\delta F_{\o {v+\ell }}\times \delta F_{v+\ell })\}_\ell$ of $\GD^{(1)}_\ell$.
\Ee

{\bbf The higher (algebraic) $KK$-theory is then given by:}
\begin{multline*} 
\mbox{\boldmath {$K_1(\FG^{(1)}(\rit\times\rit))\times K^1(\FG^{(1)}(F_{\o v}\times F_v))$}}\\
\mbox{\boldmath{$=
\Pi ^1(\BFGL_1(\rit\times\rit)^+)\times \Pi _1(\BFGL_1(F_{\o v}\times F_v)^+)$}}
\end{multline*}
 in such a way that the bilinear contracting $K$-theory
$K_1(\FG^{(1)}(\rit\times\rit))$ responsible for a differentiable biaction acts on the $K$-theory
$K^1(\FG^{(1)}(F_{\o v}\times F_v))$ of the bisemisheaf
$\FG^{(1)}(F_{\o v}\times F_v)$ in one-to-one correspondence with the biaction of the cohomotopy bisemigroup
$\Pi ^1(\BFGL_1(\rit\times\rit)^+)$ of the ``plus'' classifying bisemisheaf 
$\BFGL_1(\rit\times\rit)^+)$.
\vskip 11pt

\subsection{Proposition}

{\em
The deformation
\[
\GD^{(1)}_\ell: \qquad \GL_1( F_{\o {v+\ell }}\times F_{v+\ell }) \quad \To \quad
\GL_1(  F_{\o v  }\times F_{v  })\]
induces the following deformation:
\[
\GL_{r+\delta r_\ell\to r}: \qquad
\GL_{r+\delta r_\ell}(\rit\times\rit) \quad \To \quad
\GL_{r}(\rit\times\rit) \]
on random matrices 
$\GL_{r}(\rit\times\rit)$.
}
\vskip 11pt

\bpr
Indeed, the matrix of $\GL_{r+\delta r_\ell}(\rit\times\rit) $ of order $(r+\ell)$ constitutes a representation of the deformed bioperator or
$(( D_R + \delta D_R)\times ( D_L + \delta D_L))$ acting on the deformed bifunction
$\phi _{r+\delta r_\ell}(G_{(Q)}(F_{\o v^1}\times F_{v^1}))$ and has for spectral representation the
$(r+\delta r_\ell)$-bituple:
\[
\langle
\phi (G^{(1+\ell)}_{(Q)}( F_{\o {v^1+\ell }}\times F_{v^1+\ell })),\dots,
\phi (G^{(r+\ell,m_r+\ell)}_{(Q)}( F_{\o {v^1+\ell }}\times F_{v^1+\ell }))
\rangle\]
decomposing into the sum of the $r$-bituple and the $\delta r_\ell$-bituple according to:
\begin{multline*}
\langle
\phi (G^{(1+\ell)}_{(Q)}( F_{\o {v^1+\ell }}\times F_{v^1+\ell })),\dots,
\phi (G^{(r+\ell,m_r+\ell)}_{(Q)}( F_{\o {v^1+\ell }}\times F_{v^1+\ell }))
\rangle\\
=
\langle\phi (G^{(1)}_{(Q)}( F_{\o v^1  }\times F_{v^1 })),\dots,\phi (G^{(r,m_{r})}_{(Q)}( F_{\o v^1  }\times F_{v^1} ))\rangle \\
+
\langle\phi (G^{(\ell)}_{(Q)}( F_{\o v^1_\ell  }\times F_{v^1_\ell })),\dots,\phi (G^{(\ell)}_{(Q)}( F_{\o v^1_\ell  }\times F_{v^1_\ell }))\rangle\;.\end{multline*}
The deformation $\GL_{r+\delta r_\ell}(\rit\times\rit)$ of the random matrix of
$\GL_{r}(\rit\times\rit)$, constituting a deformation of the Galois differential representation
$\GL_{r}(\rit\times\rit)$, corresponds to the cohomotopy bisemigroup
$\Pi ^1(\BFGL_{1}(\rit\times\rit)^+)$ according to section~3.7.\epr
\vskip 11pt

\subsection{The second question}

{\bbf The second question ``What is behind the partition and correlation functions between eigenvalues of random matrices?''\/} concerns the distribution of eigenvalues of random matrices.  It will be seen that this problem {\bbf is based on the (bisemi)group of ``Galois'' automorphisms of shifted transcendental compact (bi)quanta which leads to a reevaluation of the probabilistic interpretation in quantum theories.}
\vskip 11pt

\subsection{Distribution of eigenvalues of GUE ensembles}

Wigner introduced the idea of {\bbf 
statistical mechanics of nuclei based on a Gaussian ensemble (GUE) having ``$r$'' quantum states and characterized by a Hamiltonian symmetric matrix of order ``$r$'' whose entries are Gaussian random variables and to which a Gaussian statistical weight is associated\/} \cite{Wig}.

Unsatisfied by the impossibility of defining a uniform probability distribution on an infinite range, {\bbf
F. Dyson introduced the circular unitary \resp{orthogonal} ensemble CUE \resp{COE}\/} \cite{Dys} in such a way that the Hamiltonian is now described by a unitary matrix of order ``$r$'' whose eigenvalues are complex numbers $\exp(i\theta _j)$, $1\le j\le r$, distributed around the unit circle.

This circular unitary \resp{orthogonal} ensemble corresponds to the Riemann symmetric space $U(r)/O(r)$ which ``lives'' in $\GL_r(\cit)/U(r)$ \cite{Mez}.

Let then $M=(M_{ij})^r_{i,j=1}$ be a random hermitian matrix to which is assigned the probability distribution
\[
d{\mu _r}^{\rm GUE}(M) = \frac1{Z^{\rm GUE}_r}\ e^{-r\Tr M^2}\ dM\]
where:
\Bi
\item $\Tr M^2 =\sum_{i,j=1}^r\ M_{ji}M_{ij}=\sum^r_{i=1} M^2_{ii}+2\sum_{i>j}|M_{ij}|^2$;

\item $\mu ^{\rm GUE}_r(dM) =\frac1{Z^{\rm GUE}_r}\ \prod^r_{i=1} (e^{-rM^2_{ii}})
\ \prod_{i>1} (e^{-2r|M_{ij}|^2})\ dM$.
\Ei.

{\bbf The distribution of eigenvalues of $M$ with respect to the ensemble
$\mu ^{\rm GUE}_r$ is\/} \cite{Meh}, \cite{A-VM}, \cite{B-Z}:
\[
d\mu ^{\rm GUE}_r(\lambda )=\frac1{Z^{\rm GUE}_r}\prod_{i>j}(\lambda _i-\lambda _j)^2\prod_{i=1}^re^{-r\lambda ^2_i}\ d\lambda \]
where the partition function \cite{Rue} of GUE is given by:
\begin{align*}
Z^{\rm GUE}_r
&= \int \prod_{i>j}(\lambda _i-\lambda _j)^2\prod_{i=1}^re^{-r\lambda ^2_i}\ d\lambda _i\\[11pt]
&= \frac{(2r )^{r/2}}{(2r )^{r^2/2}}\prod_{i=1}^r\ i!\;.\end{align*}
\vskip 11pt

\subsection{$m$-point correlation function for GUE and Jacobi matrix}

{\bbf
The joint probability density function for the eigenvalues of matrices from a Gaussian orthogonal, unitary or symplectic ensemble is given by\/}:
\[
P_{r\beta }(x_1,\dots,x_r)
= c_{r\beta }\exp\L( -\tfrac\beta 2\sum_{i=1}^rx_i^2\R) \prod_{i>j} (x_i-x_j)^\beta \;, \qquad -\infty <x_i<+\infty \;,\]
where $\beta =1,2$ or $4$ according as the ensemble is orthogonal, unitary or symplectic and $c_{r2}=\dfrac1{Z^{\rm GUE}_r}$.
\vskip 11pt

{\bbf The $m$-point correlation function for the Gaussian unitary ensemble is defined by} \cite{Dys}:
\[ 
R_{mr}(x_1,\dots,x_m)=\frac{r!}{(r-m)!}\int_{\rit^{r-m}}P_r(x_1,\dots,x_r)\ dx_{m+1}\ \dots\ dx_r\]
which is the probability density of finding a level around each of the points (i.e. entries of $M$) $x_1,\dots,x_m$, the positions of the remaining levels being unobserved.

{\bbf The Dyson determinantal formulas for correlation functions is} \cite{Meh}, \cite{Ble}, \cite{B-H}:
\[ R_{mr}(x_1,\dots,x_m)=\det (K_r(x_k,x_\ell))^m_{k,\ell=1}\]
where $K_r(x_k,x_\ell)$ is given by:
\[
 K_r(x_k,x_\ell) = \sum_{i=0}^{r-1}\psi _i(x_k)\ \psi _i(x_\ell)\qquad \qquad
\text{with}\quad \psi _i(x)=h_i^{-\half}P_i(x)\ e^{-rM^2(x)/2}\]
where $P_i(x)$ is an orthogonal polynomial or degree $i$ corresponding to the weight function $e^{-rM^2(x)/2}$ and verifying:
\[
\int_{-\infty }^{+\infty }P_i(x)\ P_j(x)\ e^{-rM^2(x)}\ dx\simeq \delta _{ij}\ h_i\;.\]
\vskip 11pt

{\bbf These orthogonal polynomials\/}, being sometimes Hermite polynomials, {\bbf satisfy the three term
 recurrent relation\/}:
 \begin{align*}
\beta _{i+1}\ P_{i+1}(x)
 &= (x-\alpha _i)\ P_i(x)-\beta _i\ P_{i-1}(x)\\
 \text{or} \qquad \qquad \qquad
 x\ P_i &= \beta _i\ P_{i-1} + \alpha _i\ P_i+\beta _{i+1}\ P_{i+1}\;, \qquad
 \beta _i = \frac{h_i}{h_i-1}\;, \quad \alpha _i=(x\ P_i,P_i)\;.\end{align*}
 If we set $P_{-1}(x)=0$, we get the tower:
 \begin{align*}
 x\ P_0 &= \alpha _0\ P_0 + \beta _1\ P_1\;, \\
 x\ P_1 &= \beta _1\ P_0+\alpha _1\ P_1 + \beta _2\ P_2\;, \\
 x\ P_2 &= \quad 0 \quad + \beta _2\ P_1+\alpha _2\ P_2 + \beta _3\ P_3+0\;, \\
 \vdots \ &\\
 x\ P_{i-1} &= \qquad 0 \qquad + \beta _{i-1}\ P_{i-2}+\alpha _{i-1}\ P_{i-1} + \beta _i\ P_i\;, \end{align*}
 which, in matricial form, is:
 \[ x\ \Ps = J\ \Ps + \beta _i\ P_i\;.\]
 The matrix $J$ is symmetric and is the Jacobi matrix such that the $i$ roots of $P_i(x)$ verifying
 \[ P_i(x_j) = 0\;, \qquad 1\le j\le i\;, \]
 lead to
 \[ x_j[\Ps(x_j)]=J\ [P(x_j)]\;.\]
{\bbf The $i$ roots of $P_i(x)$ are then the eigenvalues of the Jacobi matrix $J$\/}.
 \vskip 11pt
 
 \subsection{Joint probability density function for the eigenvalues of matrices from BGUE and BGOE}
 
 In the bilinear case, the random  matrix corresponding to $M$ is
 \[ G=TG^T\times TG\in \GL_r(\rit\times\rit) \qquad \text{(or\ } \in \GL_r(\cit\times\cit)\ )\]
 where $TG\in T_r(\rit)$
\resp{$TG^T\in T^t_r(\rit)$} is an upper \resp{lower} triangular matrix of order $r$ with entries in $\rit$ (or $\cit$).

The BGUE (or BGOE) probability distribution corresponding to $G$ is:
\[
d^{\rm BGUE}_{\mu _r}(G) = \frac1{Z^{\rm BGUE}_r}\ e^{-rTr(TG^T\times TG)}\ dG\]
leading to
\[ \mu ^{\rm BGUE}_r(G) = \frac1{Z^{\rm BGUE}_r}\
\prod^r_{i=1} e^{-rG^2_{ii}}
\prod_{i>j} e^{-2r|G_{ij}|^2}\;.\]
\vskip 11pt

{\bbf The joint probability density function for the eigenvalues of matrices from a Gaussian bilinear orthogonal or unitary ensemble is thus\/}:
\begin{align*}
P_{r\RL}(x_1,\dots,x_r)
&= c_r\exp \ \L(-\sum_{i=1}^rrx_i^2\R)\ \prod_{i>j}(x_i-x_j)^2\\[11pt]
&= c_r \prod_{i=1}^r\ e^{-rx_i^2} \prod_{i>j}^r\ (x_i-x_j)^2\end{align*}
which corresponds to the distribution of eigenvalues $d^{\rm GUE}_{\mu _r}(\lambda )$ of a random matrix $M$ with respect to GUE, the eigenvalues 
$(x_1,\dots,x_r)$ being in fact eigenbivalues
$(x_1^2,\dots,x_r^2)$.

{\bbf
The $m$-point correlation function for the bilinear Gaussian (``unitary'') (or real (``orthogonal'')) ensemble BGUE (or BGOE)\/}:
\[ 
R_{mr\RL}(x_1^2,\dots,x_m^2)
= \frac{r!}{(r-m)!} \int_{\rit^{r-m}}P_{r\RL}(x_1^2,\dots,x_r^2)\ dx_{m+1}\ \dots dx_r\]
then corresponds to the $m$-point correlation function for GUE $  R_{mr}(x_1,\dots,x_m)$ developed in section~3.1.1.

$R_{mr\RL}(x_1^2,\dots,x_m^2)$ is thus also given by:
\[
R_{mr\RL}(x_1^2,\dots,x_m^2)=\det (K_r(x_k,x_\ell)^m_{k,\ell=1})\]
with:
\Bi
\item $K_r(x_k,x_\ell)=\sum\limits_{i=0}^{r-1}\psi _i(x_k)\ \psi _i(x_\ell)$

and

\item $\psi _i(x)=h^{-\half}P_i(x)\ e^{-r[(TG^T\times TG)(\rit\times\rit)]/2}$.
\Ei
As in section~3.11, $P_i(x)$ is an orthogonal polynomial of degree $i$ associated with  the weight function $e^{-r[(TG^T\times TG)(\rit\times\rit)]/2}$
where $TG^T(\rit)\times TG (\rit)\equiv G(\rit\times\rit)$is the bilinear Gauss decomposition of the matrix $G$.

These orthogonal polynomials $P_i(x)$ also satisfy the three term recurrent relation leading to the Jacobi matrix.
\vskip 11pt

\subsection{Proposition}

{\em
Let
\[
K_r(x,x)=\sum_{i=0}^{r-1}\psi _i(x)\ \psi _i(x)\]
be the energy level density with $\psi _i(x)$ given by:
\[
\psi _i(x)=h^{-\half} P_i(x)\ e^{-r[(TG^T\times TG)(\rit\times\rit)]/2}\;.\]
Then, we have that:
\Bena
\item {\bbf the squares of the roots $(x_j)$ of the polynomial $P_i(x)$ correspond to the eigenbivalues of the product, right by left, $(U_{r_R}\times U_{r_L})$ of the Hecke operators\/};

\item{\bbf the weight $e^{-r[(TG^T\times TG)(\rit\times\rit)]}$ is a measure on the eigenbivalues of the matrix $G\in \GL_r(\rit\times\rit)$ being a representation of the differential bilinear Galois semigroup\/}.
\Ee
}
\vskip 11pt

\bpr
\Bena
\item According to section~3.11, the roots of the orthogonal polynomial $P_i(x)$ are the eigenvalues of the Jacobi matrix $J$.

On the other hand, the set $\{P_j(x)\}_{j=1}^i$ of orthogonal polynomials envisaged in the tower of the three term recurrent relations leading to the matricial form
\[ x\ P=J\ P+\beta _i\ P_i\]
is in one-to-one correspondence with the set
$\{\phi (G^{(j)}_{(G)}  (F_{\o v^1}\times F_{v^1})) \}_{j=1}^i$ of differentiable eigenbifunctions, being the spectral representation of the bioperator $(D_R\otimes D_L)$ according to proposition~3.4 and constituting a representation of the bilinear semigroup of automorphisms 
$\Aut_k(F_{\o v})\times\Aut_k(F_{v})$ of compact transcendental biquanta generating a bisemilattice according to section~3.3.

{\bbf This bisemilattice results from the action of Hecke bioperators
$(U_{i_R}\otimes U_{i_L})$ generating the endomorphisms of the bisemisheaf
$(\wh M^1_{\o v_R}\times \wh M^1_{v_L})$} referring to sections~1.5 and 3.6.

Consequently, {\bbf the Jacobi matrix $J$ must be a representation of the Hecke operator $U_{i_L}$\/}, and, thus, the square of the roots of the polynomial $P_i(x)$ correspond to the eigenbivalues of the Hecke bioperator 
$(U_{i_R}\otimes U_{i_L})$ and are global Hecke (extended) bicharacters referring to \cite{Pie2}.

\item The weight $e^{-r[(TG^T\times TG)(\rit\times\rit)]/2}$ is thus a representation of the differential bilinear Galois semigroup
$\Aut_k(\phi _R(\rit))\times \Aut_k(\phi _L(\rit))$ shifting the product, right by left, of automorphism semigroups
$\Aut_k(\phi _R(F_{\o v}))\times\Aut_k(\phi _L(F_{v}))$ of cofunctions and functions on compact transcendental quanta according to proposition~2.4 and developments of N. Katz \cite{Kat}.

This representation of the bilinear differential Galois semigroup is then associated with the biaction of the differential bioperator $(D_R\otimes D_L)$ of which eigenbivalues are shifts of global Hecke (extended) bicharacters referring to section~3.4 and proposition~3.5.\epr
\Ee
\vskip 11pt

\subsection{Corollary}

{\em
The Jacobi matrix $J$ is a representation of the Hecke operator.}
\vskip 11pt

\bpr
This results directly from proposition~3.13.\epr
\vskip 11pt

\subsection{Proposition}

{\em
The probabilistic interpretation of quantum (field) theories is thus related to the bilinear semigroup of automorphisms
$\Aut_k(F_{\o v})\times\Aut_k(F_{v})$
\resp{$\Aut_k(\wt F_{\o v})\times\Aut_k(\wt F_{v})$}
of compact transcendental \resp{algebraic} biquanta generating a bisemilattice of these.
}
\vskip 11pt

\bpr
The probabilistic interpretation in QFT is given by the function density
\[
P_r(x,x)=\sum_{i=0}^{r-1} P_i(x)\ P_i(x)\]
where:
\Bi
\item $P_i(x)$ is an orthogonal polynomial of degree $i$;
\item $P_r(x,x)\ dx\in K_r(x,x)$ being the (wave) function density gives the probability of finding a (bisemi)particle in a volume element $(x,x+dx)$.
\Ei

As $P_i(x)\  P_i(x)$ or more exactly $P_i(-x)\ P_i(x)$, $x\in \rit_+$, constitutes a representation of the bilinear semigroup of automorphisms
$\Aut_k(F_{\o v})\times\Aut_k(F_{v})$ according to sections~1.5 to 2.12, we get the thesis.\epr
\vskip 11pt

\subsection{Proposition}

{\em
The $m$-point correlation function for BGUE (or BGOE)
$R_{mr}(x_1^2,\dots,x_m^2)$ constitutes a representation of a the bilinear semigroup of automorphisms
$\Aut_k(\phi _R(F_{\o v}\times\rit))\times\linebreak \Aut_k(\phi _L(F_{v}\times\rit))$
of fibered or shifted compact transcendental biquanta:
\[ \Rep_{\Aut_k(\phi \RL(F_{\o v\times\rit}))}: \qquad
\Aut_k(\phi _R(F_{\o v}\times\rit))\times\Aut_k(\phi _L(F_{v}\times\rit))
\quad \To \quad R_{mr}(x_1,\dots,x_m)\;.\]
}
\vskip 11pt

\bpr
We refer to proposition~2.4 and section~3.12 showing that, in
\[ K_r(x_k,x_\ell)=\sum_{i=0}^{r-1}\psi _i(x_k)\ \psi _i(x_\ell)\]
with
\[ \psi _i(x)=h^{-\half}\ P_i(x)\ \exp (-r[TG^T\times TG)(\rit\times\rit)]/2)\;, \]
\Bean
\item the polynomials $P_i(x)$ constitute a representation of the linear semigroup of automorphisms $\Aut_k(F_v)$ of compact transcendental quanta;

\item the weight $ \exp (-r[TG^T\times TG)(\rit\times\rit)])$ is a representation of the differential bilinear Galois semigroup
$\Aut_k(\phi _R(\rit))\times\Aut_k(\phi _L(\rit))$;

\item the products $\psi _i(x_k)\ \psi _i(x_\ell)$ constitute a representation of the bilinear semigroup of automorphisms
$\Aut_k(\phi _R(F_{\o v}\times\rit))\times\Aut_k(\phi _L(F_{v}\times\rit))$
of shifted compact transcendental biquanta.\epr
\Ee
\vskip 11pt

\subsection{The third question ``What interpretation can we give to the local spacings between the eigenvalues of large random matrices?''}

This question \cite{A-B-K}, \cite{Joh} is a direct consequence of the responses given to the two first questions and shows  the central importance of the biquanta in this field as proved in the following sections.
\vskip 11pt

\subsection{Proposition}

{\em
\Bena
\item The consecutive spacings
\[ 
\delta E\RL(j) = E\RL(j+1)-E\RL(j)\;, \qquad 1\le j\le r< \infty \;, \]
between the eigenbivalues of the random matrix $G$ of
$\GL_r(\rit\times\rit)$ are the infinitesimal bigenerators on one biquantum of the Lie subbisemialgebra
$\gl_1(F_{\o v^1}\times F_{v^1})$ of the bilinear parabolic unitary semigroup
$P_1(F_{\o v^1}\times F_{v^1})\subset\GL_1(F_{\o v}\times F_{v})$.

\item The $k$-th consecutive spacings
\[ 
\delta E^{(k)}\RL(j) = E\RL(j+k)-E\RL(j)\;, \qquad k\le j\;, \]
between the eigenbivalues of the random matrix $G$ of $\GL_r(\rit\times\rit)$ are the infinitesimal bigenerators on $k$ biquanta of the Lie subbisemialgebra
$\gl_1(F_{\o v^k}\times F_{v^k})$ of the bilinear $k$-th semigroup
$\GL_1(F_{\o v^k}\times F_{v^k})\subset \GL_1(F_{\o v}\times F_{v})$,
where $F_{v^k}$ (and $F_{\o v^k}$) denotes the set of transcendental subextensions characterized by a transcendence degree
\[ \tr\cdot d\cdot F_{v^k}=k\cdot N\]
referring to $k$ biquanta included in transcendental extensions having higher or equal transcendence degree.
\Ee}
\vskip 11pt

\bpr
Referring to proposition~3.4 and corollary~3.5, it appears that the set
$\{E\RL(j)\}_{j=1}^r$ of eigenvalues of a random matrix $G$ of $\GL_r(\rit\times\rit)$, being shifts of eigenbifunctions of the eigenbivalue equation
\[
(D_R\otimes D_L)(\phi _r(G_{(Q)}(F_{\o v^1}\times F_{v^1})))
=E\RL(\phi _r(G_{(Q)}(F_{\o v^1}\times F_{v^1})))\]
is the set of infinitesimal bigenerators on $j$ biquanta of the Lie bisemialgebra
$\gl_1(F_{\o v}\times F_{v})$ of the Lie bilinear semigroup
$\GL_1(F_{\o v}\times F_{v})$.

And, thus, the consecutive spacings:
\[ 
\delta E\RL(j) = E\RL(j+1)-E\RL(j)\;,  \]
between the eigenbivalues of a random matrix $G$ of $\GL_r(\rit\times\rit)$ are the infinitesimal bigenerators on the biquantum on the Lie subbisemialgebra
$\gl_1(F_{\o v^1}\times F_{v^1})$ of the bilinear parabolic unitary semigroup
$P_1(F_{\o v^1}\times F_{v^1})$ on the product, right by left, of sets
$F_{\o v^1}$ and $F_{v^1}$ of transcendental subextensions characterized by a transcendence degree $\tr\cdot d\cdot F_{v^1}=N$.

Similarly, the $k$-th consecutive spacings
\[ 
\delta E^{(k)}\RL(j) = E\RL(j+k)-E\RL(j)  \]
refer to the infinitesimal bigenerators on $k$ biquanta of the Lie subbisemiagebra
$\gl_1(F_{\o v^k}\times F_{v^k})$ of the Lie bilinear $k$-th semigroup
$\GL_1(F_{\o v^k}\times F_{v^k})\subset\GL_1(F_{\o v}\times F_{v})$.\epr
\vskip 11pt

\subsection{Corollary}

{\em
\Bena
\item The consecutive spacings
\[ 
\delta E\RL(j) = E\RL(j+1)-E\RL(j)\;, \qquad 1\le j\le r< \infty \;, \]
between the eigenbivalues of the random matrix $G$ of
$\GL_r(\rit\times\rit)$ correspond to the energies of one free biquantum from subbisemilattices of $(j+1)$ biquanta.

\item The $k$-th consecutive spacings
\[ 
\delta E^{(k)}\RL(j) = E\RL(j+k)-E\RL(j)\;, \qquad k\le j\;, \]
between the eigenbivalues of the random matrix $G$ of $\GL_r(\rit\times\rit)$ 
correspond to the energies of $k$ free biquanta from subbisemilattices of $(j+k)$ biquanta.
\Ee}
\vskip 11pt

\bpr
As the infinitesimal bigenerators
$E\RL(j)$ of the Lie bisemiagebra
$\gl_1(F_{\o v}\times F_{v})$ are shifts of eigenbifunctions on $j$ transcendental biquanta according to section~3.3, they correspond to the energies of these $j$ transcendental biquanta and thus to $j$ transcendental biquanta of energy.\epr
\vskip 11pt

\subsection{Lemma}

{\em 
Let $\delta E\RL(j)$ denote the consecutive spacings between the eigenbivalues of the matrix $G$ of $\GL_r(\rit\times\rit)$.

Then, $\delta E\RL(j)$ decomposes into:
\[
\delta E\RL(j)= \delta EF\RL(j) +  \delta EV\RL(j)\]
where $\delta EF\RL(j)$ and $\delta EV\RL(j)$ denote respectively the fixed (or constant) and variable consecutive spacings between the $r$ eigenbivalues of the matrix $G$ of $\GL_r(\rit\times\rit)$.
}
\vskip 11pt

\bpr
There are surjective maps
\[
\delta E_{\to v}(j): \qquad \delta E\RL(j) \quad \To \quad
\delta EV\RL(j)\]
between the consecutive spacings
$\delta E\RL(j)$, referring to the infinitesimal generators of Lie subbisemialgebras or energies of one biquantum in a lattice of $(j+1)$ biquanta, and the respective ``variable'' consecutive spacing
$\delta EV\RL(j)$ in such a way that their kernels are the fixed consecutive spacings
$\delta EF\RL(j)$ being equal for every integer $j$, $1\le j\le r<\infty $.

Consequently, the energy
$\delta E\RL(j)$ of one biquantum in a lattice of $(j+1)$ biquanta decomposes into a fixed part common to all considered bisections and into a variable part differing from one bisection to another.\epr
\vskip 11pt

\subsection{Proposition}

{\em
\Bena
\item The consecutive spacings
\[
\delta EV^{\rm BCUE}\RL(j)= E^{\rm BCUE}\RL(j+1)- E^{\rm BCUE}\RL(j)\;, \qquad 1\le j\le r\;, \]
between the eigenbivalues
$E^{\rm BCUE}\RL(j+1)$ and
$E^{\rm BCUE}\RL(j)$ of the unitary random matrix
$u_r(\cit\times\cit)\in U_r(\cit\times\cit)$ (or
$o_r(\rit\times\rit)\in O_r(\rit\times\rit)$) are the variable (unitary) infinitesimal bigenerators on one biquantum of the envisaged Lie subbisemialgebra or variable (unitary) energies
$\delta EV^{\rm BCUE}\RL(j)$ of one biquantum in subbisemilattices of $(j+1)$ biquanta.

\item The $k$-th consecutive spacings
\[
\delta EV^{(k)\rm BCUE}\RL(j)= E^{\rm BCUE}\RL(j+k)- E^{\rm BCUE}\RL(j)\;, \qquad k\le j\;, \]
between the eigenbivalues
of the unitary random matrix
$ u_r(\cit\times\cit)$ (or
$o_r(\rit\times\rit)$) are the variable (unitary) infinitesimal bigenerators on $k$ biquanta of the envisaged Lie subbisemialgebra or variable (unitary) energies
$\delta EV^{(k)\rm BCUE}\RL(j)$ on $k$ biquanta in subbisemilattices of $(j+k)$ biquanta.
\Ee}
\vskip 11pt

\bpr
Similarly as it was developed in proposition~3.18, the set
$\{E^{\rm BCUE}\RL(j)\}_{j=1}^r$ of eigenbivalues of $u_r(\cit\times\cit)$, being shifts of eigenbifunctions of the eigenbivalue equation:
\[
(D_R\otimes D_L)(\phi _r(P_{(Q)}(F_{\o v^1}\times F_{v^1})))=
E^{\rm BCUE}\RL(j)(\phi _r(P_{(Q)}(F_{\o v^1}\times F_{v^1})))\]
is the set of unitary infinitesimal bigenerators of the Lie bisemialgebra
$\gl_1(F_{\o v^1}\times F_{v^1}))$ of the Lie bilinear parabolic semigroup
$P_1(F_{\o v^1}\times F_{v^1})  $.

Indeed, the bilinear semigroup of matrices $U_r(\cit\times\cit)$ (or
$O_r(\rit\times\rit)$), of the bilinear circular unitary (or orthogonal) ensemble
BCUE \resp{BCOE}, constitutes the representation of the unitary bilinear differential Galois semigroup associated with the biaction of the differential bioperator
$(D_R\otimes D_L)$ on the unitary bisemisheaf
$(\wh M^{(1)}_{v^1_R}\otimes \wh M^{(1)}_{v^1_R})\subset
(\wh M^{(1)}_{v_R}\otimes \wh M^{(1)}_{v_R})$ according to proposition~2.4 to corollary~2.7.

And, thus, the consecutive spacings
\[
 E^{\rm BCUE}\RL(j+1)- E^{\rm BGUE}\RL(j)=
\delta EV^{\rm BCUE}\RL(j)\]
between the eigenbivalues of $U_r(\cit\times\cit) $ (or $O_r(\rit\times\rit)$) are the variable (unitary) infinitesimal bigenerators of Lie subbisemialgebras $\gl_1(F_{\o v^1}\times F_{v^1})$ or (unitary) variable energies
$\delta EV^{\rm BCUE}\RL(j)$ of one biquantum in sublattices of $(j+1)$ biquanta.

The case of $k$-th consecutive spacings
$\delta EV^{(k)\rm BCUE}\RL(j)$ can be handled similarly by taking into account the content of proposition~3.18.

Remark finally that
\be
\delta EV\RL(j)=
\delta EV^{\rm BCUE}\RL(j)\equiv
\delta E^{\rm BCUE}\RL(j)\;.\tag*{\eop}\ee
\vskip 11pt

\subsection{The fourth question ``What interpretation can we give to the spacings between the nontrivial zeros of the zeta function $\zeta (s)$?''}

This question depends on the solution of the Riemann hypothesis proposed in \cite{Pie7} and is briefly recalled in the following sections \cite{Edw}, \cite{Rie}, \cite{Tit}.
\vskip 11pt

\subsection{Cuspidal representation given by global elliptic bisemimodule}

Let $S_L(2,N=1)$ \resp{$S_R(2,N=1)$} denote the (semi)algebra of cusp forms $f_L(z)$ \resp{ $f_R(z)$} of weight $2$ and level $N=1$ holomorphic in the upper \resp{lower} half plane as developed in \cite{Pie7}.
$f_L(z)$ \resp{$f_R(z^*)$}, expanded in Fourier series according to:
\[ f_L(z)=\sum_n\ a_n\ q^{n}\ ,\; q=e^{2\pi iz}\ , \; z\in \cit\;, \quad
\rresp{f_R(z^*)=\sum_n\ a_n\ q^{*n}\ ,\; q^*=e^{-2\pi iz^*}}\]
is the functional representation space of $G^{(2)}(F^T_\omega )\equiv T_2(F^T_\omega )$ \resp{$G^{(2)}(F^T_{\o \omega} )\equiv T^t_2(F^T_{\o\omega })$}
where $F^T_\omega $ \resp{$F^T_{\o\omega }$} is the set of increasing toroidal complex transcendental extensions referring to section~1.9.

Then, we have that:
\begin{align*}
f_R(z^*)\times f_L(z)
&= \FREPSP (\GL_2(F^T_{\o\omega }\times F^T_\omega )\\
&= G^{(2)}(F^T_{\o\omega }\times F^T_\omega )\\
&= M^{(2)}_{\omega ^T_R}\otimes M^{(2)}_{\omega ^T_L}
\end{align*}
is a cusp biform of weight $2$ and level $1$.

On the other hand, we can consider the map:
\[ \Ms_{f\to \zeta _L}: \quad f_L(z)\To \zeta _L(s_+)
\rresp{ \Ms_{f\to \zeta _R}: \quad f_R(z^*)\To \zeta _R(s_-)}
\]
of the cusp form $f_L(z)$ \resp{$f_R(z^*)$} into the corresponding zeta function
$\zeta _L(s_+)$ \resp{$\zeta _R(s_-)$} in such a way that
$s_+=\sigma +i\tau $
\resp{$s_-=\sigma -i\tau $} be an ``energy'' variable conjugate to the spatial variable $z$ \resp{$z^*$}.

Referring to section~2.12 introducing a $2n$-dimensional global elliptic
$(\Gamma _{\wh M^{(2n)}_{v^T\RL}})$-bisemi\-module
$\phi ^{(2n)}\RL(x)$, we see that it can be reduced to {\bbf a $1D$-pseudounramified simple global elliptic
$(\Gamma _{\wh M^{(2n)}_{v^T\RL}})$-bisemimodule\/}:
\[
\phi ^{(1)}\RL(x) = \sum_n (\lambda ^{(nr)}(n)\ e^{-2\pi inx}\otimes_D
\lambda ^{(nr)}(n)\ e^{+2\pi inx})\;, \qquad x\in\rit\;,\]
if:
\Bean
\item pseudounramification is concerned, i.e.  the conductor $N=1$;

\item simplicity is supposed, i.e. the multiplicity $m_n$ is equal to one on each level ``$n$''.
\Ee
This global elliptic
$(\Gamma _{\wh M^{(1)}_{v^T\RL}})$-bisemimodule
\[ \phi ^{(1),nr}\RL(x)=
 \phi ^{(1),nr}_R(x)\otimes_D
 \phi ^{(1),nr}_L(x) \qquad (\ \otimes_D \text{\ : diagonal tensor product})\]
 can be interpreted geometrically as {\bbf the sum of products, right by left, of semicircles of level ``$n$'' on $n$ transcendental compact quanta\/} in such a way that
$\phi ^{(1),nr}\RL(x)$ be the cuspidal automorphic representation of the complete bilinear semigroup
$\GL_2(F^{(nr)}_{\o v}\times F^{(nr)}_{v})$ according to sections~1.5 and 2.12 and cover the weight $2$ cusp biform $f_R(z^*)\otimes f_L(z)$.

On the other hand, as we are dealing with bisemiobjects, the zeta functions
$\zeta _R(s_-)$ and
$\zeta _L(s_+)$,  defined respectively on the lower and upper half planes, are considered.
\vskip 11pt

\subsection{Proposition}

{\em
Let
\[ H_{\phi \RL\to\zeta \RL} : \qquad
2\phi ^{(1),nr}_R(x)\otimes_D
 \phi ^{(1),nr}_L(x)\quad \To \quad \zeta _R(s_-)\otimes_D \zeta _L(s_+)\]
 be the map between the double $1D$-pseudounramified simple global elliptic bisemimodule
$ \phi ^{(1),nr}\RL(x)$ and the product, right by left, of zeta functions given by:
\[
\zeta _R(s_-)\otimes_D \zeta _L(s_+)=\sum_n\L(n^{-s_-}\otimes_D n^{-s_+}\R)\;.\]

Then, {\bbf the kernel $\Ker(H_{\phi \RL\to\zeta \RL})$ of the map
$H_{\phi \RL\to\zeta \RL}$ is the set of squares of trivial zeros of $\zeta _R(s_-)$ and
$\zeta _L(s_+)$\/}.
}
\vskip 11pt

\bpr The kernel
$\Ker(H_{\phi \RL\to\zeta \RL})$ of 
$H_{\phi \RL\to\zeta \RL}$ maps into the set  of ``trivial'' zeros of $\zeta _R(s_-)$ and
$\zeta _L(s_+)$ which are the negative integers $-2$, $-4$, \ldots

Consequently, {\bbf
this kernel must be the set of bipoints\/}:
\begin{align*}
\{ \sigma _{n_R}\times \sigma _{n_L}
&= 2\lambda ^{(nr)}(n) (e^{-2\pi inx}\mid x=0)\times 
 2\lambda ^{(nr)}(n) (e^{-2\pi inx}\mid x=0)\}\\
 &= 4(\lambda ^{(nr)}(n^2))^2\ \}\\
 &= 4f^2_{v_n}=4n^2\ \}\end{align*}
 where $(\lambda ^{(nr)}(n^2))^2$ is the square of the global residue degree
 $f_{v_n}=n$ as proved in \cite{Pie7},
 
 {\bbf in such a way that the \lr point
 $\sigma _{n_L}=2\lambda ^{(nr)}(n)\ e^{2\pi inx}\mid x=0)$
 \resp{$\sigma _{n_R}=2\lambda ^{(nr)}(n)\ e^{-2\pi inx}\mid x=0)$}
 corresponds to the degeneracy of the (irreducible) circle
 $2\lambda ^{(nr)}(n)\ e^{2\pi inx}$
 \resp{$2\lambda ^{(nr)}(n)\ e^{-2\pi inx}$} on $2n$ \lr transcendental quanta\/}.
 
 Remark that the one-to-one correspondence between the global elliptic semimodule\linebreak
 $2\phi ^{(1),nr}_L(x)$
 \resp{$2\phi ^{(1),nr}_R(x)$} and the \lr zeta function
 $\zeta _L(s_+)=\sum\limits_n n^{-s_+}$
 \resp{$\zeta _R(s_-)=\sum\limits_n n^{-s_-}$} is clear if it is noted that
 $n^s=n^{\sigma +i\tau }$ can be written
 \[ e^{(\sigma +i\tau )\ell nn}= e^{\sigma \ell nn}\cdot
 e^{i\tau \ell nn}\]
 where
$ e^{\sigma  \ell nn}$ can correspond to the radius of a circle with phase
$e^{i\tau \ell nn}$ in such a way that $1/e^{\sigma \ell nn}$ maps into
$2\lambda ^{(nr)}(n)$ and
$1/e^{i\tau  \ell nn}$ maps into
$e^{2\pi inx}$ for each term of $\zeta _L(s_+)$ and of $2\phi ^{(1),nr}_L(x)$.\epr
\vskip 11pt

\subsection{Proposition}

{\em
Let $D_{4n^2,i^2}\cdot \varepsilon _{4n^2}$ be a coset representative of the Lie (bisemi)algebra of the decomposition (bisemi)group
$D_{i^2}(\zit)_{|4n^2}$ and let $\alpha _{4n^2}$ be the split Cartan subgroup (bi)element.

Then, {\bbf
the products of the pairs of the trivial zeros of the Riemann zeta functions
$\zeta _R(s_-)$ and $\zeta _L(s_+)$ are mapped into the products of the corresponding pairs of the nontrivial zeros according to:\/}:
\begin{multline*}
\begin{aligned}[t]
\{D_{4n^2,i^2}\cdot \varepsilon _{4n^2}\}: \qquad
\{\det (\alpha _{4n^2})\}_n \quad &\To \quad
\{\det (D_{4n^2,i^2}\cdot \varepsilon _{4n^2}\cdot \alpha _{4n^2})_{ss}\}_n \\
\{(-2n)\times(-2n)\}_n \quad &\To \quad
\{ \lambda ^{(nr)}_+(4n^2,i^2,E_{4n^2})\times
\lambda ^{(nr)}_-(4n^2,i^2,E_{4n^2})\}_n\;, 
\end{aligned}\\
\forall\ n\in\nit\end{multline*}
where ``$ss$ denotes the semisimple form.
}
\vskip 11pt

\bpr
Let
{\bbf
\[ \alpha _{4n^2}=\BM 4n^2 & 0\\ 0&1\EM\]
be the (split) Cartan subgroup element\/} associated with the square of the global residue degree $f_{v_{2n}}=2n$.

Let
{\bbf
\[ D_{4n^2,i^2}=\BM 1&i \\ 0&1\EM \ \BM 1&0 \\ i&1\EM \]
be the coset representative of the Lie (bisemi)algebra} \cite{Pie3} {\bbf
of the decomposition (bisemi)group acting on $\alpha _{4n^2}$\/}:

It corresponds to the coset representative of an unipotent Lie (bisemi)algebra mapping in the topological Lie (bisemi)algebra 
$\gl_1(F^{(nr)}_{\o v}\times F^{(nr)}_{v})$ consisting in vector fields on the Lie subgroup
$\GL_1(F^{(nr)}_{\o v}\times F^{(nr)}_{v})$.

Let
{\bbf
\[ \varepsilon _{4n^2}=\BM E_{4n^2} & 0\\ 0&1\EM\]
be the infinitesimal (bi)generator of this Lie (bisemi)algebra} corresponding to the square of the global residue degree $f_{v_{2n}}=2n$.
\vskip 11pt

{\bbf Every root of the Lie (bisemi)algebra is determined by the (equivalent) eigenvalues
$\lambda ^{(nr)}_\pm(4n^2,i^2,E_{4n^2})$ of
\[ 
D_{4n^2,i^2}\cdot \varepsilon _{4n^2}\cdot\alpha _{4n^2}=\L[
\BM 1&i \\ 0&1\EM \ \BM 1&0\\ i&1\EM\R] \BM E_{4n^2} &0\\ 0&1\EM\ \BM {4n^2} &0\\ 0&1\EM
\]
given by
\[
\lambda ^{(nr)}_\pm(4n^2,i^2,E_{4n^2})=\frac{1\pm i\sqrt{(16n^2\cdot E_{4n^2})-1}}2\]
}
where $D_{4n^2,i^2}\cdot \varepsilon _{4n^2}$ is the coset representative of the Lie
(bisemi)algebra ${\rm Lie}(D_{4n^2}(\zit)_{|4n^2})$ of the decomposition group
$D_{i^2}(\zit)_{|4n^2}$.
\vskip 11pt

According to proposition~3.24, the squares $(-2n)^2$ of the trivial zeros of
$\zeta _R(s_-)$ and $\zeta _L(s_+)$ are the squares of the global residue degrees
$f_{v_{2n}}=2n$.

As $D_{4n^2}\cdot \varepsilon _{4n^2}$ is of Galois type, it maps squares of trivial zeros $(-2n)^2$ into products of corresponding pairs of other zeros
$\lambda ^{(nr)}_+(4n^2,i^2,E_{4n^2})=
\lambda ^{(nr)}_-(4n^2,i^2,E_{4n^2})$ which are nontrivial zeros since they have real parts localized on the line $\sigma =\half$.

So, $D_{4n^2,i^2}\cdot \epsilon _{4n^2}$ maps a lattice of transcendental biquanta on the considered Lie (bisemi) group into a lattice of energies of these biquanta on the associated Lie (bisemi)algebra.\epr
\vskip 11pt

\subsection{Corollary}

{\em
The eigenvalues 
$\lambda ^{(nr)}_+(4n^2,i^2,E_{4n^2})$ and 
$\lambda ^{(nr)}_-(4n^2,i^2,E_{4n^2})$ of
$(D_{4n^2,i^2}\cdot \varepsilon _{4n^2}\cdot\alpha _{4n^2})$ for all $n\in\nit$ are the nontrivial zeros of the Riemann zeta function
$\zeta (s)=\sum\limits_n n^{-s}$.
}
\vskip 11pt

\bpr
The set $\{-2n\}_n$, being the trivial zeros of the right and left zeta functions
$\zeta _R(s_-)$ and $\zeta _L(s_+)$, constitutes also the set of trivial zeros of the classical Riemann zeta function $\zeta (s)=\sum\limits_n n^{-s}$.

So, the eigenvalues
$\lambda ^{(nr)}_+(4n^2,i^2,E_{4n^2})$ and 
$\lambda ^{(nr)}_-(4n^2,i^2,E_{4n^2})$ of
$(D_{4n^2,i^2}\cdot \varepsilon _{4n^2}\cdot\alpha _{4n^2})$, generated from the corresponding trivial zeros ``$-2n$'', are:
\Bean
\item  the nontrivial zeros of $\zeta (s)$ since they are localized on the line $\sigma =\half$ and disposed symmetrically on this line with respect to $\tau =0$, $s=\sigma +i\tau $;
\item the relevant zeros of the energy density function $\zeta (s)$ which is an inverse space function, i.e. an ``energy'' function, on transcendental quanta labelled by the integers ``$n$''.\epr
\Ee
\vskip 11pt

\subsection{Proposition}

{\em
Let the nontrivial zeros of $\zeta (s)$, $\lambda ^{(nr)}_+(4n^2,i^2,E_{4n^2})$ and 
$\lambda ^{(nr)}_-(4n^2,i^2,E_{4n^2})$ be written compactly and classically according to
$\half+i\gamma _j$ and 
$\half-i\gamma _j$, $j\hookleftarrow n$.

\Bena
\item Then, {\bbf 
the consecutive spacings
\[ 
\delta \gamma _j=\gamma _{j+1}-\gamma_j\;, \qquad j=1,2,\dots\]
between the nontrivial zeros of $\zeta (s)$ are equivalently:
\Be
\item the infinitesimal generators on one quantum of the Lie subsemialgebra
$\gl_1(F^{(nr)}_{v^1})$ (or
$\gl_1(F^{(nr)}_{\o v^1})$) of the linear parabolic unitary semigroup
$P_1(F^{(nr)}_{v^1})\subset \GL_1(F^{(nr)}_{v})\equiv T_1(F^{(nr)}_{v})$
(or $P_1(F^{(nr)}_{\o v^1})\subset \GL_1(F^{(nr)}_{\o v})\equiv T^t_1(F^{(nr)}_{\o v})$).

\item the energies of one (free) quantum from subsemilattice of $(j+1)$ quanta
in such a way that $\delta _{\gamma _j}>\F{\gamma _{j+1}}{j+1}$ where $\F{\gamma _{j+1}}{j+1}$ is the energy of one quantum bound to $j$ quanta.
\Ee

\item {\bbf 
The $k$-th consecutive spacings
\[ 
\delta^{(k)} _j=\gamma _{j+k}-\gamma_j\;, \]
between the nontrivial zeros of $\zeta (s)$ are equivalently:
\Be
\item the infinitesimal generators on $k$ quanta of the Lie subsemialgebra
$\gl_1(F^{(nr)}_{v^k})$ (or
$\gl_1(F^{(nr)}_{\o v^k})$) of the $k$-th  semigroup
$ \GL_1(F^{(nr)}_{v^k})\subset \GL_1(F^{(nr)}_{v})$
(or $ \GL_1(F^{(nr)}_{\o v^k})\subset \GL_1(F^{(nr)}_{\o v})$).

\item the energies of $k$ (free) quanta from subsemilattice of $(j+k)$ quanta.
\Ee}
}
\Ee}
\vskip 11pt

\bpr
According to proposition~3.25, every nontrivial zero
$\L(\half\pm\gamma _{j+1}\R)$
 of $\zeta _R(s_-)$, $\zeta _L(s_+)$ or $\zeta (s)$ is the infinitesimal generator of the Lie semialgebra 
 $\gl_1(F^{(nr)}_{v_{j+1}})$ (or
 $\gl_1(F^{(nr)}_{\o v_{j+1}})$) on $(j+1)$ compact transcendental quanta or, physically, the energy of $(j+1)$ compact transcendental quanta.
 
 So, the consecutive spacing
 \[ 
\delta \gamma _j=\gamma _{j+k}-\gamma_j\;, \]
between $\gamma _{j+1}$ and $\gamma _j$, is the infinitesimal generator on one quantum on the Lie subsemialgebra  $\gl_1(F^{(nr)}_{v^1})$ of the parabolic unitary semigroup
 $P_1(F^{(nr)}_{v^1})$ or the energy of one free quantum in a subsemilattice of $(j+1)$ transcendental compact quanta.
 
 The $k$-th consecutive spacings 
\[ 
\delta^{(k)} _j=\gamma _{j+k}-\gamma_j\]
can be handled similarly.\epr
\vskip 11pt

\subsection{The fifth question ``What is the curious connection between the spacings of the nontrivial zeros of $\zeta (s)$ and the corresponding spacings between the eigenvalues of random matrices?''}

This question finds a response in the following propositions \cite{B-K}, \cite{G-M}.
\vskip 11pt

\subsection{Proposition}

{\em
The consecutive spacings
\[ 
\delta \gamma _j=\gamma _{j+1}-\gamma_j\;, \qquad j\in\nit\;,\]
between the nontrivial zeros of the Riemann zeta function $\zeta (s)$ as well as the consecutive spacings
\[ 
\delta E^{(nr)}\Rl(j)= E^{(nr)}\Rl(j+1)-E^{(nr)}\Rl(j)\]
between the square roots of the pseudounramified eigenbivalues of a random matrix (of) $\GL_r(\rit\times\rit)$ (or of
$\GL_r(\cit\times\cit)$) or between the pseudounramified eigenvalues of a random matrix of $\GL_r(\rit)$, are equivalently:
\Bean
\item the infinitesimal generators on one quantum of the Lie subsemialgebra
$\gl_1(F^{(nr)}_{v^1})$ (or 
$\gl_1(F^{(nr)}_{\o v^1})$) of the linear parabolic unitary semigroup
$P_1(F^{(nr)}_{v^1})$ (or 
$P_1(F^{(nr)}_{\o v^1})$);

\item the energies of one transcendental pseudounramified $(N=1)$ quantum in subsemilattices of $(j+1)$ transcendental pseudounramified quanta.
\Ee
( $R,L$ means $R$ (right) or $L$ (left)).
}
\vskip 11pt

\bpr
The equivalence between the consecutive spacings $\delta \gamma _j$ of $\zeta (s)$ and the consecutive spacings $\delta E^{(nr)}\Rl(j)$ between the square roots of pseudounramified eigenbivalues of a random matrix of 
$\GL_r(\rit\times\rit)$ or between the pseudounramified eigenvalues of a random matrix of $\GL_r(\rit)$, results from propositions~3.18 and 3.27.

The equivalence is exact if the eigenbivalues of a matrix of $\GL_r(\rit\times\rit)$ are pseudounramified, i.e. if they are eigenbivalues of the eigenbivalue equation (see proposition~3.18):
\[ (D_R\otimes D_L)
(\phi _r(G_{(Q)}(F^{(nr)}_{\o v^1}\times F^{(nr)}_{v^1})))=E^{(nr)}\RL(j)
(\phi _r(G_{(Q)}(F^{(nr)}_{\o v^1}\times F^{(nr)}_{v^1})))\]
where $F^{(nr)}_{v^1}$ (and  $F^{(nr)}_{\o v^1}$) are unitary transcendental pseudounramified extensions (case $N=1$) referring to section~1.5.

Finally, remark that
\[ 
|\delta E^{(nr)}\Rl| = \L|\sqrt{\delta E^{(nr)}\RL(j)}\R|\]
where
$\delta E^{(nr)}\Rl(j)$ is a consecutive spacing between pseudounramified eigenvalues and\linebreak 
$\sqrt{\delta E^{(nr)}\RL(j)}$ is a consecutive spacing between square roots of pseudounramified eigenbivalues.\epr
\vskip 11pt

\subsection{Corollary}

{\em
The set $\{\delta E^{(nr)}\Rl(j)\}_j$ of consecutive spacings between the square roots of the eigenbivalues of a random matrix of 
$\GL_r(\rit\times\rit)$ (or of
$\GL_r(\cit\times\cit)$) or between the eigenvalues of a random matrix of $\GL_r(\rit)$, as well as the set $\{\delta \gamma _j\}_j$ of  consecutive spacings between the nontrivial zeros of $\zeta (s)$ constitutes a representation of the differential inertia Galois (semi)group associated with the action of the differential operator $D_L$ or $D_R$.
}
\vskip 11pt

\bpr
This is immediate since the sets $\{\delta E^{(nr)}\Rl(j)\}_j$ and 
$\{\delta \gamma _j\}_j$ are infinitesimal generators of the Lie subsemialgebra of the linear parabolic unitary semigroup
$P_1(F^{(nr)}_{v^1})$ (or 
$P_1(F^{(nr)}_{\o v^1})$) and as the unitary parabolic bilinear semigroup
$P_r(\rit\times\rit)\subset \GL_r(\rit\times\rit)$ corresponds to the ``bilinear'' representation of the product, right by left, 
$\Int_k(\phi _R(\rit))\times\Int_k(\phi _L(\rit))$ of differential inertia Galois semigroups according to proposition~2.4 and corollaries~2.5 and 3.5, we have that:
\be P_r(\rit\times\rit)=\Rep
[(\Int_k(\phi _R(\rit )) _{|F_{\o v^1_r}})\times(\Int_k(\phi _L(\rit )) _{|F_{\o v^1_r}})]\;.\tag*{\eop}\ee

\subsection{Proposition}

{\em
Let $\{\delta \gamma _j\}_j$ be the set of consecutive spacings between the nontrivial zeros of $\zeta (s)$ and let
$\{\delta E\Rl(j)\}_j$ be the set of consecutive spacings between the square roots of the eigenbivalues of a random matrix of
$\GL_r(\rit\times\rit)$ (or of
$\GL_r(\cit\times\cit)$) or between the eigenbivalues of a random matrix of
$\GL_r(\rit)$.

Then, {\bbf
there is a surjective map:
\[ IM_{E\to \gamma }: \qquad \{\delta E\Rl(j)\}_j \quad \To \quad \{\delta \gamma _j\}_j\]
of which kernel $\Ker[IM_{E\to \gamma }]$ is equivalently the set
\[ 
\{ \delta E\Rl(j)-\delta E^{(nr)}\Rl(j) \}_j\;, \qquad \forall\ j\ , \quad 1\le j\le r\;, \]
\Bean
\item of differences of consecutive spacings between the square roots of the pseudoramified and pseudounramified eigenbivalues of a random matrix of
$\GL_r(\rit\times\rit)$)(or of
$\GL_r(\cit\times\cit)$), or between the pseudoramified and pseudounramified eigenvalues of a random matrix of $\GL_r(\rit)$;

\item of the energies of one compact transcendental pseudoramified $(N>2)$ quantum in subsemilattices of $(j+1)$ transcendental pseudoramified  quanta.
\Ee
}}
\vskip 11pt

\bpr
First, remark that the surjective map
\[ IM_{E\to \gamma }: \qquad \{\delta E\Rl(j)\}_j \quad \To \quad \{\delta \gamma _j\}_j\]
leads precisely to the map (or to the equality) between the spacing distribution between eigenvalues of a random matrix and the pair correlation of the nontrivial zeros of $\zeta (s)$ \cite{Mon}.
\vskip 11pt

If the kernel $\Ker[IM_{E\to \gamma }]$ of the map 
$IM_{E\to \gamma }$ is null, then
\[
\delta E\Rl(j) = \delta E^{(nr)}\Rl(j)\;,\qquad j\in\nit\;, \]
i.e. the consecutive spacings between the eigen(bi)values of a random matrix are pseudounramified, implies the thesis of proposition~3.29, i.e. the equality between the consecutive spacings $\delta \gamma _j$ of $\zeta (s)$ and the consecutive spacings
$\delta E^{(nr)}\Rl(j)$ between eigenvalues.\epr
\vskip 11pt

\subsection{Proposition}

{\vskip 11pt\em
Let
\begin{align*}
\delta EV^{(nr),\rm BCOE}\Rl(j)
&= E^{(nr),\rm BCOE}\Rl(j+1)-E^{(nr),\rm BCOE}\Rl(j)\\
&= \delta E^{(nr),\rm BCOE}\Rl(j)\;, && j\in\nit\;, \end{align*}
be the consecutive spacings between the square roots of the pseudounramified eigenbivalues of a random unitary matrix of
$O_r(\rit\times\rit)$ or of $U_r(\cit\times\cit)$) or between the pseudounramified eigenvalues of a random unitary matrix of $O_r(\rit)$.

Then, there is a surjective map:
\[
IM_{\gamma \to E^{(nr)}_{\rm BCOE}}: \qquad
\{\delta \gamma _j\}_j \quad \To \quad \{\delta EV^{(nr),\rm BCOE}\Rl(j)\}_j\]
between the set $\{\delta \gamma _j\}_j$ of consecutive spacings between the nontrivial zeros of $\zeta (s)$ and the set
$\{\delta EV^{(nr),\rm BCOE}\Rl(j)\}_j$.
}
\vskip 11pt

\bpr
According to lemma~3.20, the consecutive spacings
$\delta E\RL(j)$ between the eigenbivalues of $\GL_r(\rit\times\rit)$ decomposes into fixed and variable consecutive spacings
\[ \delta E\RL(j)=\delta EF\RL(j)+\delta EV\RL(j)\;, \]
which is also the case for the consecutive spacings
$\delta E^{(nr)}\RL(j)$ between pseudounramified eigenbivalues.

As
\[ \delta \gamma _j=\delta  E^{(nr)}\Rl(j)\]
 according to proposition~3.29, it results that the consecutive spacings between nontrivial zeros of $\zeta (s)$ also decomposes according to fixed and variable consecutive spacings
\[
\delta \gamma _j=\delta \gamma F_j+\delta \gamma V_j\;, \qquad \forall\ j\;, \quad 1\le j\le r\;, \]
in such a way that
\[
\delta \gamma V_j=\delta EV\Rl^{(nr),\rm BCOE}(j)\equiv\delta E\Rl^{(nr),\rm BCOE}(j)\;.\]
And, the variable consecutive spacings
$\{\delta \gamma V_j\}_j$ constitute a representation of the differential variable inertia Galois subgroup according to corollary~3.30, {\bbf
the differential variable inertia Galois subgroup being a subgroup of the differential inertia Galois subgroup\/}.\epr
\vskip 11pt

\subsection{Proposition (main)}

{\em
{\bbf
Let
\[ \delta \gamma _j^{(k)}=\gamma _{j+k}-\gamma _j\]
denote the $k$-th consecutive spacings between the nontrivial zeros of $\zeta (s)$.

Let:
\Bi
\item $\delta E\Rl^{(k)}(j)=E\Rl(j+k)-E\Rl(j)$;
\item $\delta E\Rl^{(nr),(k)}(j)=E^{(nr)}\Rl(j+k)-E^{(nr)}\Rl(j)$;

and

\item $\delta EV^{(k),(nr),\rm BCOE}\Rl(j)=E^{(nr),\rm BCOE}\Rl(j+k)-E^{(nr),\rm BCOE}\Rl(j)$, 

$1\le j\le r$, $k\le j$,
\Ei
be the $k$-th consecutive spacings between respectively}
\Bi 
\item the pseudoramified eigenvalues of a random matrix of $\GL_r(\rit)$;
\item the pseudounramified eigenvalues of a random matrix of $\GL_r(\rit)$;
\item the pseudounramified eigenvalues of a random unitary matrix of $O_r(\rit)$.
\Ei

Then, we have:
\Bena
\item {\bbf
$\delta \gamma _j^{(k)}=\delta E^{(nr),(k)}\Rl(j)$ which are equivalently:
\Be
\item the infinitesimal generators on $k$ quanta of the Lie subsemi\-algebra
$\gl_1(F^{(nr)}_{v^k})$ of the linear $k$-th semigroup
$\GL_1(F^{(nr)}_{v^k})\subset 
\GL_1(F^{(nr)}_{v})$;

\item the energies of $k$ transcendental pseudounramified $(N=1)$ quanta in subsemilattices in $(j+k)$ transcendental pseudounramified quanta;

\item a representation of the differential Galois (semi)group associated with the action of the differential operator $D_L$ or $D_R$ on a function on $k$ transcendental pseudounramified quanta.
\Ee
}

\item {\bbf
a surjective map:
\[
IM_{E\to \gamma }^{(k)}: \qquad \{\delta E^{(k)}\Rl(j) \}_j\quad \To \quad
\{\delta \gamma _j^{(k)}\}_j\]
of which kernel $\Ker[IM_{E\to \gamma }^{(k)}]$ is the set
$\{\delta E^{(k)}\Rl(j)-\delta E^{(nr),(k)}\Rl(j)\}_j$
}
of difference of $k$-th consecutive spacings between the pseudoramified and pseudounramified eigenvalues of a random matrix of $\GL_r(\rit)$.

\item {\bbf
a bijective map:}
\[\mbox{\boldmath{$ IM^{(k)}_{\gamma \to E^{(nr)}_{\rm BCOE}}: \qquad
\{ \delta \gamma _j^{(k)}\}_j \quad \To \quad
\{\delta EV^{(k),(nr),\rm BCOE}\Rl(j)\}_j$}}\]
where $\delta \gamma ^{(k)}_j$ denotes a $k$-th (variable) consecutive spacing verifying
\[ \delta \gamma ^{(k)}_j=\delta E^{(k),(nr),\rm BCOE}\Rl(j)\;.\]

\Ee }
\vskip 11pt

\bpr
This proposition is a generalisation of propositions~3.29, 3.31 and 3.32 to $k$-th consecutive spacings.\epr
\vskip 11pt

\subsection{Physical interpretation of the nontrivial zeros of $\zeta (s)$}

It was suggested for a long time that the nontrivial zeros of the Riemann zeta function are probably related to the eigenvalues of some wave dynamical system of which (Hamiltonian) operator is unknown \cite{Kna}.

Considering the new mathematical framework presented here taking into account the solution of the Riemann hypothesis, the connection between these two fields is rather evident.

Indeed, referring to propositions~3.29 and 3.31, we see that there exists a surjective map:
\[
IM_{E^{(nr)}\Rl\to \gamma }: \qquad
E\Rl(j) \quad \To \quad \gamma _j\;, \qquad \forall\ j\; , \quad 1\le j\le r\;,\]
between square roots of eigenbivalues of a random matrix of 
$\GL_r(\rit\times\rit)$ (or $\GL_r(\cit\times\cit)$) and nontrivial zeros of $\zeta (s)$ in such a way that the kernel
$\Ker [IM_{E^{(nr)}\Rl\to \gamma }]$ of $IM_{E^{(nr)}\Rl\to \gamma }$ is given by the set
$\Ker [IM_{E^{(nr)}\Rl\to \gamma }]=\{E\Rl(j)-E^{(nr)}\Rl(j)\}_j$
of differences between the square roots of the pseudoramified and pseudounramified eigenbivalues of a matrix of $\GL_r(\rit\times\rit)$.

And thus, {\bbf
if $E\Rl(j)=E^{(nr)}\Rl(j)$, then we have that:}
\[\mbox{\boldmath{$ E^{(nr)}\Rl(j)=\gamma _j\;.$}}\]
The squares $E^{(nr)}\RL(j)$ of the pseudounramified eigenbivalues $E^{(nr)}\Rl(j)$ of $\GL_r(\rit\times\rit)$ are also eigenbivalues of the eigenbivalue equation:
\[
(D_R\otimes D_L)(\phi _r(G_{(Q)}(F^{(nr)}_{\o v^1}\times F^{(nr)}_{v^1})))=
E^{(nr)}\RL(j)(\phi _r(G_{(Q)}(F^{(nr)}_{\o v^1}\times F^{(nr)}_{v^1})))\]
(according to proposition~3.18),

of which eigenbifunctions is the set of $r$-bituple
\[
\langle \phi (G^{(1)}_{(Q)}(F^{(nr)}_{\o v^1}\times F^{(nr)}_{v^1})),\dots,
 \phi (G^{(j)}_{(Q)}(F^{(nr)}_{\o v^1}\times F^{(nr)}_{v^1} )) ,\dots,
 \phi (G^{(r)}_{(Q)}(F^{(nr)}_{\o v^1}\times F^{(nr)}_{v^1}))\rangle\]
 referring to proposition~3.4.
 
 And {\bbf these eigenbifunctions are the set
 $\Gamma (\wh M^{(1)}_{v_R}\otimes \wh M^{(1)}_{v_L})$ of sections of the bisemisheaf
 $(\wh M^{(1)}_{v_R}\otimes \wh M^{(1)}_{v_L})$ which was interpreted as the internal string field of an elementary (bisemi)particle according to section~3.6\/}.
 \vskip 11pt
 
 These considerations lead to the following proposition.
 \vskip 11pt

\subsection{Proposition}

{\em
The squares of the nontrivial zeros $\gamma _j$ of the Riemann zeta function
$\zeta (s)$ are the pseudounramified eigenbivalues of the eigenbivalue (biwave) equation:
\[
(D_R\otimes D_L)(\phi _r(G_{(Q)}(F^{(nr)}_{\o v^1}\times F^{(nr)}_{v^1})))=
\gamma _j^2(\phi _r(G_{(Q)}(F^{(nr)}_{\o v^1}\times F^{(nr)}_{v^1})))\]
of which eigenbifunctions are the sections of the bisemisheaf
$(\wh M^{(1)}_{v_R}\otimes \wh M^{(1)}_{v_L})$ being the internal string field of an elementary (bisemi)particle.}
\vskip 11pt

\bpr
Consequently, the squares $\gamma _j^2$ of the nontrivial zeros of $\zeta (s)$ are the eigenbivalues of a matrix of $\GL_r(\rit\times\rit)$ constituting a representation of the bilinear differential Galois semigroup associated with the action of the differential bioperator $(D_R\otimes D_L)$.

And, thus, {\bbf
each nontrivial zero $\gamma _j$ of $\zeta (s)$ is the infinitesimal generator on $j$ quanta of the Lie subsemialgebra $\gl_1(F^{(nr)}_{v_j})$ of the Lie subsemigroup\linebreak
 $\GL_1(F^{(nr)}_{v_j})$ or the energy of $j$ compact transcendental pseudounramified\linebreak quanta $(N=1)$\/}.
 
 This can also be seen from propositions~3.25 and 3.27.\epr

\vfill

C. Pierre\\
Universit\'e de Louvain\\
Chemin du Cyclotron, 2\\
B-1348 Louvain-la-Neuve,  Belgium\\
pierre.math.be@gmail.com
\eject
}

\end{document}